\documentclass[11pt,british]{article}
\usepackage{hyperref}
\usepackage{ae,aecompl}
\usepackage[authoryear,round]{natbib}
\usepackage{url}  
\usepackage{natbib}
\usepackage{setspace}
\usepackage{amsfonts}
\usepackage{amssymb}
\usepackage{float}
\usepackage{multirow}
\usepackage[labelfont=bf]{caption}
\usepackage{aecompl}
\usepackage{url}
\usepackage{adjustbox}
\usepackage{verbatim}
\usepackage{lipsum}
\usepackage{booktabs}
\usepackage{longtable}
\usepackage{amssymb}
\usepackage{tikz}
\usepackage{pgf,tikz}
\usepackage{amsmath}
\usepackage{threeparttable}
\usepackage{makecell}
\usepackage{bm}
\usepackage{geometry}
\usepackage{setspace}
\usepackage{subcaption} 
\usepackage{theorem}
\newtheorem{thm}{Theorem}
\newtheorem{proposition}{Proposition}
\newtheorem{assu}{Assumption}
\newtheorem{lem}{Lemma}
\newtheorem{defn}{Definition}
\newtheorem{exmp}{Example}

\newcommand{\appref}[1]{\hyperref[#1]{Appendix~\ref*{#1}}}
\usepackage{titletoc}

\usepackage{lscape}
\usepackage{stackengine}

\hypersetup{
colorlinks,
linktoc=all,
linktocpage=true,
citecolor=blue,
breaklinks=true,
filecolor=blue,
linkcolor=red,
urlcolor=blue,
linktocpage=true
}
\usepackage{shuffle}
\usepackage{bm}

\usepackage[affil-it]{authblk}

\usepackage{colortbl}
\usepackage{multirow}
\usepackage{xcolor}
\newcommand{\rmk}{\noindent{\textit{\textbf{Remark}}}. }
\newcommand{\prf}{\noindent{\textit{\textbf{Proof}}}. }
\newcommand{\QED}{\hfill $\square$}

\usepackage[table]{xcolor}

\begin{document}

\def\spacingset#1{\renewcommand{\baselinestretch}%
		{#1}\small\normalsize} \spacingset{1.5}


\title{\vspace{-2.5em}
\textbf{How Fast Do Signatures Learn?\\Statistical Theory and Applications for Path Regression}}

\author{\normalsize
  Blanka Horvath\textsuperscript{1,4}\thanks{Email: \texttt{blanka.horvath@maths.ox.ac.uk}},
  Wen Su\textsuperscript{1}\thanks{Email: \texttt{wen.su@maths.ox.ac.uk}},
  Wu Su\textsuperscript{2}\thanks{Email: \texttt{wusu@stu.pku.edu.cn}},
  Binnan Wang\textsuperscript{3}\thanks{Email: \texttt{wangbinnan@stu.pku.edu.cn}},
  Ruixun Zhang\textsuperscript{3,5}\thanks{Email: \texttt{zhangruixun@pku.edu.cn}}
  \\
  {\small
  \textsuperscript{1}Mathematical Institute, University of Oxford, UK\\
  \textsuperscript{2}Center for Data Science, Peking University, China\\
  \textsuperscript{3}School of Mathematical Sciences, Peking University, China\\
  \textsuperscript{4}Oxford Man Institute, University of Oxford, UK\\
  \textsuperscript{5}Center for Statistical Science, Peking University, China
  }
}

\date{\vspace{-1em} \today}

\maketitle

\vspace{-1.5em}
\noindent \textsc{\textbf{Abstract}}:
Many modern statistical learning and prediction problems involve path-valued covariates---data that evolve over time---for which path signatures have become a canonical feature representation.
Their use is justified by a universal approximation theorem, but this is an existence result: it guarantees that a finite-level signature can approximate any continuous path functional, without quantifying how fast the approximation error decreases as the truncation level grows.
This paper develops approximation and statistical theory for signature-based path regression.
We establish a minimax-optimal \(L^2\) approximation rate for smooth functionals of It\^{o} diffusions and then propagate the truncation error through three statistical learning procedures---Signature-OLS, Signature-LASSO, and Signature-Logistic---to establish their consistency.
Three real-data applications show that signatures provide informative finite-dimensional representations of path-valued covariates and can improve prediction relative to handcrafted features, in the context of finance---foreign exchange realized volatility forecasting from intraday price paths; energy---battery end-of-life prediction from early diagnostic current--voltage pulse paths; and medicine---epileptic seizure detection from short electroencephalogram windows.


\newpage

\section{Introduction}
\label{sec_introduction}

Many modern statistical learning and prediction problems involve information that evolves over time.
Foreign exchange realized volatility forecasting, battery end-of-life prediction, epileptic seizure detection, and other complex dynamic problems often include \textit{path-valued covariates}.
In these settings, the response may depend not only on the terminal state of the trajectory, but also on the path by which that state is reached, including the timing of changes, local fluctuations, cumulative variation, and ordered interactions among multiple coordinates.
We refer to regression with such path-valued covariates as \textit{path regression}.

Path-valued covariates are inherently infinite-dimensional and therefore difficult to analyze directly.
A natural approach is to transform the path into finite-dimensional features.
Handcrafted summaries can be effective in specific domains, but they are often problem-dependent and may fail to capture important pathwise interactions.
Classical functional data methods provide general basis expansions, but they do not automatically encode the noncommutative temporal order of a trajectory.
The \textit{path signature}, introduced by \citet{ChenAoM1957} and developed in rough path theory \citep{Lyons1998RMI,FirzandVictoir2010}, offers a canonical hierarchy of ordered path features through iterated integrals.
Its lower levels record path increments, while higher levels describe ordered interactions across time and coordinates.
The seminal \textit{universal approximation theorem} for signatures states that any continuous functional of a path can be approximated by a linear functional of its signature.
This universal approximation property provides the main theoretical justification for using truncated signatures in path regression \citep{KiralyOberhauser2019Kernels,Fermanian2022FunctionalSignature,chevyrev2025primer,GWZZ2025OR}.

One major theoretical gap in the signature-based statistical learning framework is that the classical universal approximation theorem alone is not sufficient to ensure consistency.
It is qualitative: it ensures approximation as the truncation order increases, but does not quantify how fast the approximation error decreases.
Once truncated signatures are used in regression, this missing rate becomes crucial.

This paper develops the approximation error rate for truncated signatures and propagates it through three widely used path regression frameworks.
First, focusing on a class of smooth functionals of It\^{o} diffusions, we prove that the minimax-optimal squared approximation error is of order \(K^{-2\gamma}\), where \(K\) is the truncation level and \(\gamma\) characterizes the smoothness of the functional.
We then use this rate to analyze ordinary least squares (OLS), least absolute shrinkage and selection operator (LASSO), and logistic regression based on truncated signatures for path-valued data, and illustrate the resulting framework in three real-data prediction problems.
In this sense, the paper provides a quantitative foundation for signature-based path regression and demonstrates its effectiveness in real-data applications.

Our contribution is threefold.
First, the established approximation rates fill the current gap in the signature-based learning framework.
We begin with Brownian motion as a benchmark to show how signatures play the role of polynomials on path space.
We then establish the main approximation result for a class of smooth coefficient functionals of It\^{o} diffusions generated through time and Stratonovich integrals.
We further derive a matching lower bound, showing that the resulting rate is optimal over the corresponding smoothness class.
Second, we develop statistical theory for regression with truncated signature features.
Using the established approximation bounds, we derive consistency results for OLS, LASSO, and logistic regression.
Our analysis also provides a general template for studying the same approximation--estimation tradeoff in other statistical learning methods based on signatures.
Third, we illustrate the framework through simulations and three prediction problems with path-valued covariates: foreign exchange realized volatility forecasting, battery end-of-life prediction, and electroencephalogram (EEG) seizure detection.
The simulations verify the theoretical results while the empirical studies demonstrate improvements over methods based on handcrafted features.

\paragraph*{Related Literature.}
This paper is related to several strands of literature.
The first strand concerns the mathematical foundations of signatures.
Signatures originate from Chen's theory of iterated integrals \citep{ChenAoM1957} and became a central object in rough path theory after \citet{Lyons1998RMI}.
Their injectivity and law-characterization properties have been studied in \citet{HamblyLyons2010Annals}, \citet{BoedihardjoGengLyonsYang2016AIM}, and \citet{ChevyrevLyons2016AOP}.
A key reason signatures are useful for statistics and machine learning is their universal approximation property: after time augmentation, linear functionals of signatures can approximate a large class of continuous path functionals.
This property turns nonlinear path dependence into linear learning on a structured hierarchy of iterated-integral features.
However, the usual universal approximation theorem is qualitative.
It guarantees the existence of accurate finite-level approximations, but does not quantify how the truncation level controls the approximation error.

The second strand concerns signature-based learning methods.
Signatures have been used as features for stream regression \citep{LevinLyonsNi2013Learning}, kernel methods \citep{KiralyOberhauser2019Kernels}, deep learning architectures \citep{KBASL2019NIPS}, neural rough differential equations \citep{Morrilletal2021ICML}, and moment-based inference for laws of stochastic processes \citep{ChevyrevOberhauser2022SignatureMoments}.
They have also been applied in a range of path-valued prediction problems, including nonparametric pricing and hedging \citep{lyons2020nonparametric}, model calibration \citep{CGS2023SIAMFIN}, path-dependent option pricing \citep{BFZ2024SIAMFM}, and optimal stopping \citep{BayerPelizzariSchoenmakers2025StoppingSignature}.
Related path-based prediction tasks also arise in atmospheric prediction \citep{FSK2024GRL}, battery prognostics \citep{Severson2019NatureEnergy,Attia2020Nature,Ibraheem2025}, and EEG seizure detection \citep{ShoebGuttag2010ICML}.
These studies show that signature features can be effective in practice.
Our paper studies the statistical question behind these applications: when a path functional is approximated by a truncated signature representation, under what conditions does the resulting estimator remain statistically consistent?

The paper is also related to the literature on functional data regression.
A common approach in this literature is to represent functional covariates through the leading eigenfunctions of their covariance operator and regress the response on the resulting functional principal component scores \citep{DauxoisPousseRomain1982,CardotFerratySarda1999,RamsaySilverman2005FDA,HallHorowitz2007FunctionalLinear}.
Unlike functional principal component regression, signatures are organized by iterated-integral order and capture nonlinear ordered interactions across time and coordinates.
\citet{Fermanian2022FunctionalSignature} studied truncated signatures in functional linear regression and showed that they are competitive with traditional functional linear methods, particularly when the functional covariates are high-dimensional.

The statistical frameworks most closely related to ours are those developed in the seminal works of \citet{Fermanian2022FunctionalSignature} and \citet{GWZZ2025OR}.
\citet{Fermanian2022FunctionalSignature} studied functional linear regression using truncated signatures, while \citet{GWZZ2025OR} established LASSO consistency for signature regression under sparse finite-level representations.
Our contribution differs from these works in three main ways.
First, we derive an \(L^2\) approximation rate for truncated signatures, whereas these studies treat the truncated signature model as exact and do not model the finite-level approximation error.
Second, \citet{Fermanian2022FunctionalSignature} studied path functionals of bounded-variation paths, and \citet{GWZZ2025OR} studied Brownian motion and discussed the Ornstein--Uhlenbeck (OU) process.
In contrast, this paper studies path functionals of general It\^{o} diffusions and identifies smoothness conditions under which quantitative signature approximation rates can be obtained.
Third, this paper develops Signature-OLS, Signature-LASSO, and Signature-Logistic frameworks.
The general strategy of balancing truncation error and statistical error can also be adapted to other statistical learning methods based on truncated signature features.

\vspace{1em}
The rest of the paper is organized as follows.
\autoref{sec:background} reviews the definitions and basic properties of signatures, the universal approximation theorem, and the signature-based path regression.
\autoref{sec:signature-approximation} develops \(L^2\) approximation rates for smooth coefficient functionals of It\^{o} diffusions.
\autoref{sec_Reg} establishes the statistical theory for Signature-OLS, Signature-LASSO, and Signature-Logistic regression.
\autoref{sec_simulation} reports Monte Carlo simulation results, and \autoref{sec:real-data} presents three real-data applications: foreign exchange realized volatility forecasting, battery end-of-life prediction, and EEG seizure detection.
\autoref{sec_conclusion} concludes.
The appendices contain technical proofs and additional numerical details.

\vspace{1em}
\noindent \textbf{Notation}:
For a vector \(\bm v\), \(\|\bm v\|\), \(\|\bm v\|_1\), and \(\|\bm v\|_\infty\) denote the Euclidean, \(\ell_1\), and maximum norms, respectively. 
For a matrix \(\bm A\), \(\bm A^\top\) denotes its transpose, \(\|\bm A\|\) denotes the operator norm, \(\|\bm A\|_\infty\) denotes the maximum absolute row-sum norm, and \(\|\bm A\|_F\) denotes the Frobenius norm.
When \(\bm A\) is symmetric, \(\lambda_{\min}(\bm A)\) and \(\lambda_{\max}(\bm A)\) denote its smallest and largest eigenvalues, respectively. 
We write \(\bm A\succeq 0\) if \(\bm A\) is positive semidefinite, and write \(\bm I_d\) for the \(d\)-dimensional identity matrix.
For a set \(A\), \(A^c\) denotes its complement and \(|A|\) denotes its cardinality.
For a path \(\bm X=\{\bm X_t\}_{t\in[0,T]}\), \(\mathcal F_T^{\bm X}\) denotes the sigma-field generated by the path up to time \(T\). 
For a domain \(D\), \(H^\gamma(D)=W^{\gamma,2}(D)\) denotes the usual \(L^2\)-Sobolev space of order \(\gamma\). 
We write \(W^{\gamma,\infty}(D)\) for the Sobolev space whose derivatives up to order \(\gamma\) are uniformly bounded.
For sequences \(a_n\) and \(b_n\), \(a_n=O(b_n)\), \(a_n=o(b_n)\), and \(a_n\asymp b_n\) have their usual asymptotic meanings. 
The notation \(O_P(\cdot)\) and \(o_P(\cdot)\) denotes stochastic order, and \(\xrightarrow{d}\) denotes convergence in distribution.

\section{Preliminaries}
\label{sec:background}

This section reviews the signature transform used throughout the paper.
Further background on signatures can be found in \cite{BayerDosReisHorvathOberhauser2026SignatureFinance}.

\subsection{Signature}
\label{subsec:signature}

Let \(\{\bm X_t\}_{t\in[0,T]}\) be an \(\mathbb R^d\)-valued continuous stochastic process, written as $\bm X_t=\left(X_t^1,X_t^2,\ldots,X_t^d\right)^\top$.
For each \(i\in\{1,2,\ldots,d\}\), the first-level signature component is defined by $S(\bm X)^i_{0,t}=\int_0^t \mathrm{d}X_s^i=X_t^i-X_0^i$, which is the increment of the \(i\)-th coordinate. 
For \(1\le i,j\le d\), we define second-level signature $S(\bm X)^{i,j}_{0,t}=\int_0^t S(\bm X)^i_{0,s}\,\mathrm{d}X_s^j=\int_{0<r<s<t} \mathrm{d}X_r^i\,\mathrm{d}X_s^j$.
More generally, for any \(1\le i_1,\ldots,i_k\le d\), we define the \(k\)-fold iterated integral
\begin{equation}
S(\bm X)^{i_1,\ldots,i_k}_{0,t}
=
\int_{0<t_1<\cdots<t_k<t}
\mathrm{d}X_{t_1}^{i_1}\cdots \mathrm{d}X_{t_k}^{i_k}.
\end{equation}
We call \(S(\bm X)^{i_1,\ldots,i_k}_{0,t}\) a \(k\)-th level signature component.
When the starting time is \(0\), we write \(S(\bm X)^{i_1,\ldots,i_k}_t\) for \(S(\bm X)^{i_1,\ldots,i_k}_{0,t}\).

\begin{defn}[Signature]
The signature of \(\bm X\) over \([0,T]\) is the collection of all iterated integrals:
\(
S(\bm X)_T
=
\left(
1,
S(\bm X)^1_T,\ldots,S(\bm X)^d_T,
S(\bm X)^{1,1}_T,\ldots,S(\bm X)^{d,d}_T,
S(\bm X)^{1,1,1}_T,\ldots
\right),
\)
where the level-0 term equals \(1\) by convention.
\end{defn}

The signature can be defined using either It\^{o} or Stratonovich iterated integrals, depending on the application.
Throughout this paper, since the Stratonovich signature satisfies the shuffle identity and is more widely used, we use it unless otherwise stated. 
Since \(S(\bm X)_T\) is infinite-dimensional, one often works with the truncated signature \(S(\bm X)^{\le K}_T\), defined to be all signature components with levels not exceeding \(K\). 
Also, the signature of time-augmented path $\widehat{\bm X}_t = (t,\bm X_t^\top)^\top$ is usually used, which removes the ambiguity caused by time reparametrization.

The signature is rich in representing path information.
For bounded-variation paths or continuous semimartingales \(\bm{X}\) and \(\bm{Y}\) starting at the same initial value $\bm{x}_0$, one has
\begin{equation}\label{eqn:time_augmented_uniqueness}
S(\widehat{\bm{X}})_T=S(\widehat{\bm{Y}})_T
\quad \Leftrightarrow\quad
\bm{X}_t=\bm{Y}_t,\,\, 0\le t\le T,\,\,\mathrm{a.s.}
\end{equation}
This property is often referred to as the uniqueness of signature; see \citet{HamblyLyons2010Annals} and Lemma 2.6 of \citet{CGS2023SIAMFIN} for related statements.
Hence, the full signature of the time-augmented path preserves the information of the whole path.

\subsection{Universal Approximation}
\label{subsec:path_regression}

We next recall the signature universal approximation theorem. 
It states that, on compact sets of time-augmented paths, every continuous path functional can be uniformly approximated by a finite linear combination of signature coordinates. 
Specifically, for any compact set \(\mathcal C\) of time-augmented paths, any continuous functional \(f:\mathcal C\to\mathbb R\), and any \(\epsilon>0\), there exists \(\bm{\ell}\) with only finitely many nonzero coordinates such that
\begin{equation}
\sup_{\widehat{\bm X}\in\mathcal C}
\left|
f(\widehat{\bm X})
-
\left\langle \bm{\ell},S(\widehat{\bm X})_T\right\rangle
\right|
\le
\epsilon .
\end{equation}
Thus, linear functionals of the full time-augmented signature form a rich class of functions on path space.
In this sense, signature coordinates play the role of monomials for path-valued inputs; see, for example, \citet{HamblyLyons2010Annals}, Proposition A.6 of \citet{KBASL2019NIPS}, and Theorem 2.12 of \citet{BFZ2024SIAMFM}.

Motivated by the universal approximation theorem, signatures are widely used in path regression.
We use signature-based linear regression as a leading example.
Let \(Z\) be a scalar response variable and let \(\mathcal F_T^{\bm X}:=\sigma(\bm X_t:0\le t\le T)\) be the sigma-field generated by the path $\bm{X}$ up to time \(T\).
We write the conditional mean of \(Z\) given the path as
$Y:=\mathbb E[Z\mid \mathcal F_T^{\bm X}]
=
f_0(\widehat{\bm X})$,
where \(f_0\) is an unknown path functional and \(f_0(\widehat{\bm X})\) is shorthand for \(f_0\bigl((\widehat{\bm X}_t)_{t\in[0,T]}\bigr)\).
Equivalently,
\begin{equation}\label{eqn:path_regression_model}
Z
=
f_0(\widehat{\bm X})+\varepsilon,
\quad
\mathbb E[\varepsilon\mid \mathcal F_T^{\bm X}]=0 .
\end{equation}
The case \(Z\in\mathcal F_T^{\bm X}\) corresponds to the path-determined setting, where \(\varepsilon=0\).

For a fixed truncation level \(K\), let
\(\bm s_K(\widehat{\bm X})\in\mathbb R^{d_K}\)
denote the vector collecting all (or a subset of all) coordinates of the truncated signature
\(S(\widehat{\bm X})^{\le K}_T\).
The best linear approximation is then defined by
\begin{equation}
\bm L_{0K}
:=
\arg\min \limits_{\bm L_K\in\mathbb R^{d_K}}
\mathbb E\left[
\left|
f_0(\widehat{\bm X})
-
\bm s_K(\widehat{\bm X})^\top \bm L_K
\right|^2
\right].
\end{equation}
The corresponding projection residual is $\xi_K^{\bm X}(f_0)
:=
f_0(\widehat{\bm X})
-
\bm s_K(\widehat{\bm X})^\top \bm L_{0K}$.
Thus, the regression model can be decomposed as
\begin{equation}
Z
=
\bm s_K(\widehat{\bm X})^\top \bm L_{0K}
+
\xi_K^{\bm X}(f_0)
+
\varepsilon.
\end{equation}
This decomposition, motivated by the universal approximation theorem, separates the population approximation problem from the finite-sample estimation problem.
However, the traditional universal approximation theorem is qualitative: it guarantees approximation but does not quantify how fast \(\xi_K^{\bm X}(f_0)\) decays as \(K\) increases.

\rmk Related approximation rates can be derived from stochastic Taylor expansions by fixing \(K\) and considering the short-time regime \(T\to0\); see \citet{CGS2023SIAMFIN}. In contrast, we fix \(T\) and quantify the approximation error as \(K\to\infty\).

\section{Signature Approximation of Path Functionals}
\label{sec:signature-approximation}

\subsection{A Motivating Example: Signature-Chaos Relation}

The universal approximation theorem explains why signatures can approximate path functionals, but it is not sufficient for statistical analysis by itself.
For statistical learning, one needs to quantify how the finite-level approximation residual \(\xi_K^{\bm X}(f_0)\) decays as the truncation level \(K\) increases.
This section studies this approximation error in \(L^2\).

We begin with Brownian motion as a benchmark case.
Since signature coordinates play the role of monomials on path space, the Brownian setting gives the cleanest way to see what finite-level signatures are doing.
In this case, the Wiener--It\^{o} chaos expansion decomposes a square-integrable path functional into deterministic kernels, while time-augmented signature coordinates generate polynomial approximations to these kernels.

Let \(B=\{B_t\}_{0\le t\le T}\) be a standard Brownian motion and let \(\widehat B_t=(t,B_t)\).
For \(k\ge1\), let \(\Delta_k(0,T):=\{0<t_1<\cdots<t_k<T\}\).
For any \(Y\in L^2(\mathcal F_T^B)\), it is well known that \(Y\) admits a Wiener--It\^{o} chaos expansion
\(Y=\sum_{k=0}^{\infty}I_k(h_k)\), where \(I_0(h_0)=\mathbb E[Y]\), and, for \(k\ge1\), \(h_k\in L^2(\Delta_k(0,T))\) and
\begin{equation}\label{eqn:brownian-multiple-ito-integral}
I_k(h_k)
:=
\int_{\Delta_k(0,T)}
h_k(t_1,\ldots,t_k)
\,\mathrm{d}B_{t_1}\cdots \mathrm{d}B_{t_k}.
\end{equation}
In \appref{app:proof-brownian-case}, we show that the level-\(K\) time-augmented Brownian signature approximates \(Y\) by approximating each low-order chaos kernel \(h_k\) with polynomials.
Thus, the approximation error is governed by the polynomial approximation error for the low-order kernels, which is controlled by the regularity of the kernels \(\{h_k\}\) through classical polynomial approximation theory.

The Brownian case is useful as a transparent benchmark, but most observed paths in applications are not Brownian motion itself.
They may exhibit drift, mean reversion, volatility structure, or other dynamics.
We therefore turn next to observable It\^{o} diffusions and define a smooth class of path functionals for which explicit signature approximation rates can be obtained.

\subsection{Smooth Functionals of It\^{o} Diffusions}

It\^{o} diffusions provide a flexible and widely used model for observed continuous-time paths in physics, economics, and finance \citep{AïtSahalia2002ECMT}.
We focus on a tractable class of path functionals generated by smooth coefficient functions evaluated along the time-augmented diffusion path.
Although this class is narrower than the class of all continuous path functionals, it is broad enough to cover many time-averaged and accumulated path-dependent quantities that arise in applications.

Let \(\bm X=\{\bm X_t\}_{0\le t\le T}\) be a \(d\)-dimensional It\^{o} diffusion satisfying
\begin{equation}
    \mathrm d\bm X_t
    =
    \bm b(t,\bm X_t)\,\mathrm dt
    +
    \bm\sigma(t,\bm X_t)\,\mathrm d\bm B_t,
    \qquad
    \bm X_0=\bm x_0,
    \label{eq:smooth-coefficient-X-SDE}
\end{equation}
where \(\bm B=\{\bm B_t\}_{0\le t\le T}\) is a \(q\)-dimensional standard Brownian motion.
Let \(D:=[0,T]\times\mathbb R^d\).
We assume that \(\bm b:D\to\mathbb R^d\) and \(\bm\sigma:D\to\mathbb R^{d\times q}\) are globally Lipschitz and uniformly bounded, i.e., $\|\bm b(t,x)\|+\|\bm\sigma(t,x)\|_F\le R$, $(t,x)\in D$, for some finite constant \(R\).

Define the time-augmented path
\(
    \widehat{\bm X}_t=(t,\bm X_t^\top)^\top .
\)
We write \(\widehat X_t^0=t\) and \(\widehat X_t^a=X_t^a\) for \(a=1,\ldots,d\).
We consider path functionals of the form
\begin{equation}
    Y
    =
    F_{\bm h}(\bm X)
    :=
    \sum_{a=0}^d
    \int_0^T
    h_a(t,\bm X_t)
    \circ\mathrm d\widehat X_t^a,
    \label{eq:smooth-coefficient-functional}
\end{equation}
where \(\bm h=(h_0,\ldots,h_d)\) is a collection of deterministic coefficient functions,
\(\circ\mathrm d\widehat X_t^0=\mathrm dt\), and
\(\circ\mathrm d\widehat X_t^a=\circ\mathrm dX_t^a\) for \(a=1,\ldots,d\).
The standard functional linear regression target \(\int_0^1 \beta(t)X_t\,\mathrm{d}t\) is also a special case of this.

This formulation is closely related to the functional It\^{o} calculus literature, which studies smooth non-anticipative functionals of semimartingale paths through horizontal and vertical derivatives \citep{ContFournie2010JFA, ContFournie2013AOP, Dupire2019FunctionalIto}.
Our class is more concrete: the functional is built from smooth coefficient functions evaluated along the observed diffusion path.
This structure is restrictive enough to allow explicit signature approximation rates, while retaining the accumulated and path-dependent quantities that are common in applications.

Let \(\gamma\in\mathbb N_+\).
We define the mixed smoothness norm
\begin{equation}
    \mathcal H_\gamma(\bm h)
    = 
    \sum_{\ell=0}^{\gamma}
    \sum_{|\alpha|\le 2\gamma+1}
    \sup_{(t,x)\in D}
    \left|
        \partial_t^\ell\partial_x^\alpha h_0(t,x)
    \right| +
    \sum_{a=1}^d
    \sum_{\ell=0}^{\gamma}
    \sum_{|\alpha|\le 2\gamma+2}
    \sup_{(t,x)\in D}
    \left|
        \partial_t^\ell\partial_x^\alpha h_a(t,x)
    \right|.
    \label{eq:Hgamma-smooth-coefficient}
\end{equation}

\begin{defn}
Define the smooth diffusion-functional class \(\mathcal U_{\gamma,R}\) as the collection of pairs \((\bm X,\bm h)\) such that: (i) \(\bm X\) solves \autoref{eq:smooth-coefficient-X-SDE} with globally Lipschitz and bounded coefficients with $\sup\limits_{(t,x)\in D}
\left\{
\|\bm b(t,x)\|+\|\bm\sigma(t,x)\|_F
\right\}
\le R$; and (ii) \(\bm h\) is smooth and satisfies \(\mathcal H_\gamma(\bm h)\le R\).
\end{defn}

\rmk
The smoothness condition is anisotropic.
The time direction requires \(\gamma\) derivatives, while the state direction requires additional smoothness.
The coefficient of the time channel \(h_0\) requires \(2\gamma+1\) spatial derivatives, whereas each spatial channel \(h_a\), \(a\ge1\), requires \(2\gamma+2\) spatial derivatives.
The additional spatial derivative for the spatial channels is used to control the Stratonovich--It\^{o} correction.
We use the same \(R\) to bound the SDE coefficients and the functional smoothness norm for simplicity.

\subsection{Signature Approximation Rate}

\begin{thm}
\label{thm:smooth-coefficient-functional}
Suppose \((\bm X,\bm h)\in \mathcal U_{\gamma,R}\), and let \(Y=F_{\bm h}(\bm X)\).
Then, for every \(K\ge2\),
\begin{equation}
    \mathbb E|\xi_K^{\bm X}(Y)|^2
    \le
    C\,\mathcal H_\gamma(\bm h)^2 K^{-2\gamma}
    =
    O(K^{-2\gamma}),
    \label{eq:smooth-coefficient-rate}
\end{equation}
where \(C>0\) does not depend on \(K\).
Moreover, $\sup \limits_{(\bm X,\bm h)\in\mathcal U_{\gamma,R}}
    \mathbb E|\xi_K^{\bm X}(F_{\bm h}(\bm X))|^2
    \asymp
    K^{-2\gamma}$.
\end{thm}

\autoref{thm:smooth-coefficient-functional} is the main approximation result of this section.
It turns the qualitative signature universal approximation theorem into a quantitative statement for a class of diffusion-driven path functionals.
Under \(\gamma\)-th order smoothness, the upper bound of the squared finite-level signature bias is of order \(K^{-2\gamma}\), and this order is sharp over \(\mathcal U_{\gamma,R}\).
Thus the truncation level \(K\) plays the same role as a polynomial degree in classical nonparametric approximation: smoother functionals admit faster approximation, while the dependence on \(K\) cannot be uniformly improved without imposing stronger regularity.

In summary, the result yields an approximation bound of the generic form
\(\mathbb E|\xi_K^{\bm X}(Y)|^2=O(K^{-Q})\)
for a positive exponent \(Q\) determined by the regularity of the target path functional.
For the smooth diffusion-functional class above, \(Q=2\gamma\).
This polynomial-rate bound is the key population input for the statistical analysis in the next section.
When stronger regularity is available, \appref{app:analytic-rate} shows that analytic coefficient functions lead to exponential signature approximation rate $O(\rho^{-2K})$ for some $\rho$.

\section{Statistical Theory for Path Regression via Signature}
\label{sec_Reg}

This section develops limit theory for three finite-level signature learning procedures: Signature-OLS, Signature-LASSO, and Signature-Logistic.
The approximation results in \autoref{sec:signature-approximation} provide the population input for this analysis.
Our goal is to show how the finite-level approximation residual affects statistical estimation.
Throughout this section, we continue to work with the case where
\(f_0(\widehat{\bm X})=F_{\bm h}(\bm X)\) for
\((\bm X,\bm h)\in\mathcal U_{\gamma,R}\).

In all three frameworks, the key issue is the approximation-estimation tradeoff induced by the truncation level \(K\), which is summarized in \autoref{tab:tradeoff}.

\begin{table}[h!]
\centering
\caption{Approximation-estimation tradeoff induced by the truncation level \(K\).}\label{tab:tradeoff}
\small
\setlength{\tabcolsep}{15pt}
\renewcommand{\arraystretch}{0.95}
\begin{tabular}{ccc}
\toprule
Truncation level & Estimation side & Approximation side \\
\midrule
\(K \uparrow\) 
& \(d_K\uparrow\), statistical error \(\uparrow\) 
& \(K^{-Q}\downarrow\), approximation error \(\downarrow\) \\

\(K \downarrow\) 
& \(d_K\downarrow\), statistical error \(\downarrow\) 
& \(K^{-Q}\uparrow\), approximation error \(\uparrow\) \\
\bottomrule
\end{tabular}
\end{table}

\subsection{Signature-OLS: Path Regression Consistency}

We first recall the signature OLS procedure.
Suppose \(Y:=\mathbb E[Z\mid\mathcal F_T^{\bm X}]=f_0(\widehat{\bm X})\), where \(\widehat{\bm X}:=(\widehat{\bm X}_t)_{t\in[0,T]}\) is the time-augmented path, and write \(\varepsilon:=Z-f_0(\widehat{\bm X})\). 
For Signature-OLS, let \(\bm s_K(\widehat{\bm X})\in\mathbb R^{d_K}\) be a   subvector of \(S(\widehat{\bm X})_{0,T}^{\le K}\), obtained by deleting exact linear redundancies.
Define the best linear coefficient and its approximation residual by
\begin{equation}
\label{eq:population-best-linear-coefficient}
\begin{aligned}
&\bm L_{0K}
:=
\arg\min_{\bm L_K\in\mathbb R^{d_K}}
\mathbb E\left[
\left|
f_0(\widehat{\bm X})
-
\bm s_K(\widehat{\bm X})^\top\bm L_K
\right|^2
\right],\\
&\xi_K^{\bm X}(f_0)
:=
f_0(\widehat{\bm X})
-
\bm s_K(\widehat{\bm X})^\top\bm L_{0K}.
\end{aligned}
\end{equation}
Then \(Z=\bm s_K(\widehat{\bm X})^\top\bm L_{0K}+\xi_K^{\bm X}(f_0)+\varepsilon\). By construction, \(\bm s_K(\widehat{\bm X})^\top\bm L_{0K}\) is the \(L^2\)-projection of \(f_0(\widehat{\bm X})\) onto the level-\(K\) signature approximation space, and \(\mathbb E[\bm s_K(\widehat{\bm X})\xi_K^{\bm X}(f_0)]=\bm 0\).

Suppose that \(\{(Z_i,\bm X_i)\}_{i=1}^n\) are independent and identically distributed,\footnote{This might be extended to blocks of strong mixing processes; see, for example, \citet{LucchesePakkanenVeraart2025ICML}.} and let \(\mathcal X_n:=\sigma(\widehat{\bm X}_1,\ldots,\widehat{\bm X}_n)\).
Collect the responses, features, approximation residuals, and regression errors as
\begin{equation}
\label{eq:sample-path-regression-model}
\begin{aligned}
&\mathbf Z
:=(Z_1,\ldots,Z_n)^\top,
&
&\mathbf S
:=
(\bm s_K(\widehat{\bm X}_1),\ldots,
\bm s_K(\widehat{\bm X}_n))^\top,\\
&\bm\xi_K
:=
(\xi_K^{\bm X_1}(f_0),\ldots,
\xi_K^{\bm X_n}(f_0))^\top,
& &\bm\varepsilon
:=(\varepsilon_1,\ldots,\varepsilon_n)^\top,\\
&\mathbf Z
=
\mathbf S\bm L_{0K}
+
\bm\xi_K
+
\bm\varepsilon.
\end{aligned}
\end{equation}
Whenever \(\mathbf S^\top\mathbf S\) is invertible, the Signature-OLS estimator is \(\widehat{\bm L}_K=(\mathbf S^\top\mathbf S)^{-1}\mathbf S^\top\mathbf Z\), and its prediction at a fixed path \(\widehat{\bm x}:=(\widehat{\bm x}_t)_{t\in[0,T]}\) is \(\widehat f_K(\widehat{\bm x})=\bm s_K(\widehat{\bm x})^\top\widehat{\bm L}_K\).

\begin{assu}
\label{assu: PathReg}
(i) 
We assume \(\bm\Sigma_K
:=
\mathbb E[
\bm s_K(\widehat{\bm X})
\bm s_K(\widehat{\bm X})^\top]\) is positive definite for every \(K\), and $\frac{1}{n}
\mathbb E\left[
\left\{
\bm s_K(\widehat{\bm X})^\top
\bm\Sigma_K^{-1}
\bm s_K(\widehat{\bm X})
\right\}^2
\right]
\to 0$.
(ii) There exists \(Q>0\) such that
\(\mathbb E[\{\xi_K^{\bm X}(f_0)\}^2]=O(K^{-Q})\).
(iii) Conditional on \(\mathcal X_n\), the errors are independent, and almost surely,
$\mathbb E[\varepsilon_i\mid\mathcal X_n]=0$,
$\mathbb E[\varepsilon_i^2\mid\mathcal X_n]=\sigma^2>0$,
and
$\max \limits_{1\le i\le n}
\mathbb E[|\varepsilon_i|^4\mid\mathcal X_n]
\le C_{\varepsilon}<\infty$.
\end{assu}

\rmk
Condition (i) ensures that the least-squares coefficient is uniquely defined for every finite \(K\).
By contrast, the subsequent LASSO theory does not require the full Gram matrix to be nonsingular.
The condition $\frac{1}{n}
\mathbb E\left[
\left\{
\bm s_K(\widehat{\bm X})^\top
\bm\Sigma_K^{-1}
\bm s_K(\widehat{\bm X})
\right\}^2
\right] \to 0$ requires the sample size to be large enough to assure consistency.
This rate condition also implies \(d_K^2/n\to0\).
The exponent \(Q\) in condition (ii) reflects the smoothness of the target functional, while condition (iii) is standard in regression analysis.

For pointwise prediction, let
\(\widehat{\bm x}:=(\widehat{\bm x}_t)_{t\in[0,T]}\)
be a deterministic time-augmented test path in the domain of \(f_0\) that has a well-defined signature.
We assume $c\le \bm s_K(\widehat{\bm x})^\top
\bm\Sigma_K^{-1}
\bm s_K(\widehat{\bm x}) \le  Cd_K$ for some $c,C>0$, and \(|\xi_K^{\bm x}(f_0)|^2=O(K^{-Q})\).
This condition is imposed because we restrict attention to the probability space of the sample and do not wish to treat the test path as random.
If the test path is drawn from the same distribution as the sample paths, these rate conditions hold with probability approaching one.

\begin{thm}
\label{thm: OLSCLT}
Let \(K=K(n)\to\infty\), and let \(\widehat{\bm x}\) be a test path satisfying the preceding conditions.
Under \autoref{assu: PathReg}, the following statements hold.\\
(i) If \(d_K/K^Q\to0\), then
$
\left|
\widehat f_K(\widehat{\bm x})
-
f_0(\widehat{\bm x})
\right|^2
=
O_P\left(
\frac{d_K}{n}
+
\frac{d_K}{K^Q}
\right)
=
o_P(1).
$\\
(ii) If \(n/K^Q\to0\), then
$
\frac{
\sqrt n
\{\widehat f_K(\widehat{\bm x})-f_0(\widehat{\bm x})\}
}{
\{\bm s_K(\widehat{\bm x})^\top
\bm\Sigma_K^{-1}
\bm s_K(\widehat{\bm x})\}^{1/2}
}
\xrightarrow{d}
\mathcal N(0,\sigma^2).
$
\end{thm}

\rmk 
Generally, \autoref{assu: PathReg} needs $d_K^2/n \to 0$, which is relatively mild: it requires the sample size \(n\) to be large enough relative to the number of parameters being estimated.
This condition ensures that the sample Gram matrix \(\mathbf S^\top\mathbf S/n\) concentrates around its population counterpart \(\bm\Sigma_K\) and remains invertible with high probability.
The condition \(d_K/K^Q\to0\) makes the approximation contribution vanish, while the stronger condition \(n/K^Q\to0\) in part (ii) makes it negligible on the standard-error scale.

We illustrate why \(\bm s_K(\widehat{\bm X})\) is chosen to be a retained signature vector.
Specifically, if all time-augmented signature coordinates up to level \(K\) are included, then the dictionary dimension grows exponentially, \(d_K\asymp(d+1)^K\).
Under the polynomial approximation rate \(O(K^{-Q})\), one has \(d_K/K^Q\to\infty\) for every finite \(Q\), so the consistency condition \(d_K/K^Q\to0\) cannot hold. 
There are two natural ways to address this dimensional limitation:
\\
\textit{1. Analytic approximation.}
For analytic coefficient functionals, \appref{app:analytic-rate} gives an approximation error of \(O(\rho^{-2K})\), where \(\rho\) is related to the analytic radius, and in the case \(\rho>d+1\), $\left|
\widehat f_K(\widehat{\bm x})
-
f_0(\widehat{\bm x})
\right|^2
=
O_P\left(
\frac{(d+1)^K}{n}
+
\frac{(d+1)^K}{\rho^{2K}}
\right)$.
Balancing the two terms gives the rate-optimal choice $K\asymp\frac{\log n}{2\log\rho}$.
Since \(\rho\) is unknown, one may search \(K<\frac{\log n}{2\log(d+1)}\).
\\
\textit{2. Signature selection.}
A more practical approach is to apply OLS to a post-selection subset of signature coordinates.
The Signature-LASSO procedure introduced below provides such a feature-selection mechanism.
This approach can reduce the retained dimension to \(d_K\ll n\), although an explicit optimal choice of \(K(n)\) is more difficult to derive because the selected feature dimension is data-dependent.

\subsection{Signature-LASSO: Sparse Recovery}

In many applications, a path functional may depend only on a small subset of signature coordinates.
This motivates Signature-LASSO, which selects relevant terms from a large  signature dictionary.
The following OU example illustrates how such sparsity may arise.

\begin{exmp}
\label{exmp: OU}
Suppose \(Y\) satisfies the OU SDE \(\mathrm dY_t=-\theta Y_t\,\mathrm dt+\sigma\,\mathrm dB_t\).
Then \(Y_t\) admits a sparse representation in terms of the signature of the time-augmented Brownian path.
Let \(\bm X_t=(X_t^0,X_t^1)=(t,B_t)\).
A third-order truncation yields
\begin{equation}
\begin{aligned}
Y_t
&=Y_0e^{-\theta t}+\int_0^t\sigma e^{-\theta(t-s)}\,\mathrm dB_s\\
&\approx
Y_0e^{-\theta t}
+\sigma e^{-\theta t}
\int_0^t
\left(1+\theta s+\frac{\theta^2s^2}{2}\right)\mathrm dB_s\\
&=
Y_0e^{-\theta t}S(\bm X)_t^{\emptyset}
+\sigma e^{-\theta t}S(\bm X)_t^1
+\theta\sigma e^{-\theta t}S(\bm X)_t^{0,1}
+\theta^2\sigma e^{-\theta t}S(\bm X)_t^{0,0,1}.
\end{aligned}
\end{equation}
Thus, in this third-order approximation, only the signature coordinates indexed by
\(\emptyset\), \(1\), \((0,1)\), and \((0,0,1)\) enter the representation.
\end{exmp}

The example suggests that the population representation may be sparse.
Throughout this subsection, \(\bm s_K\) can be the full truncated signature dictionary and need not coincide with the retained vector used for Signature-OLS.
But for each \(K\), we suppose the population coefficient vector \(\bm L_{0K}\in\mathbb R^{d_K}\) is sparse.
Specifically, let \(A_0:=\operatorname{supp}(\bm L_{0K})\) and \(q_K:=|A_0|\), where \(1\le q_K\ll d_K\).
The Signature-LASSO estimator is defined as
\begin{equation}
\label{eq:lasso}
\widehat{\bm L}_K
\in
\arg\min_{\bm L\in\mathbb R^{d_K}}
\left\{
\frac{1}{2n}\|\mathbf Z-\mathbf S\bm L\|^2
+\lambda_n\sum_{j=1}^{d_K}w_{Kj}|L_j|
\right\},
\end{equation}
where, if the \(j\)-th coordinate is indexed by the word \(I_j\), we set \(k_j:=|I_j|\) and
\(w_{Kj}:=T^{k_j/2}/\sqrt{k_j!}\).
By \autoref{lem:sigmoment}, \(w_{Kj}\) is the deterministic time-and-factorial factor in the uniform moment bound, with \(1/\sqrt{k_j!}\) reflecting factorial decay across signature levels.
Let \(\overline{\bm s}_K\) have coordinates \(\overline s_{Kj}:=s_{Kj}/w_{Kj}\).
Since \(s_{Kj}=w_{Kj}\overline s_{Kj}\), the original and rescaled dictionaries span the same linear space.

Define the (weighted) population Gram matrix and its active--inactive partition by
\begin{equation}
\bm\Gamma_K
:=
\mathbb E\left[
\overline{\bm s}_K(\widehat{\bm X})
\overline{\bm s}_K(\widehat{\bm X})^\top
\right]
=
\begin{pmatrix}
(\bm\Gamma_K)_{A_0A_0} & (\bm\Gamma_K)_{A_0A_0^c}\\
(\bm\Gamma_K)_{A_0^cA_0} & (\bm\Gamma_K)_{A_0^cA_0^c}
\end{pmatrix}.
\end{equation}

\begin{assu}
\label{assu:lasso}
(i) The active block \((\bm\Gamma_K)_{A_0A_0}\) is positive definite for every \(K\), and $\sup\limits_K\left\|(\bm\Gamma_K)_{A_0A_0}^{-1}\right\|_\infty<\infty$, $\sup \limits_K \left\|(\bm\Gamma_K)_{A_0^cA_0}(\bm\Gamma_K)_{A_0A_0}^{-1}\right\|_\infty<1$.
(ii) For some \(Q>0\), \(\mathbb E[\{\xi_K^{\bm X}(f_0)\}^2]=O(K^{-Q})\).
(iii) Conditional on \(\mathcal X_n\), the errors \(\{\varepsilon_i\}_{i=1}^n\) are independent and satisfy, almost surely, \(\mathbb E[\varepsilon_i\mid\mathcal X_n]=0\) and \(\mathbb E[\exp(t\varepsilon_i)\mid\mathcal X_n]\le\exp(\sigma_0^2t^2/2)\) for every \(t\in\mathbb R\) and \(i=1,\ldots,n\), where \(\sigma_0>0\) is independent of \(n\) and \(K\).
\end{assu}

\rmk
Condition (i) combines a uniform inverse bound for the weighted active Gram matrix with the classical irrepresentable condition, a standard sufficient condition for support recovery \citep{ZhaoYu2006JMLR,Wainwright2009IT}.
It does not require the full Gram matrix to be nonsingular.
Condition (ii) connects the sparse population representation to the approximation theory in \autoref{sec:signature-approximation}, while condition (iii) provides the conditional sub-Gaussian error bound used in the high-dimensional analysis.

For pointwise prediction, let \(\widehat{\bm x}\) be a deterministic time-augmented test path and we assume that
$\bm s_{K,A_0}(\widehat{\bm x})^\top
(\bm\Sigma_K)_{A_0A_0}^{-1}
\bm s_{K,A_0}(\widehat{\bm x})
=O(1)$ and
$|\xi_K^{\bm x}(f_0)|^2=O(K^{-Q})$.

\begin{thm}
\label{thm:lasso_finite_sample}
Under the preceding i.i.d. sampling setup, let \(K=K(n)\to\infty\) and \(d_K\to\infty\), and suppose that \autoref{assu:lasso} holds.
For a sufficiently large fixed constant \(C>0\) depending only on \(R\) and \(T\), suppose that
\begin{equation}
\begin{aligned}
&q_K(CK)^K
\left\{
\sqrt{\frac{\log(e d_K)}{n}}
+\frac{\{\log(e d_K)\}^K}{n}
\right\}
\to0,\\
&\lambda_n
\asymp
\sqrt{q_K}(CK)^K
\left\{
\sqrt{\frac{\log(e d_K)}{n}}
+\frac{\{\log(e d_K)\}^K}{n}
\right\},
\end{aligned}
\end{equation}
where \(\lambda_n\) is chosen at the displayed order with a sufficiently large fixed leading constant.
Then the solution of \autoref{eq:lasso} is unique with probability approaching one, and the following statements hold.\\
(i)
\(\mathbb P\{\operatorname{supp}(\widehat{\bm L}_K)\subseteq A_0\}\to1\).\\
(ii) $
\max\limits_{1\le j\le d_K}
w_{Kj}|\widehat L_{K,j}-L_{0K,j}|
=O_P(\lambda_n)$, and
$\left\{
\sum_{j=1}^{d_K}
w_{Kj}^2(\widehat L_{K,j}-L_{0K,j})^2
\right\}^{1/2}
=O_P(\sqrt{q_K}\lambda_n)$.
\\
(iii)
For test path $\widehat{\bm x}$,
$\left|
\widehat f_K(\widehat{\bm x})
-f_0(\widehat{\bm x})
\right|^2
=
O_P\left(q_K\lambda_n^2+K^{-Q}\right)
=o_P(1)$.
\\
(iv)
If
\(\min \limits_{j\in A_0}w_{Kj}|L_{0K,j}|>C_\beta\lambda_n\) for sufficiently large \(C_\beta\), then
\(\mathbb P\{\operatorname{sign}(\widehat{\bm L}_K)=\operatorname{sign}(\bm L_{0K})\}\to1\).
\end{thm}

\rmk
In contrast to the OLS condition \(d_K/\sqrt n\to0\), the LASSO rate condition is driven by the active-set size \(q_K\), with the ambient dimension entering only through \(\log(e d_K)\).
In particular,
\(q_K/\sqrt n
\le
q_K(CK)^{K}\sqrt{\log(e d_K)/n}
\to0\).
At the stated penalty level, \(\lambda_n=o(1)\), \(\sqrt{q_K}\lambda_n=o(1)\).
The beta-min condition is imposed on the coefficients in the rescaled coordinates.
Notably, the exponent \(Q\) does not enter the weighted estimation rates for \(\bm L_{0K}\), so unlike the OLS case, a polynomial approximation rate does not by itself preclude LASSO selection consistency even when we include all truncated signature coordinates.

\subsection{Signature-Logistic: Classification Consistency}

We now consider a binary response \(Z\in\{0,1\}\) and suppose that
\begin{equation}
\label{eq:classification-logit-model}
\mathbb P(Z=1\mid\widehat{\bm X})
=
\Lambda\{f_0(\widehat{\bm X})\},
\qquad
\Lambda(u):=\frac{1}{1+e^{-u}}.
\end{equation}
Here \(f_0(\widehat{\bm X})\) is the conditional log-odds. 
We use the same finite-dimensional setting and notation as in the Signature-OLS analysis. 
Conditional on \(\mathcal X_n\), the responses are independent and satisfy
\(\mathbb E[Z_i\mid\mathcal X_n]
=\Lambda\{f_0(\widehat{\bm X}_i)\}\).
Whenever a finite minimizer exists, define the Signature-Logistic estimator by
\begin{equation}
\label{eq:logistic-signature-estimator}
\begin{aligned}
&\widehat{\bm L}_K
=
\arg\min_{\bm L\in\mathbb R^{d_K}}
\mathcal L_n(\bm L),\\
&\mathcal L_n(\bm L)
:=
\frac{1}{n}
\sum_{i=1}^n
\left[
\log\left\{1+\exp\bigl(
\bm s_K(\widehat{\bm X}_i)^\top\bm L
\bigr)\right\}
-
Z_i\bm s_K(\widehat{\bm X}_i)^\top\bm L
\right].
\end{aligned}
\end{equation}

Recall that
\(\bm\Sigma_K
=\mathbb E[\bm s_K(\widehat{\bm X})
\bm s_K(\widehat{\bm X})^\top]\), and define the population logistic
Hessian by
\begin{equation}
\bm H_K(\bm L)
:=
\mathbb E
\left[
\Lambda\{\bm s_K(\widehat{\bm X})^\top\bm L\}
\left[
1-\Lambda\{\bm s_K(\widehat{\bm X})^\top\bm L\}
\right]
\bm s_K(\widehat{\bm X})
\bm s_K(\widehat{\bm X})^\top
\right].
\end{equation}

\begin{assu}
\label{assu:classification}
(i) \(\bm\Sigma_K\) is positive definite for every \(K\), and
$\frac{1}{n}
\mathbb E\left[
\left\{
\bm s_K(\widehat{\bm X})^\top
\bm\Sigma_K^{-1}
\bm s_K(\widehat{\bm X})
\right\}^2
\right]
\to0$.
(ii) The population Hessian satisfies
$\inf\limits_{
\left\|
\bm\Sigma_K^{1/2}(\bm L-\bm L_{0K})
\right\|\le r
}
\lambda_{\min}
\left\{
\bm\Sigma_K^{-1/2}
\bm H_K(\bm L)
\bm\Sigma_K^{-1/2}
\right\}
\ge c$ for some $r,c>0$.
(iii) The population approximation error satisfies
\(\mathbb E[\{\xi_K^{\bm X}(f_0)\}^2]=O(K^{-Q})\).
\end{assu}

\rmk
Condition (i) is the relative Gram condition used in the Signature-OLS analysis and implies \(d_K^2/n\to0\).
Condition (ii) imposes local population curvature in the same Gram geometry and is therefore invariant to the scale of the signature coordinates. 
As before, \(Q\) in condition (iii) reflects the smoothness of \(f_0\).

For pointwise prediction, we use the same deterministic test-path
setting as in the Signature-OLS analysis. For a test path
\(\widehat{\bm x}\) satisfying those conditions, define
\(\widehat f_K(\widehat{\bm x})
:=\bm s_K(\widehat{\bm x})^\top\widehat{\bm L}_K\).

\begin{thm}
\label{thm:classification-consistency}
Let \(K=K(n)\to\infty\), and let \(\widehat{\bm x}\) satisfy the
preceding test-path conditions. Suppose that
\autoref{assu:classification} holds and \(d_K/K^Q\to0\). Then
\(\mathcal L_n\) admits a finite and unique minimizer with probability
approaching one, and
\begin{equation}
\begin{aligned}
&\left\|
\bm\Sigma_K^{1/2}
(\widehat{\bm L}_K-\bm L_{0K})
\right\|
=
O_P\left(
\sqrt{\frac{d_K}{n}}
+
K^{-Q/2}
\right),\\
&\left|
\widehat f_K(\widehat{\bm x})
-
f_0(\widehat{\bm x})
\right|
=
O_P\left(
\frac{d_K}{\sqrt n}
+
\frac{\sqrt{d_K}}{K^{Q/2}}
\right)
=o_P(1).
\end{aligned}
\end{equation}
\end{thm}

\rmk
Since \(\Lambda\) is Lipschitz, the same rate holds for the conditional success probability.

\vspace{1em}

Overall, the three results in this section show how the approximation rate \(\mathbb E|\xi_K^{\bm X}(f_0)|^2=O(K^{-Q})\) connects signature approximation with statistical learning.
The same idea can also be extended to other signature-based learning procedures, such as ridge regression and elastic-net regularization.

\section{Simulation Study}
\label{sec_simulation}

This section verifies the theory in controlled settings where the true target
functional is known.
Across three designs, we examine population approximation, the approximation--estimation tradeoff that governs the choice of truncation level, and sparse recovery under Signature-LASSO.
All implementation details, including the data-generating processes, tuning choices, and evaluation metrics, are collected in \appref{app:simulation-details}.

\begin{figure}[h!]
    \centering
    \includegraphics[width=0.96\linewidth]
    {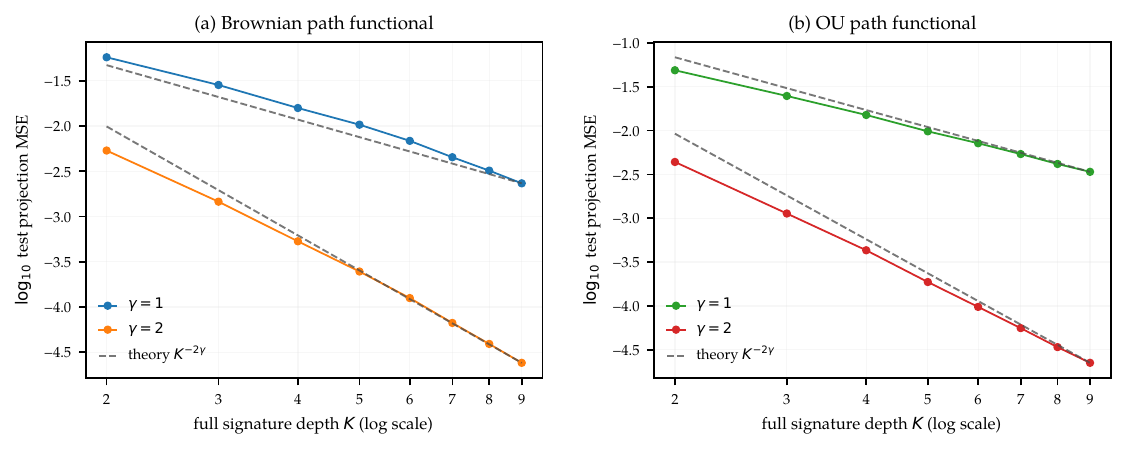}
    \caption{Population signature approximation under Brownian motion (a) and OU diffusion (b). Dashed lines show \(K^{-2\gamma}\) reference rates.}
    \label{fig:sim-summary-approx}
\end{figure}

\paragraph*{Population approximation.}
The first experiment constructs coefficient functionals with prescribed smoothness orders \(\gamma\in\{1,2\}\) under Brownian motion and an OU diffusion.
\autoref{fig:sim-summary-approx} shows that the projection error decreases as the signature depth \(K\) increases.
The smoother target, indexed by \(\gamma=2\), is approximated faster than the rougher target, indexed by \(\gamma=1\), in both benchmarks.
Overall, as the truncation order $K$ increases, the approximation error converges to 0 quickly.

\begin{figure}[h!]
    \centering
    \includegraphics[width=0.96\linewidth]
    {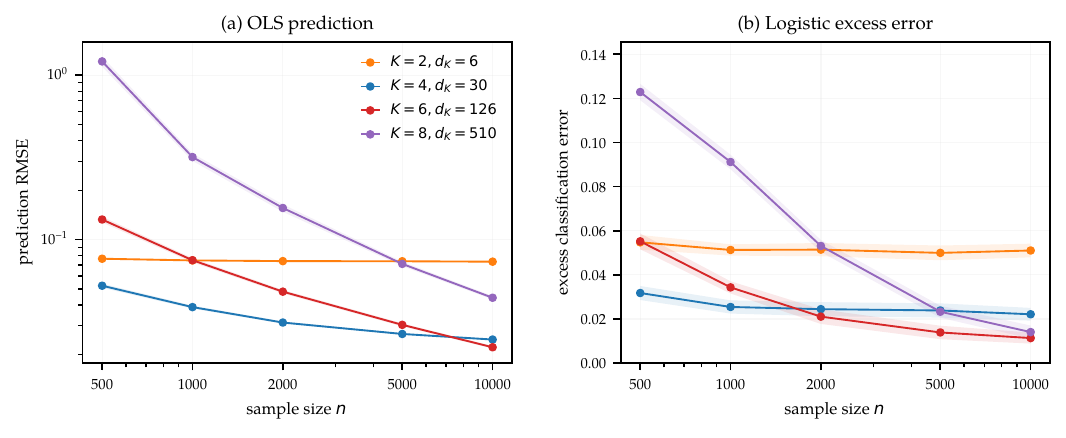}
    \caption{Finite-sample signature regression: OLS (a) and logistic regression (b). Shaded bands are 95\% Monte Carlo intervals.}
    \label{fig:sim-summary-regression}
\end{figure}

\paragraph*{OLS and logistic regression.}
The second experiment uses the smooth Brownian functional from the first
design as the latent signal.
We generate a continuous response by adding Gaussian noise to the functional
value and a binary response by applying a logistic link to the same latent
signal.
We then fit Signature-OLS and Signature-Logistic across a range of truncation
depths.
\autoref{fig:sim-summary-regression} illustrates the finite-sample
approximation--estimation tradeoff.
Small truncation levels are stable but underfit, whereas deeper signatures
reduce approximation bias at the cost of estimating more coefficients.
In the one-dimensional time-augmented case, the full nonconstant dictionary
has \(d_K=2^{K+1}-2\) coordinates, so \(K=8\) already corresponds to \(510\)
regressors.
In both experiments, an intermediate depth gives the best finite-sample
balance.

\begin{figure}[h!]
    \centering
    \includegraphics[width=0.96\linewidth]{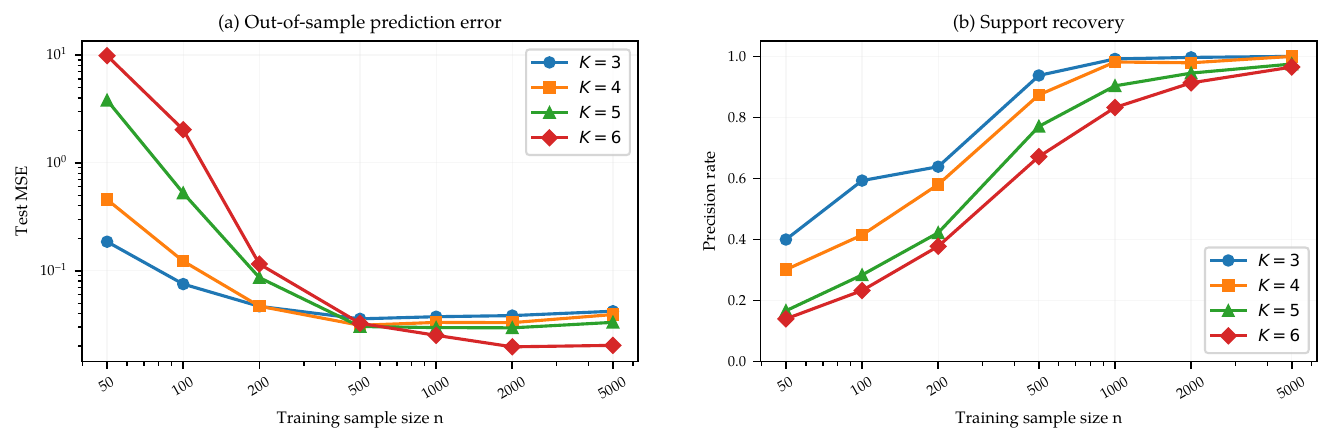}
    \caption{Signature-LASSO under the OU target: prediction error (a) and selection precision rate (b). Results are averaged over 100 independent replications.}
    \label{fig:sim-summary-lasso}
\end{figure}

\paragraph*{Signature-LASSO.}
The third experiment is motivated by the sparse OU structure in
\autoref{exmp: OU}.
We generate responses from an OU-driven integral target whose finite-level
signature expansion is concentrated on the OU-relevant family
\((1),(0,1),(0,0,1),\ldots\).
This design creates a setting in which the full signature dictionary is large, but the relevant signal is carried by a small structured subset of coordinates.
\autoref{fig:sim-summary-lasso} reports prediction and selection performance.
The test error decreases with the sample size, and deeper signatures become
useful once there is enough data to estimate the enlarged dictionary.
Precision for the OU-relevant family also improves with the sample size, but
selection becomes harder at larger \(K\) because the dictionary expands
rapidly and higher-order OU coordinates carry weaker signal.
Thus the appropriate interpretation is recovery of a finite-level OU-relevant family, rather than exact support recovery of an infinite integral functional.

\section{Real Data Applications}
\label{sec:real-data}

This section illustrates the finite-level signature framework in three applications that match the statistical regimes studied above: path regression, high-dimensional sparse regression, and binary classification. We apply the framework to forecasting foreign exchange realized volatility from intraday price paths, predicting battery end of life from short diagnostic pulse trajectories, and detecting epileptic seizures from multichannel EEG windows. This section presents the empirical design and the evidence needed to interpret the main findings; full data construction, implementation details, and sensitivity analysis are reported in \appref{app:additional-experiments}.





\subsection{Foreign Exchange Realized Volatility Forecasting}
\label{subsec:fx_rv_application}

We first apply signature regression to forecast next-day realized volatility in foreign exchange markets. We use one-minute mid-price series for 12 major U.S. dollar currency pairs from the HistData database over 2013--2023. If \(r_{t,j}\) is the \(j\)-th intraday log return on day \(t\), then \(\mathrm{RV}_t=(\sum_j r_{t,j}^2)^{1/2}\), and the response is \(Z_{t+1}=\log \mathrm{RV}_{t+1}\). Each currency-pair sample was split chronologically, with the first 80\% for estimation and the last 20\% for out-of-sample evaluation.

For each day, we form a cumulative-return path and apply a time-augmented lead--lag transformation, which preserves temporal order while encoding quadratic-variation information. We compute its depth-3 signature and concatenate the signatures from the most recent 22 trading days, producing 858 candidate features. Estimation then proceeded in two steps: a five-fold cross-validated Signature-LASSO selected coordinates on the training sample, and Signature-OLS was refitted on the selected coordinates. This specification was compared with a naive forecast, AR(1), and four HAR-family models based on realized volatility, semivariance, and jump summaries. The complete protocol and pair-level feature-selection results are reported in \appref{Appen_exchangerates}.

\begin{figure}[h!]
    \centering
    \includegraphics[width=\linewidth]
    {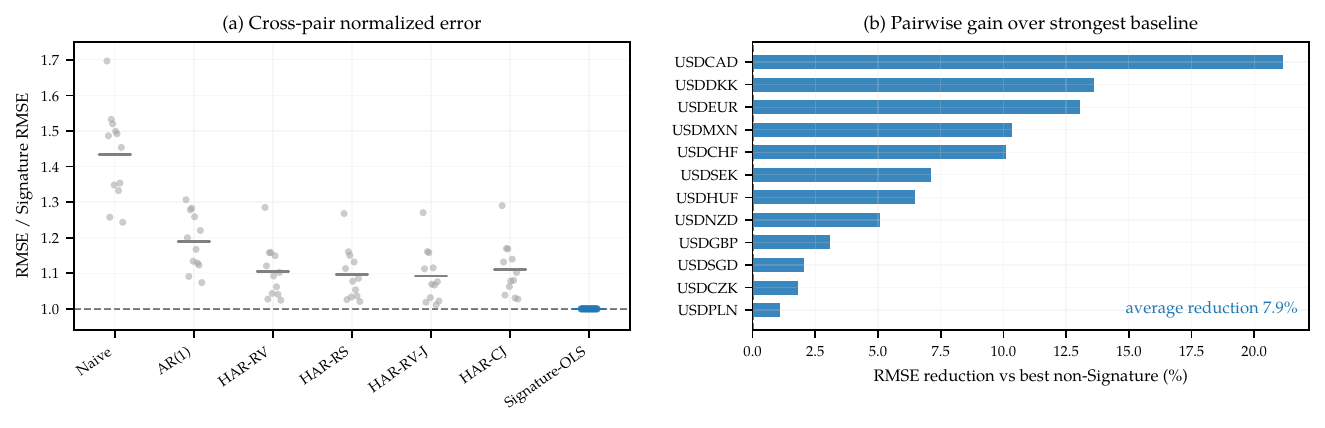}
    \caption{Relative root mean squared error (RMSE) performance across currency pairs. Panel (a) reports
    each benchmark's RMSE relative to Signature-OLS; the dashed line marks parity.
    Panel (b) reports the percentage RMSE reduction of Signature-OLS relative to
    the best non-signature baseline for each pair.}
    \label{fig:USDEUR_forecast_errors_main}
\end{figure}

\autoref{fig:USDEUR_forecast_errors_main} shows that Signature-OLS has lower RMSE than every competing benchmark for all 12 pairs and reduces RMSE by 7.9\%, on average, relative to the best non-signature model. The improvement is therefore not driven by one currency or one unusually weak comparator.

\begin{table}[h!]
\centering
\caption{Cross-pair summary of out-of-sample FX forecasting performance. RMSE and
MAE denote root mean squared error and mean absolute error for log realized volatility;
\(R^2_{\mathrm{oos}}\) denotes the out-of-sample \(R^2\); and QLIKE denotes
quasi-likelihood loss on the variance scale. Entries are averages across the 12 currency
pairs, and the last four columns count pairwise wins; arrows indicate the preferred direction.}
\label{tab:fx_performance_summary}
\small
\setlength{\tabcolsep}{5.2pt}
\renewcommand{\arraystretch}{0.85}
\resizebox{\linewidth}{!}{
\begin{tabular}{lrrrrrrrr}
\toprule
& \multicolumn{4}{c}{\textbf{Cross-pair average}}
& \multicolumn{4}{c}{\textbf{Number of pairwise wins}} \\
\cmidrule(lr){2-5}\cmidrule(lr){6-9}
\textbf{Method}
& \textbf{RMSE $\downarrow$} & \textbf{MAE $\downarrow$} & \(\bm{R^2_{\mathrm{oos}}}\) \(\bm{\uparrow}\) & \textbf{QLIKE $\downarrow$}
& \textbf{RMSE} & \textbf{MAE} & \(\bm{R^2_{\mathrm{oos}}}\) & \textbf{QLIKE} \\
\midrule
Signature-OLS & \textbf{0.749} & \textbf{0.565} & \textbf{0.361} & \textbf{-9.375} & \textbf{12} & \textbf{11} & \textbf{12} & \textbf{7} \\
HAR-RS        & 0.820 & 0.612 & 0.237 & -9.360 & 0 & 0 & 0 & 3 \\
HAR-RV-J      & 0.817 & 0.617 & 0.243 & -9.352 & 0 & 0 & 0 & 0 \\
HAR-RV        & 0.827 & 0.611 & 0.224 & -9.362 & 0 & 1 & 0 & 2 \\
HAR-CJ        & 0.830 & 0.615 & 0.218 & -9.344 & 0 & 0 & 0 & 0 \\
AR(1)         & 0.891 & 0.685 & 0.106 & -9.260 & 0 & 0 & 0 & 0 \\
Naive         & 1.076 & 0.792 & -0.306 & -8.891 & 0 & 0 & 0 & 0 \\
\bottomrule
\end{tabular}
}
\end{table}

The selected features, reported in \appref{Appen_exchangerates}, are dominated
by higher-order lead--lag terms and display a multi-horizon pattern across
recent trading days.
This suggests that the gains reflect both multi-horizon dependence and
intraday temporal order and lead--lag interactions that scalar volatility
summaries discard.

\subsection{Battery End-of-Life Prediction}
\label{sec:battery}

We next study early prediction of battery end-of-life (EOL) from short diagnostic pulse paths using the Argonne Cell Analysis, Modeling, and Prototyping lithium-ion dataset \citep{Paulson2022}. Following \citet{Ibraheem2025}, we retain the first usable 120-second Hybrid Pulse Power Characterization (HPPC) pulse from each cell and define EOL as the first cycle at which discharge capacity falls to 80\% of its initial value. After preprocessing and chemistry filtering, the sample contained 241 cells: 166 for training and 75 for testing.

For cell \(j\), the path contained time, current, voltage, and the temporal derivatives of current and voltage. The main depth-\(3\) signature therefore had 155 non-constant coordinates. We compare this representation with 21 handcrafted summaries computed from the same pulse and apply OLS, ridge, and LASSO to both feature sets. Holding the learner fixed across representations allowed us to separate the value of the path representation from the effect of regularization. Full feature definitions, split diagnostics, and sensitivity analyses are reported in \appref{Appen_battery}.

\autoref{tab:battery_main} reports the training and test performance of the six feature--estimator combinations.
The results reveal a clear tradeoff between representation richness and
estimation stability.
Signature-OLS achieves an almost perfect in-sample fit, substantially
outperforming Handcrafted-OLS and showing that signature features capture
detailed variation in the HPPC pulse path.
However, its poor test performance reveals severe overfitting, while
Signature-Ridge substantially reduces this instability but performs similarly
to Handcrafted-Ridge.
Among the six specifications, Signature-LASSO achieves the best out-of-sample
accuracy, with test MAE \(188.9\), test RMSE \(331.9\), and test \(R^2=0.529\).
Relative to Handcrafted-LASSO, it reduces test MAE by \(39.1\%\) and test RMSE
by \(27.7\%\).
This pattern is consistent with the sparse-regression message in
\autoref{sec_Reg}: the signature dictionary can be highly informative, but
regularization is essential for stable estimation when the feature dimension
is large relative to the sample size.

\begin{table}[h!]
\centering
\caption{Training and test performance of signature and handcrafted features. MAE and
RMSE denote mean absolute error and root mean squared error, both measured in cycles;
\(R^2\) denotes the coefficient of determination.}
\label{tab:battery_main}
\small
\setlength{\tabcolsep}{14.2pt}
\renewcommand{\arraystretch}{0.5}
 \resizebox{\linewidth}{!}{
\begin{tabular}{llcccccc}
\toprule
\multirow{2}{*}{\textbf{Feature set}}
& \multirow{2}{*}{\textbf{Method}}
& \multicolumn{3}{c}{\textbf{Training set}}
& \multicolumn{3}{c}{\textbf{Test set}} \\
\cmidrule(lr){3-5}\cmidrule(lr){6-8}
& & \textbf{MAE} & \textbf{RMSE} & \(\bm{R^2}\)
& \textbf{MAE} & \textbf{RMSE} & \(\bm{R^2}\) \\
\midrule
\multirow{3}{*}{Signature}
& OLS   & 14.4  & 29.5  & 0.997 & 616.9 & 917.8 & -2.599 \\
& Ridge & 290.1 & 443.7 & 0.409 & 314.1 & 463.1 & 0.084 \\
& LASSO & 121.3 & 210.0 & 0.868 & \textbf{188.9} & \textbf{331.9} & \textbf{0.529} \\
\cmidrule(lr){1-2}
\multirow{3}{*}{Handcrafted}
& OLS   & 219.1 & 303.8 & 0.723 & 358.4 & 638.1 & -0.740 \\
& Ridge & 354.7 & 545.6 & 0.107 & 316.2 & 464.2 & 0.079 \\
& LASSO & 349.5 & 539.3 & 0.127 & 310.4 & 459.2 & 0.099 \\
\bottomrule
\end{tabular}
}
\end{table}

\begin{figure}[h!]
\centering
\includegraphics[width=\linewidth]{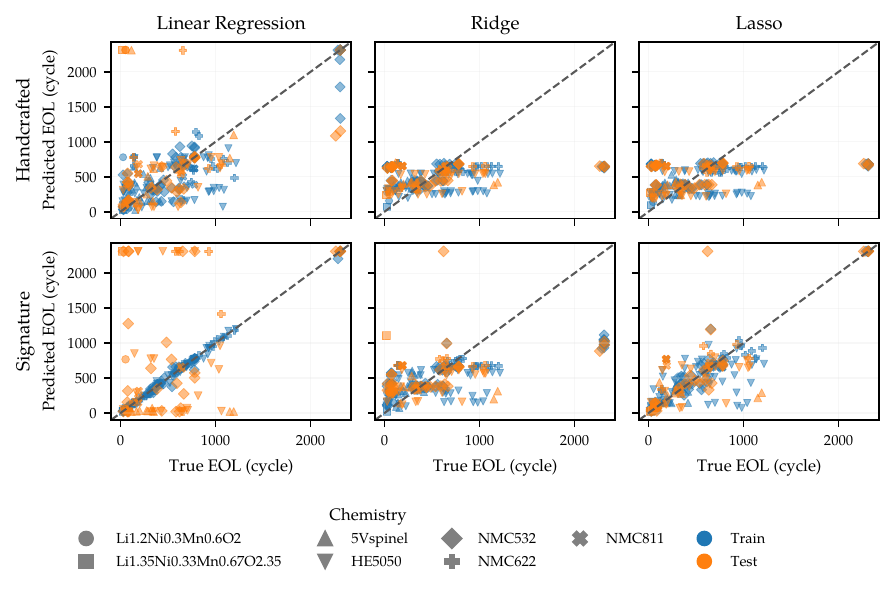}
\caption{Observed versus predicted battery EOL for the six main linear specifications.}
\label{fig:battery_parity_all}
\end{figure}

\autoref{fig:battery_parity_all} complements \autoref{tab:battery_main} by
plotting observed versus predicted EOL for the six linear specifications. Each panel
compares fitted values against the \(45^\circ\) reference line, with training and test
observations shown separately. The figure highlights the main pattern in the table:
Signature-OLS fits the training sample extremely well but does not generalize, whereas
Signature-LASSO gives the closest test-set predictions to the \(45^\circ\) line.

The contrast between Signature-OLS and Signature-LASSO is the substantive result: the pulse path is informative, but the useful information is sparse relative to the full signature dictionary. Sensitivity analyses in \appref{Appen_battery} show that prediction improves through \(K=3\) and then deteriorates as the dictionary expands, and that substantial predictive signal is already present in the early part of the retained pulse. These diagnostics support the approximation--estimation tradeoff developed in \autoref{sec_Reg}, while the modest sample size cautions against interpreting the selected coordinates as stable electrochemical mechanisms.

\subsection{Epileptic Seizure Detection}
\label{sec:seizure}

We finally consider epileptic seizure detection from three-second multichannel EEG windows using the CHB-MIT scalp EEG database \citep{ShoebGuttag2010ICML}. A window was labeled positive if it overlapped an annotated seizure interval. For each patient, we augment the EEG channels with time and temporal derivatives and compute a depth-\(2\) signature; depending on channel availability, this produced 1,406--4,032 coordinates. We fit a logistic regression classifier (LR) with balanced class weights and compare it with logistic regression and random forest (RF) based on line length, nonlinear energy, and Shannon entropy features.
The modeled cohort contained 22 subjects, and seizure windows accounted for only 0.338\% of all retained windows. We therefore use patient-specific stratified 25\%/75\% train--test splits and averaged metrics across patients rather than pooling all windows. The classification threshold for Signature + LR was fixed at 0.10 using the threshold analysis in \appref{app:seizure_appendix}. This threshold affected only the conversion of fitted probabilities to labels; AUC remains threshold-free. Complete cohort construction, feature dimensions, metric definitions, threshold sensitivity, and patient-level results are reported in \appref{app:seizure_appendix}.

Because the positive class is exceptionally rare, raw accuracy can be misleading: a classifier can attain high accuracy by predicting almost every window as non-seizure. We therefore emphasized sensitivity and balanced accuracy, which give the seizure class meaningful weight, and use AUC as a threshold-free ranking diagnostic. The literature rows in \autoref{tab:seizure_main} provide context under their reported 25\% or 50\% training rates; only the last three rows form a controlled comparison under our common patient-specific 25\% training protocol.

\begin{table}[h!]
\centering
\caption{Seizure detection performance under different training rates (TR). ACC, B-ACC,
SENS, SPEC, and AUC denote accuracy, balanced accuracy, sensitivity,
specificity, and area under the receiver operating characteristic curve, respectively;
LR and RF denote logistic regression and random forest. Performance metrics are percentages.
Literature rows use their reported protocols; our rows report patient means with patient-level
standard deviations in parentheses. Bold entries identify the best result among our three implementations.}
\label{tab:seizure_main}
\small
\setlength{\tabcolsep}{9pt}
\renewcommand{\arraystretch}{0.5}
 \resizebox{\linewidth}{!}{
\begin{tabular}{lrrrrrr}
\toprule
\multirow{2}{*}{\textbf{Method}}
& \multirow{2}{*}{\textbf{TR}}
& \multicolumn{5}{c}{\textbf{Performance Metric}} \\
\cmidrule(lr){3-7}
& & \textbf{ACC} & \textbf{B-ACC} & \textbf{SENS} & \textbf{SPEC} & \textbf{AUC} \\
\midrule
\citet{Zabihi2016}    & 50\% & 94.7 & 92.0 & 89.1 & 94.8 & -- \\
\citet{Selvakumari2019} & 50\% & 95.6 & 96.1 & 95.7 & 96.6 & -- \\
\citet{Zabihi2020}    & 25\% & 95.1 & 93.2 & 91.2 & 95.2 & 93.0 \\
\citet{Kiranyaz2014}  & 25\% & --   & 91.9 & 89.0 & 94.7 & -- \\
\midrule
Handcrafted + LR      & 25\% & 88.2 (16.0) & 89.4 (8.0)  & 90.8 (8.3)  & 88.1 (16.2) & 95.9 (3.8) \\
Handcrafted + RF      & 25\% & \textbf{99.7 (0.4)}  & 82.1 (10.1) & 64.4 (20.3) & \textbf{99.8 (0.2)}  & 96.8 (3.9) \\
Signature + LR           & 25\% & 97.5 (2.6)  & \textbf{95.8 (4.0)}  & \textbf{94.1 (7.6)}  & 97.5 (2.6)  & \textbf{98.6 (1.9)} \\
\bottomrule
\end{tabular}
}
\end{table}

Among our implementations, Handcrafted + RF has the highest accuracy and specificity, but its sensitivity is only 64.4, meaning that its apparently excellent accuracy is obtained by missing many seizure windows. Signature + LR instead attains balanced accuracy 95.8, sensitivity 94.1, specificity 97.5, and AUC 98.6. It improves sensitivity over the random forest by 29.7 percentage points while preserving high specificity. The external rows suggest that this performance is competitive with published CHB-MIT results. 

\begin{figure}
    \centering
    \includegraphics[width=\linewidth]{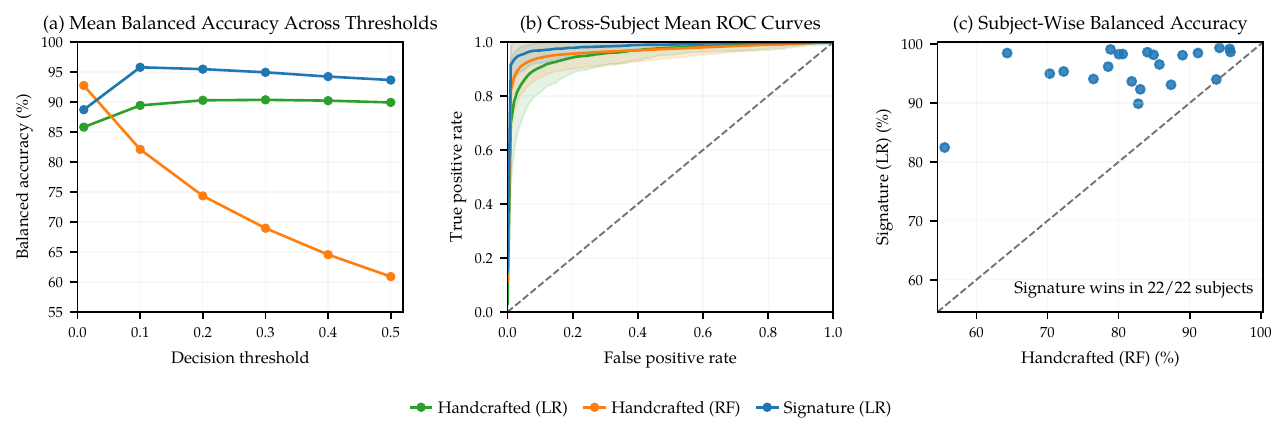}
    \caption{Diagnostic results for seizure detection. Panel (a) reports patient-averaged balanced accuracy across decision thresholds. Panel (b) reports cross-subject mean receiver operating characteristic (ROC) curves with \(\pm1\) standard-deviation bands. Panel (c) compares subject-wise balanced accuracy for the handcrafted random-forest classifier and Signature + LR at threshold 0.10.}
    \label{fig:seizure_case_study_main}
\end{figure}

\autoref{fig:seizure_case_study_main} shows that the main conclusion is not an artifact of a single aggregate score. The signature classifier reaches its maximum balanced accuracy at the reported threshold 0.10 and remains strong over a broader threshold range than the random forest. Its mean ROC curve is also highest, and its balanced accuracy exceeds that of Handcrafted + RF for every modeled subject. The patient-level result is important because cohort-level pooling would allow subjects with more windows to dominate the comparison.

\vspace{1em}

Across the three applications, the empirical findings are consistent with the theoretical message of the paper.
Truncated signatures provide rich finite-dimensional representations of path-valued covariates, but their practical performance depends on matching the learning procedure to the statistical regime.
Signature-OLS is effective when the selected feature dimension is stable; Signature-LASSO is useful when the signature dictionary is rich, but the effective representation is sparse, and Signature-Logistic extends the same idea to binary-response decisions.
These examples suggest that the proposed framework can serve as a practical guide for applying signatures to path regression problems.

\section{Conclusion}
\label{sec_conclusion}

This paper develops approximation and statistical theory for signature-based path regression.
On the approximation side, we establish \(L^2\) rates for truncated signatures.
For It\^{o} diffusions, we prove a minimax-optimal rate of order \(K^{-2\gamma}\) for a structured class of smooth coefficient functionals.
Thus, the finite-level approximation error is explicitly governed by the regularity of the target path functional.
On the statistical side, we study three procedures built on truncated signature features: Signature-OLS, Signature-LASSO, and Signature-Logistic.
This extends the existing statistical literature on signature regression \citep{Fermanian2022FunctionalSignature,GWZZ2025OR} by explicitly carrying the truncation residual through the estimation theory.
Across three real-data applications, signatures provide informative finite-dimensional representations of ordered path data and often improve upon handcrafted feature representations.

Several directions remain open.
First, the approximation theory could be extended to more general path functionals, including barriers, stopping rules, occupation times, and other nonsmooth path-dependent quantities.
Such objects arise naturally in simulation-based pricing and optimal stopping \citep{LongstaffSchwartz2001RFS}.
Second, because the number of signature coordinates grows quickly with the path dimension and truncation level, sharper theory is needed for structured, selected, and sparsified signature dictionaries.
Third, extending the statistical analysis to dependent samples, alternative losses, adaptive feature selection, and other decision-oriented objectives would further broaden the role of signatures in path-valued prediction and analytics.
Related starting points include weakly dependent functional data and high-dimensional generalized linear models \citep{HormannKokoszka2010AOS,VanDeGeer2008AOS}.

\vspace{1em}

\noindent
{\bf Acknowledgements}. The authors thank Terry Lyons, Rama Cont, Jan Obloj, and Zhuangyan Li for their helpful comments. 
We are also grateful to the DataSig group and the Oxford Mathematical and Computational Finance group for organizing seminars where this paper was presented, and to all 
participants for their 
feedback.
Blanka Horvath acknowledges financial support from UKRI1010: High order mathematical and computational infrastructure for streamed data that enhance contemporary generative and large language models. 
Wen Su acknowledges financial support from the EPSRC Centre for Doctoral Training in Mathematics of Random Systems (EP/S023925/1).
Ruixun Zhang acknowledges financial support from the National Natural Science Foundation of China (72342004,12271013), the Scientific Research Innovation Capability Support Project for Young Faculty (SRICSPYF-ZY2025170), and the National Key R\&D Program of China (2022YFA1007900).

\bibliographystyle{abbrvnat.bst}
\bibliography{reference}

\clearpage

\appendix

\clearpage
\appendix

\startcontents[appendix]

\setcounter{page}{1}
\renewcommand{\thepage}{A-\arabic{page}}

\renewcommand{\theequation}{A.\arabic{equation}}
\renewcommand{\thetable}{A.\arabic{table}}
\renewcommand{\thefigure}{A.\arabic{figure}}
\renewcommand{\thelem}{A.\arabic{lem}}
\renewcommand{\thethm}{A.\arabic{thm}}
\renewcommand{\theproposition}{A.\arabic{proposition}}

\setcounter{proposition}{0}
\setcounter{equation}{0}
\setcounter{table}{0}
\setcounter{figure}{0}
\setcounter{lem}{0}
\setcounter{thm}{0}

\begin{center}
\textbf{\LARGE{Supplementary Appendices}}
\end{center}

\vspace{1em}
\section*{Appendix Contents}
\printcontents[appendix]{}{1}{\setcounter{tocdepth}{2}}

\section{Theoretical Details for \autoref{sec:signature-approximation}}
\label{app:proof-brownian-case}

This appendix proves the approximation results in \autoref{sec:signature-approximation}.
We use $\mathbb{S}_K(\bm{X})$ to be the space spanned by the signature terms of $\bm{\widehat{X}}$ up to order $K$.

\subsection{Technical Lemmas}
\label{app:brownian-technical-lemmas}

We first record the finite-level equivalence between It\^{o} and Stratonovich Brownian signature spans.
The point is that, for Brownian motion, converting between It\^{o} and Stratonovich iterated integrals only introduces lower-level correction terms.
Hence the two finite collections of coordinates carry the same linear information at any fixed truncation level.

\begin{lem}
\label{lem:ito-stratonovich-signature-spans}
Let \(\bm B=\{\bm B_t\}_{0\le t\le T}\) be a \(d\)-dimensional standard Brownian motion and let \(\widehat{\bm B}_t=(t,\bm B_t^\top)^\top\).
For every \(K\ge0\), $\operatorname{span}
    \left\{
    S(\widehat{\bm B})^{I,\mathrm{Str}}_T:\ 0\le |I|\le K
    \right\}
    =
    \operatorname{span}
    \left\{
    S(\widehat{\bm B})^{I,\mathrm{It\hat{o}}}_T:\ 0\le |I|\le K
    \right\}$.
\end{lem}

\prf
We work on the alphabet \(\{0,1,\ldots,d\}\), where the \(0\) denotes the time coordinate and the letters \(1,\ldots,d\) denote the Brownian coordinates.
For a word \(I=(i_1,\ldots,i_\ell)\), the corresponding iterated integral is formed by integrating successively with respect to \(\mathrm d\widehat B^{i_1},\ldots,\mathrm d\widehat B^{i_\ell}\).
The It\^{o}--Stratonovich correction formula converts each Stratonovich differential into an It\^{o} differential plus quadratic-variation corrections.
When this conversion is applied repeatedly inside an iterated integral, each correction replaces two adjacent Brownian differentials by a quadratic-variation differential.

For the time-augmented Brownian path, the only non-zero quadratic covariations are
\begin{equation}
    \mathrm{d}[\widehat B^a,\widehat B^b]_t
    =
    \mathbf 1\{a=b,\ 1\le a\le d\}\,\mathrm{d}t,
    \label{eqn:brownian-quadratic-covariation}
\end{equation}
while \(\mathrm{d}[\widehat B^0,\widehat B^a]_t=0\) for all \(0\le a\le d\).
Therefore, each non-zero correction term can only arise from a pair of identical Brownian letters, and the correction replaces this pair by a single time differential.
In particular, every correction term has strictly smaller word length than the original word.

Consequently, repeated application of the It\^{o}--Stratonovich correction formula gives a triangular relation between the two coordinate systems.
More precisely, for every word \(I\) with \(|I|\le K\),
\begin{equation}
    S(\widehat{\bm B})^{I,\mathrm{Str}}_T
    =
    S(\widehat{\bm B})^{I,\mathrm{It\hat{o}}}_T
    +
    \sum_{|J|<|I|}
    c_{I,J}S(\widehat{\bm B})^{J,\mathrm{It\hat{o}}}_T,
    \label{eqn:strato-ito-triangular}
\end{equation}
where the coefficients \(c_{I,J}\) are deterministic.
The leading term is the It\^{o} iterated integral with the same word \(I\), and all remaining terms have strictly smaller length.

Now order all words with \(|I|\le K\) by nondecreasing length, and choose any fixed order inside each level.
With respect to this ordered list of coordinates, \autoref{eqn:strato-ito-triangular} shows that the linear map from the level-\(\le K\) It\^{o} signature coordinates to the level-\(\le K\) Stratonovich signature coordinates is represented by a deterministic lower-triangular matrix whose diagonal entries are equal to one.
Such a matrix is invertible.
Hence every level-\(\le K\) Stratonovich coordinate is a linear combination of level-\(\le K\) It\^{o} coordinates and, conversely, every level-\(\le K\) It\^{o} coordinate is a linear combination of level-\(\le K\) Stratonovich coordinates.
Thus the two families of coordinates span the same finite-dimensional linear space.
\QED

\subsection{Signature-Chaos Relation and Brownian Signature Approximation Rate}
\label{app:multi-brownian}

We now establish the full signature-chaos relation and the error rate of using Brownian signature to approximate Brownian functional.
We state the results in the multidimensional version.
Let \(\bm B=\{ \bm B_t\}_{0\le t\le T}\), with \(\bm B_t=(B_t^1,\ldots,B_t^d)^\top\), be a standard \(d\)-dimensional Brownian motion.
The main difference from the one-dimensional case is that a \(k\)-th chaos kernel now has \(d^k\) component functions, one for each ordered sequence of Brownian component labels.

For any \(Y\in L^2(\mathcal F_T^{\bm B})\), write its Wiener--It\^{o} chaos expansion as \(Y=\sum_{k=0}^{\infty}I_k(\bm h_k)\), where \(I_0(\bm h_0)=\mathbb E[Y]\), and for \(k\ge1\),
\begin{equation}
\bm h_k=(h_k^{i_1,\ldots,i_k})_{1\le i_1,\ldots,i_k\le d} \in L^2(\Delta_k(0,T);\mathbb R^{d^k})
\end{equation}
where $I_k(\bm h_k)
    :=
    \sum\limits_{1\le i_1,\ldots,i_k\le d}
    \int_{\Delta_k(0,T)}
    h_k^{i_1,\ldots,i_k}(t_1,\ldots,t_k)
    \,\mathrm{d}B_{t_1}^{i_1}\cdots \mathrm{d}B_{t_k}^{i_k}$.
The corresponding norm is
{\small
\begin{equation}
    \|\bm h_k\|_{L^2(\Delta_k(0,T);\mathbb R^{d^k})}^2
    :=
    \sum_{1\le i_1,\ldots,i_k\le d}
    \|h_k^{i_1,\ldots,i_k}\|_{L^2(\Delta_k(0,T))}^2 .
    \label{eqn:multi-kernel-norm}
\end{equation}}
In the multidimensional case, the \(k\)-th chaos space is
\begin{equation}
    \mathcal C_k
    :=
    \left\{
    I_k(\bm h):\bm h\in L^2(\Delta_k(0,T);\mathbb R^{d^k})
    \right\},
    \,\, k\ge1 .
    \label{eqn:multi-chaos-space}
\end{equation}

The following lemma shows that signature space is a polynomial subspace of the chaos space. 
Here \(\mathcal P_r(\Delta_k(0,T);\mathbb R^{d^k})\) denotes the space of vector-valued polynomial kernels whose components have total degree at most \(r\).
And for $k=0$, we use the convention that \(\mathcal C_0=\mathbb R\), \(I_0(c)=c\), and \(\mathcal P_r(\Delta_0(0,T);\mathbb R)=\mathbb R\).

\begin{lem}
\label{lem:multi-brownian-signature-polynomial-kernels}
For every \(0\le k\le K\), one has $\mathbb S_K(\bm B)\cap \mathcal C_k
    =
    \left\{
    I_k(\bm p):\bm p\in
    \mathcal P_{K-k}(\Delta_k(0,T);\mathbb R^{d^k})
    \right\}$.
\end{lem}

\prf
By \autoref{lem:ito-stratonovich-signature-spans}, it is enough to work with the It\^{o} signature.
For \(k\ge1\), fix a word \(I=(i_1,\ldots,i_\ell)\in\{0,1,\ldots,d\}^\ell\) with \(\ell\le K\).
Suppose that \(I\) contains exactly \(k\) nonzero letters, denoted by \(j_1,\ldots,j_k\in\{1,\ldots,d\}\) in their order of appearance.
Let \(q_0,\ldots,q_k\) be the numbers of zeros before the first nonzero letter, between consecutive nonzero letters, and after the last nonzero letter.
Then \(q_0+\cdots+q_k=\ell-k\).

After the \(k\) Brownian integration times are fixed, the zero letters are integrated out over the gaps between those times, exactly as in the one-dimensional proof.
Therefore, the corresponding It\^{o} signature coordinate can be written as
\begin{equation}
    S(\widehat{\bm B})^{I,\mathrm{It\hat{o}}}_T
    =
    \int_{\Delta_k(0,T)}
    p_I(t_1,\ldots,t_k)
    \,\mathrm{d}B_{t_1}^{j_1}\cdots \mathrm{d}B_{t_k}^{j_k},
    \label{eqn:multi-word-to-polynomial-kernel}
\end{equation}
where $p_I(t_1,\ldots,t_k)
    =
    \frac{t_1^{q_0}}{q_0!}
    \frac{(t_2-t_1)^{q_1}}{q_1!}
    \cdots
    \frac{(t_k-t_{k-1})^{q_{k-1}}}{q_{k-1}!}
    \frac{(T-t_k)^{q_k}}{q_k!}$.
Thus this coordinate belongs to the \(k\)-th chaos and has a vector-valued polynomial kernel whose only nonzero component is indexed by \((j_1,\ldots,j_k)\), with total degree \(\ell-k\le K-k\).

Conversely, take any \(\bm p\in\mathcal P_{K-k}(\Delta_k(0,T);\mathbb R^{d^k})\).
It is enough to consider one component and one monomial, because the vector-valued polynomial space is spanned by such componentwise monomials.
Fix a component \((j_1,\ldots,j_k)\) and introduce the gap variables \(y_0=t_1\), \(y_r=t_{r+1}-t_r\) for \(1\le r\le k-1\), and \(y_k=T-t_k\).
Every polynomial in \((t_1,\ldots,t_k)\) of degree at most \(K-k\) is a linear combination of monomials \(y_0^{q_0}\cdots y_k^{q_k}\) with \(q_0+\cdots+q_k\le K-k\).
Such a monomial is generated, up to a deterministic nonzero constant, by the word
\begin{equation}
    (\underbrace{0,\ldots,0}_{q_0},j_1,
    \underbrace{0,\ldots,0}_{q_1},j_2,\ldots,
    j_k,\underbrace{0,\ldots,0}_{q_k}),
    \label{eqn:multi-generating-word}
\end{equation}
whose length is at most \(K\).
Hence every vector-valued polynomial kernel of degree at most \(K-k\) is generated by level-\(K\) time-augmented Brownian signature coordinates.

Finally, if a finite linear combination of level-\(K\) signature coordinates lies in \(\mathcal C_k\), then all its chaos components of orders different from \(k\) vanish by the orthogonality of Wiener chaoses.
Its \(k\)-th chaos component is therefore generated by words with exactly \(k\) nonzero letters, and hence by vector-valued polynomial kernels of degree at most \(K-k\).
This gives the result.
\QED

By \autoref{lem:multi-brownian-signature-polynomial-kernels} and the orthogonality of different Wiener chaoses, the level-\(K\) residual \(\xi_K^{\bm B}(Y):=Y-\operatorname{Proj}_{\mathbb S_K(\bm B)}Y\) satisfies
\begin{equation}
\begin{aligned}
    \mathbb E|\xi_K^{\bm B}(Y)|^2
    = 
    \sum_{k>K}
    \|\bm h_k\|_{L^2(\Delta_k(0,T);\mathbb R^{d^k})}^2 
    +
    \sum_{k=0}^{K}
    \inf_{\bm p\in\mathcal P_{K-k}}
    \|\bm h_k-\bm p\|_{L^2(\Delta_k(0,T);\mathbb R^{d^k})}^2.
\end{aligned}
    \label{eqn:multi-brownian-error-decomposition}
\end{equation}

For \(k\ge1\), define the normalized kernel
\begin{equation}
    \widetilde{\bm h}_k
    :=
    \begin{cases}
    \|\bm h_k\|^{-1}_{L^2(\Delta_k(0,T);\mathbb R^{d^k})}\bm h_k,
    & \text{if } \bm{h}_k \neq \bm{0},\\
    0,
    & \text{if } \bm{h}_k=\bm{0}.
    \end{cases}
    \label{eqn:multi-normalized-kernel}
\end{equation}

\begin{thm}
\label{thm:multi-brownian-signature-rate}
Let \(Y\in L^2(\mathcal F_T^{\bm B})\) admit the chaos expansion \(Y=\sum_{k=0}^{\infty}I_k(\bm h_k)\).
Assume that there exist constants \(m\in\mathbb N_+\), \(\epsilon>0\), and \(\gamma>1/2\) such that, for \(k\ge1\), \(\|\bm h_k\|_{L^2(\Delta_k(0,T);\mathbb R^{d^k})}^2 =O\{(1+k)^{-(m+1+\epsilon)}\}\), and \(\|\widetilde{\bm h}_k\|_{H^\gamma(\Delta_k(0,T);\mathbb R^{d^k})}=O(a_k)\), where the constants in vector-valued polynomial approximation on \(\Delta_k(0,T)\) are absorbed into \(a_k\), and \(\{a_k\}_{k\ge1}\) satisfies \(\sum_{k\ge1}a_k^2(1+k)^{-(m+1+\epsilon)}<\infty\).
Let \(A_K:=\max_{1\le k\le K}a_k\).
Then
\begin{equation}
    \mathbb E|\xi_K^{\bm B}(Y)|^2
    =
    O\left(
    K^{-2\gamma}
    +
    (A_K^2+K)K^{-(m+1+\epsilon)}
    \right).
    \label{eqn:multi-brownian-signature-rate}
\end{equation}
\end{thm}

\prf
By \autoref{eqn:multi-brownian-error-decomposition}, it suffices to control the chaos tail and the polynomial approximation error.
The tail condition gives
\begin{equation}
    \sum_{k>K}
    \|\bm h_k\|_{L^2(\Delta_k(0,T);\mathbb R^{d^k})}^2
    =
    O\left((K+1)^{-(m+\epsilon)}\right).
    \label{eqn:multi-chaos-tail-bound}
\end{equation}
Then, since \(\bm h_k=\|\bm h_k\|_{L^2(\Delta_k(0,T);\mathbb R^{d^k})}\widetilde{\bm h}_k\), polynomial approximation on the simplex gives $\inf \limits_{\bm p\in\mathcal P_{K-k}}
    \|\bm h_k-\bm p\|_{L^2(\Delta_k(0,T);\mathbb R^{d^k})}^2
    \le
    C\|\bm h_k\|_{L^2(\Delta_k(0,T);\mathbb R^{d^k})}^2
    a_k^2(K-k+1)^{-2\gamma}$.
Therefore,
\begin{equation}
\begin{aligned}
    \sum_{k=0}^{K}
    \inf \limits_{\bm p\in\mathcal P_{K-k}}
    \|\bm h_k-\bm p\|_{L^2(\Delta_k(0,T);\mathbb R^{d^k})}^2
      \le
    C\sum_{k=1}^{K}
    (1+k)^{-(m+1+\epsilon)}a_k^2(K-k+1)^{-2\gamma}.
\end{aligned}
    \label{eqn:multi-low-order-bound}
\end{equation}
Splitting the last sum into \(1\le k\le K/2\) and \(K/2<k\le K\), the first part is \(O(K^{-2\gamma})\) by the summability assumption, while the second part is \(O(A_K^2K^{-(m+1+\epsilon)})\) because \(\gamma>1/2\).
Combining these estimates with \autoref{eqn:multi-chaos-tail-bound} gives \autoref{eqn:multi-brownian-signature-rate}.
\QED

\subsection{Proof of \autoref{thm:smooth-coefficient-functional}}

\prf
We split the proof into five steps.
The first step localizes the diffusion path and constructs polynomial approximations on a compact set.
The second step shows that the resulting polynomial coefficient functional is a finite-level signature functional.
The third step controls the error on the localization event.
The fourth step controls the tail event using the boundedness of the diffusion coefficients.
The last step verifies the sharpness claim by reducing to the Brownian benchmark.

\noindent \textbf{Step 1: localization and polynomial approximation.}
Choose a constant \(\beta\) such that $\frac12<\beta<1-\frac{\gamma}{2\gamma+1}$.
Such a \(\beta\) exists because \( 1-\frac{\gamma}{2\gamma+1} = \frac{\gamma+1}{2\gamma+1} > \frac12 .
\) Let
\begin{equation}
    m_K:=K^\beta,
    \quad
    R_T:=\sup_{0\le t\le T}\|\bm X_t\|_\infty,
    \quad
    A_K:=\{R_T\le m_K\}.
\end{equation}
We first consider \(K\) sufficiently large that \(m_K\ge\|\bm x_0\|_\infty\); the finitely many remaining values are handled at the end.

On the localization event \(A_K\), the path stays in \(D_K:=[0,T]\times[-m_K,m_K]^d\).
We first note that the boundedness of \(\bm b\) and \(\bm\sigma\) gives a sub-Gaussian tail bound for \(R_T\).
Indeed,
\begin{equation}
    \bm X_t
    =
    \bm x_0
    +
    \int_0^t
    \bm b(u,\bm X_u)\,\mathrm du
    +
    \int_0^t
    \bm\sigma(u,\bm X_u)\,\mathrm d\bm B_u .
\end{equation}
Since \(\bm b\) is uniformly bounded,
\begin{equation}
    R_T
    \le
    \|\bm x_0\|_\infty+RT
    +
    \sup_{0\le t\le T}
    \left\|
        \int_0^t
        \bm\sigma(u,\bm X_u)\,\mathrm d\bm B_u
    \right\|_\infty .
\end{equation}
For each coordinate of the martingale term, the quadratic variation is bounded by \(R^2T\).
The exponential martingale inequality and a union bound imply that there exist constants \(C,c>0\) such that
\begin{equation}
    \mathbb P(A_K^c)
    =
    \mathbb P(R_T>m_K)
    \le
    C e^{-cm_K^2}
    =
    Ce^{-cK^{2\beta}}.
    \label{eq:AK-tail-proof}
\end{equation}
The same tail bound implies the moment estimate $ \mathbb E R_T^p \le C^p p^{p/2}$, for $p\ge 1$.

We now approximate \(h_0,\ldots,h_d\) on the compact rectangle \(D_K\).
Rescale this rectangle to \([0,1]\times[-1,1]^d\) by setting \( \widetilde f_K(u,y):=f(Tu,m_Ky).
\) The standard Jackson approximation theorem \citep{BeutelGonska2003} on the fixed rectangle gives polynomial approximants in terms of the derivatives of \(\widetilde f_K\).
Time derivatives only produce powers of \(T\), while spatial derivatives of order \(r\) produce the scaling factor \(m_K^r\).
For \(h_0\), there exists a polynomial \(p_{0,K}\) in \((t,x)\), of total degree at most \(K-1\), such that
\begin{equation}
\begin{aligned}
    \|h_0-p_{0,K}\|_{L^\infty(D_K)}
     \le
    C
    \left\{
        K^{-\gamma}
        +
        m_K^{2\gamma+1}K^{-(2\gamma+1)}
    \right\}
    \mathcal H_\gamma(\bm h).
\end{aligned}
    \label{eq:h0-local-approx-proof}
\end{equation}
For \(a=1,\ldots,d\), applying the same approximation argument simultaneously to \(h_a\) and its first spatial derivatives gives a polynomial \(p_{a,K}\), of total degree at most \(K-1\), such that
\begin{equation}
\begin{aligned}
    &\|h_a-p_{a,K}\|_{L^\infty(D_K)}
    +
    \|\nabla_x h_a-\nabla_xp_{a,K}\|_{L^\infty(D_K)}
    \\
    &\qquad\le
    C
    \left\{
        K^{-\gamma}
        +
        m_K^{2\gamma+1}K^{-(2\gamma+1)}
    \right\}
    \mathcal H_\gamma(\bm h).
\end{aligned}
    \label{eq:ha-local-approx-proof}
\end{equation}

By the choice of \(\beta\), $m_K^{2\gamma+1}K^{-(2\gamma+1)}
    =
    K^{-(2\gamma+1)(1-\beta)}
    =
    o(K^{-\gamma})$.
Therefore, $K^{-\gamma} + m_K^{2\gamma+1}K^{-(2\gamma+1)} \le C K^{-\gamma}$.
Let \(\delta_{a,K}:=h_a-p_{a,K}\) for $a=0,\cdots,d$ and define
\begin{equation}
    Y_K
    :=
    \sum_{a=0}^d
    \int_0^T
    p_{a,K}(t,\bm X_t)
    \circ\mathrm d\widehat X_t^a .
    \label{eq:YK-proof}
\end{equation}

\noindent \textbf{Step 2: polynomial coefficient functionals are finite-level signature functionals.}
We show that \(Y_K\in\mathbb S_K(\bm X)\).
Fix \(a\in\{0,\ldots,d\}\).
Since \(p_{a,K}\) is a polynomial in \((t,x)\) of total degree at most \(K-1\), and since \(\bm X_0=\bm x_0\) is deterministic, each monomial in \(p_{a,K}(t,\bm X_t)\) can be written as a polynomial in the first-level prefix signature coordinates of \(\widehat{\bm X}\).
Indeed, $t
    =
    \left\langle e_0,S(\widehat{\bm X})_{0,t}\right\rangle$, and
\begin{equation}
    X_t^b
    =
    x_0^b+
    \left\langle e_b,S(\widehat{\bm X})_{0,t}\right\rangle,
    \quad b=1,\ldots,d .
    \label{eq:first-level-prefix-proof-new}
\end{equation}
By the shuffle identity, every product of such prefix signature coordinates can be rewritten as a linear functional of the prefix signature.
Hence, writing \(T^{(K-1)}(\mathbb R^{d+1}):=\bigoplus_{r=0}^{K-1}(\mathbb R^{d+1})^{\otimes r}\), there exists \(q_{a,K}\in T^{(K-1)}(\mathbb R^{d+1})\) such that $p_{a,K}(t,\bm X_t)
    =
    \left\langle
        q_{a,K},
        S(\widehat{\bm X})_{0,t}
    \right\rangle$.
Using the recursive definition of the Stratonovich signature,
\begin{equation}
    \int_0^T
    p_{a,K}(t,\bm X_t)
    \circ\mathrm d\widehat X_t^a
    =
    \left\langle
        q_{a,K}\otimes e_a,
        S(\widehat{\bm X})_{0,T}
    \right\rangle .
    \label{eq:append-letter-proof-new}
\end{equation}
The tensor \(q_{a,K}\otimes e_a\) has level at most \(K\).
Summing over \(a=0,\ldots,d\), we obtain $Y_K\in\mathbb S_K(\bm X)$.
Therefore,
\begin{equation}
    \mathbb E|\xi_K^{\bm X}(Y)|^2
    =
    \inf_{Z\in\mathbb S_K(\bm X)}
    \mathbb E|Y-Z|^2
    \le
    \mathbb E|Y-Y_K|^2 .
    \label{eq:projection-proof-new}
\end{equation}

\noindent \textbf{Step 3: error on the localization event.}
We split
\begin{equation}
\begin{aligned}
    \mathbb E|Y-Y_K|^2
    = 
    \mathbb E\left[|Y-Y_K|^2\mathbf 1_{A_K}\right] +
    \mathbb E\left[|Y-Y_K|^2\mathbf 1_{A_K^c}\right]
    =:
    I_K^{\mathrm{loc}}+I_K^{\mathrm{tail}}.
\end{aligned}
\end{equation}
On \(A_K\), the path stays in the localization box.
For the time channel, one has
\begin{equation}
\begin{aligned}
    \mathbb E
    \left[
        \left|
            \int_0^T
            \delta_{0,K}(t,\bm X_t)\,\mathrm dt
        \right|^2
        \mathbf 1_{A_K}
    \right]
    \le
    T^2
    \|\delta_{0,K}\|_{L^\infty(D_K)}^2          
    \le
    C K^{-2\gamma}\mathcal H_\gamma(\bm h)^2 .
\end{aligned}
    \label{eq:local-time-channel-proof}
\end{equation}

For \(a=1,\ldots,d\), the Stratonovich--It\^{o} conversion gives
\begin{equation}
\begin{aligned}
    \int_0^T
    \delta_{a,K}(t,\bm X_t)
    \circ\mathrm dX_t^a
    = 
    \int_0^T
    \delta_{a,K}(t,\bm X_t)
    \,\mathrm dX_t^a 
    +
    \frac12
    [\delta_{a,K}(\cdot,\bm X_\cdot),X^a]_T .
\end{aligned}
    \label{eq:strat-ito-proof-new}
\end{equation}
Let
\(\eta_K:=\inf\{t\in[0,T]:\|\bm X_t\|_\infty>m_K\}\wedge T\),
where \(\inf\varnothing=\infty\).
By continuity, the stopped path remains in \(D_K\), and
\(\eta_K=T\) on \(A_K\).
After decomposing \(\mathrm dX_t^a\) into its drift and martingale
parts, Cauchy--Schwarz and It\^{o} isometry give
\begin{equation}
\begin{aligned}
    \mathbb E
    \left[
        \left|
            \int_0^T
            \delta_{a,K}(t,\bm X_t)\,\mathrm dX_t^a
        \right|^2
        \mathbf 1_{A_K}
    \right]
    &\le
    \mathbb E
    \left|
        \int_0^{\eta_K}
        \delta_{a,K}(t,\bm X_t)\,\mathrm dX_t^a
    \right|^2 \\
    &\le
    C\|\delta_{a,K}\|_{L^\infty(D_K)}^2
    \le
    CK^{-2\gamma}\mathcal H_\gamma(\bm h)^2 .
\end{aligned}
\label{eq:local-ito-proof}
\end{equation}

For the quadratic covariation term, only spatial derivatives contribute:
{
\begin{equation}
    [\delta_{a,K}(\cdot,\bm X_\cdot),X^a]_T
    =
    \int_0^T
    \sum_{r=1}^d
    \partial_{x_r}\delta_{a,K}(t,\bm X_t)
    \,\mathrm d[X^r,X^a]_t .
    \label{eq:qv-chain-proof-new}
\end{equation}}
Moreover, $ \mathrm d[X^r,X^a]_t = \sum_{\ell=1}^q \sigma_{r\ell}(t,\bm X_t) \sigma_{a\ell}(t,\bm X_t) \,\mathrm dt$.
Therefore, by boundedness of \(\bm\sigma\), and \autoref{eq:ha-local-approx-proof}, one has
\begin{equation}
    \mathbb E
    \left[
        \left|
            [\delta_{a,K}(\cdot,\bm X_\cdot),X^a]_T
        \right|^2
        \mathbf 1_{A_K}
    \right]
    \le
    C K^{-2\gamma}\mathcal H_\gamma(\bm h)^2 .
    \label{eq:local-qv-proof}
\end{equation}
Combining \autoref{eq:local-time-channel-proof},
\autoref{eq:strat-ito-proof-new},
\autoref{eq:local-ito-proof}, and
\autoref{eq:local-qv-proof}, and using
\(\left|\sum_{a=0}^d z_a\right|^2
\le (d+1)\sum_{a=0}^d|z_a|^2\), yields
\begin{equation}
    I_K^{\mathrm{loc}}
    \le
    C K^{-2\gamma}\mathcal H_\gamma(\bm h)^2 .
    \label{eq:local-error-final-proof}
\end{equation}

\noindent \textbf{Step 4: tail error.}
It remains to control \(I_K^{\mathrm{tail}}\). 
Since \(h_0,\ldots,h_d\) and \(\nabla_x h_1,\ldots,\nabla_x h_d\) are bounded by \(\mathcal H_\gamma(\bm h)\), boundedness of \(\bm b\) and \(\bm\sigma\), the Stratonovich--It\^{o} conversion, and the Burkholder--Davis--Gundy inequality give $\left(\mathbb E|Y|^4\right)^{1/2}\le C\mathcal H_\gamma(\bm h)^2$ and
\begin{equation}
    \mathbb E\left[|Y|^2\mathbf 1_{A_K^c}\right]
    \le
    \left(\mathbb E|Y|^4\right)^{1/2}
    \mathbb P(A_K^c)^{1/2}
    \le
    C\mathcal H_\gamma(\bm h)^2
    e^{-cK^{2\beta}},
\end{equation}
where the last inequality follows from \autoref{eq:AK-tail-proof}.

We now bound \(Y_K\) on the tail.
From \autoref{eq:h0-local-approx-proof} and
\autoref{eq:ha-local-approx-proof},
\begin{equation}
    \max_{0\le a\le d}
    \|p_{a,K}\|_{L^\infty(D_K)}
    +
    \max_{1\le a\le d}
    \|\nabla_xp_{a,K}\|_{L^\infty(D_K)}
    \le
    C\mathcal H_\gamma(\bm h).
\end{equation}
After rescaling the localization box to the unit box, a
tensor-product Chebyshev expansion and the fact that
\(p_{a,K}\) has total degree at most \(K-1\) imply that,
for all \((t,x)\in D\),
\begin{equation}
\begin{aligned}
    &\max_{0\le a\le d}|p_{a,K}(t,x)|
    +
    \max_{1\le a\le d}\|\nabla_xp_{a,K}(t,x)\| \\
    &\qquad\le
    C\mathcal H_\gamma(\bm h)
    \exp(CK\log K)
    \left(1+\frac{\|x\|_\infty}{m_K}\right)^K .
\end{aligned}
\label{eq:pK-growth-proof}
\end{equation}
Using \autoref{eq:pK-growth-proof}, boundedness of \(\bm b,\bm\sigma\), the Stratonovich--It\^{o} conversion, and BDG,
{
\begin{equation}
    \left(\mathbb E|Y_K|^4\right)^{1/2}
    \le
    C\mathcal H_\gamma(\bm h)^2
    e^{CK\log K}
    \left[
        \mathbb E
        \left(1+\frac{R_T}{m_K}\right)^{4K}
    \right]^{1/2}.
    \label{eq:YK-fourth-before-mgf-proof}
\end{equation}}
By \(\log(1+u)\le u\) and the sub-Gaussian moment estimate above,
\begin{equation}
    \mathbb E
    \left(1+\frac{R_T}{m_K}\right)^{4K}
    \le
    \mathbb E\exp\left(\frac{4KR_T}{m_K}\right)
    \le
    C\exp\left(C\frac{K^2}{m_K^2}\right).
\end{equation}
Thus $\left(\mathbb E|Y_K|^4\right)^{1/2} \le C\mathcal H_\gamma(\bm h)^2 \exp\left(CK\log K+C\frac{K^2}{m_K^2}\right)$.
By Cauchy--Schwarz and \autoref{eq:AK-tail-proof},
\begin{equation}
\begin{aligned}
    &\mathbb E\left[
        |Y_K|^2\mathbf 1_{A_K^c}
    \right]
    \le
    \left(\mathbb E|Y_K|^4\right)^{1/2}
    \mathbb P(A_K^c)^{1/2}                                      \\
    &\quad \le
    C\mathcal H_\gamma(\bm h)^2
    \exp\left(
        CK\log K
        +
        C K^{2-2\beta}
        -
        cK^{2\beta}
    \right).
\end{aligned}
    \label{eq:YK-tail-proof}
\end{equation}
Since \(\beta>1/2\), the negative term \(K^{2\beta}\) dominates both \(K\log K\) and \(K^{2-2\beta}\).
Therefore,
\begin{equation}
    \mathbb E\left[
        |Y_K|^2\mathbf 1_{A_K^c}
    \right]
    =
    o(K^{-2\gamma})\mathcal H_\gamma(\bm h)^2 .
    \label{eq:YK-tail-final-proof}
\end{equation}
Using \(|Y-Y_K|^2\le 2|Y|^2+2|Y_K|^2\), the preceding bound for \(Y\), and \autoref{eq:YK-tail-final-proof}, we obtain
\(I_K^{\mathrm{tail}}
=o(K^{-2\gamma})\mathcal H_\gamma(\bm h)^2\).
Combining this with the decomposition above, \autoref{eq:projection-proof-new}, and \autoref{eq:local-error-final-proof} gives \(\mathbb E|\xi_K^{\bm X}(Y)|^2 \le C K^{-2\gamma}\mathcal H_\gamma(\bm h)^2\) for all sufficiently large \(K\).
Increasing \(C\) covers the finitely many remaining values \(K\ge2\).

\noindent \textbf{Step 5: sharpness.}
It remains to establish the pair-class minimax lower bound.
For the fixed dimensions \(d,q\ge1\), initial value \(\bm x_0\), and class radius \(R>0\), let \(\bm e_1^{(m)}\) denote the first canonical basis vector of \(\mathbb R^m\), and take
\(\bm b=0\) and
\(\bm\sigma=R\bm e_1^{(d)}(\bm e_1^{(q)})^\top\).
Then
\(\bm X_t=\bm x_0+R\bm e_1^{(d)}B_t^1\), and this diffusion satisfies the coefficient bound defining \(\mathcal U_{\gamma,R}\).

Consider \(h_1(t,x)=g(t)\) and \(h_a=0\) for \(a\ne1\), where \(g\in C^\infty([0,T])\) satisfies
\(\sum_{\ell=0}^{\gamma}
\|g^{(\ell)}\|_{L^\infty(0,T)}\le R\).
Since \(g\) does not depend on \(x\),
\(\mathcal H_\gamma(\bm h)
=\sum_{\ell=0}^{\gamma}
\|g^{(\ell)}\|_{L^\infty(0,T)}\le R\), so every such pair belongs to \(\mathcal U_{\gamma,R}\).
Moreover,
\(F_{\bm h}(\bm X)
=R\int_0^Tg(t)\,\mathrm dB_t^1\).

The deterministic initial translation does not affect the signature, the constant spatial coordinates contribute no signature information, and every nonzero signature coordinate containing \(m\) occurrences of the first spatial letter is multiplied by \(R^m\).
Since \(R>0\), these rescalings do not change the signature span, and hence
\(\mathbb S_K(\bm X)=\mathbb S_K(B^1)\).
Therefore, applying
\autoref{lem:multi-brownian-signature-polynomial-kernels} to the first Wiener chaos gives
\begin{equation}
    \mathbb E
    \left[
        \left|
            \xi_K^{\bm X}
            \bigl(F_{\bm h}(\bm X)\bigr)
        \right|^2
    \right]
    =
    R^2
    \inf_{p\in\mathcal P_{K-1}(0,T)}
    \|g-p\|_{L^2(0,T)}^2 .
    \label{eq:proof-minimax-reduction-polynomial}
\end{equation}

We next establish the required polynomial approximation lower bound within this smooth subclass.
Fix a nonzero function
\(\psi\in C_c^\infty(0,1)\), extended by zero outside \((0,1)\), and set
\(\overline T:=\min\{T,1\}\),
\(\delta_K:=\overline T/(2K)\), and
\(A_\gamma:=
\sum_{\ell=0}^{\gamma}
\|\psi^{(\ell)}\|_{L^\infty(\mathbb R)}\).
For \(j=1,\ldots,2K\), define
\(\phi_{j,K}(t):=
R A_\gamma^{-1}\delta_K^\gamma
\psi((t-(j-1)\delta_K)/\delta_K)\).
These functions have disjoint supports.
Thus, for every
\(\bm\varepsilon=(\varepsilon_1,\ldots,\varepsilon_{2K})
\in\{-1,1\}^{2K}\), the function
\(g_{\bm\varepsilon,K}
:=\sum_{j=1}^{2K}\varepsilon_j\phi_{j,K}\)
belongs to \(C^\infty([0,T])\) and satisfies
\(\sum_{\ell=0}^{\gamma}
\|g_{\bm\varepsilon,K}^{(\ell)}\|_{L^\infty(0,T)}
\le R\).

Let \(\Pi_K\) be the \(L^2(0,T)\)-orthogonal projection onto
\(\mathcal P_{K-1}(0,T)\), and let
\(b_K:=\|\phi_{j,K}\|_{L^2(0,T)}^2
=R^2A_\gamma^{-2}\delta_K^{2\gamma+1}
\|\psi\|_{L^2(0,1)}^2\).
Since
\(\{\phi_{j,K}/\sqrt{b_K}\}_{j=1}^{2K}\)
is an orthonormal family and
\(\dim\mathcal P_{K-1}(0,T)=K\), Rademacher averaging and Bessel's inequality yield
\begin{equation}
\begin{aligned}
    \sup_{\substack{g\in C^\infty([0,T])\\
    \sum_{\ell=0}^{\gamma}
    \|g^{(\ell)}\|_{L^\infty(0,T)}\le R}}
    \inf_{p\in\mathcal P_{K-1}(0,T)}
    \|g-p\|_{L^2(0,T)}^2                                      
    &\ge
    \mathbb E_{\bm\varepsilon}
    \left[
        \inf_{p\in\mathcal P_{K-1}(0,T)}
        \|g_{\bm\varepsilon,K}-p\|_{L^2(0,T)}^2
    \right]                                                     \\
    &=
    \sum_{j=1}^{2K}
    \|(I-\Pi_K)\phi_{j,K}\|_{L^2(0,T)}^2 \\
    &\ge
    Kb_K
    =
    c_{\gamma,T}R^2K^{-2\gamma},
\end{aligned}
\label{eq:proof-minimax-lower-bound}
\end{equation}
where
\(c_{\gamma,T}
:=
A_\gamma^{-2}\|\psi\|_{L^2(0,1)}^2
\overline T^{\,2\gamma+1}2^{-(2\gamma+1)}>0\).

Combining
\autoref{eq:proof-minimax-reduction-polynomial} and
\autoref{eq:proof-minimax-lower-bound} gives
\(\sup_{(\bm X,\bm h)\in\mathcal U_{\gamma,R}}
\mathbb E|
\xi_K^{\bm X}(F_{\bm h}(\bm X))|^2
\ge c_{\gamma,T}R^4K^{-2\gamma}\).
Since the upper bound established above is uniform over
\(\mathcal U_{\gamma,R}\), with a constant depending only on
\(T,d,q,\bm x_0,\gamma\), and \(R\), it follows that
\(\sup_{(\bm X,\bm h)\in\mathcal U_{\gamma,R}}
\mathbb E|
\xi_K^{\bm X}(F_{\bm h}(\bm X))|^2
\asymp K^{-2\gamma}\).
This proves \autoref{thm:smooth-coefficient-functional}.
\QED

\subsection{Exponential Rates for Analytic Coefficient Functionals}
\label{app:analytic-rate}

We now give an analytic analogue of \autoref{thm:smooth-coefficient-functional}.
Because the state space is unbounded, ordinary analyticity on a fixed complex neighborhood is not sufficient for a geometric rate on the expanding localization boxes used below.
We instead impose analyticity on a fixed complex neighborhood of the time interval and a controlled entire extension in the state variables.
Throughout this subsection, \(T\), \(d\), \(q\), and the deterministic initial value \(\bm x_0\) are fixed, and \(\mathbb S_K(\bm X)\) denotes the span of all time-augmented signature coordinates through level \(K\).
The result below concerns approximation by this full space; it does not automatically apply to an arbitrary selected subset of signature coordinates.

For \(r>1\), let \(\mathcal E_r\) denote the filled closed Bernstein ellipse for \([-1,1]\) with parameter \(r\), and let \(\mathcal E_r^{[0,1]}:=\{(1+z)/2:z\in\mathcal E_r\}\) be its affine image corresponding to \([0,1]\).
For \(\rho_0>1\), \(\tau>0\), and \(R>0\), define \(\mathcal U_{\rho_0,\tau,R}^{\mathrm{an}}\) as the collection of pairs \((\bm X,\bm h)\) satisfying the following conditions.
The process \(\bm X\) solves \autoref{eq:smooth-coefficient-X-SDE} from \(\bm x_0\), where \(\bm b\) and \(\bm\sigma\) are globally Lipschitz and
\(\sup_{(t,x)\in D}\{\|\bm b(t,x)\|+\|\bm\sigma(t,x)\|_F\}\le R\).
Each \(h_a\) is real-valued on \(D\), and \(h_a(T\cdot,\cdot)\) admits a jointly holomorphic extension to \(V_{\rho_0}\times\mathbb C^d\), where \(V_{\rho_0}\) is an open set containing \(\mathcal E_{\rho_0}^{[0,1]}\), such that
\begin{equation}
\label{eq:scale-analytic-class}
    \max_{0\le a\le d}
    \sup_{\substack{z_0\in\mathcal E_{\rho_0}^{[0,1]}\\
                    z\in\mathbb C^d}}
    e^{-\tau\|\operatorname{Im}z\|_1}
    |h_a(Tz_0,z)|
    \le
    R.
\end{equation}
The condition is asymmetric in time and state because the time interval remains fixed whereas the state-space localization box expands with \(K\).
On the real state space it gives \(|h_a(t,x)|\le R\).
Cauchy estimates on smaller product neighborhoods show that, for every fixed \(\gamma\), \(\mathcal H_\gamma(\bm h)\le C_{\gamma,\rho_0,\tau,d,T}R\); hence, after changing the class bound by a fixed multiplicative constant, these coefficients also satisfy the smoothness conditions above.
The coefficient condition is nondegenerate.
For example, it contains \(h_a(t,x)=g_a(t)\sin(\theta_a^\top x)\) and \(h_a(t,x)=g_a(t)\cos(\theta_a^\top x)\) when \(g_a(T\cdot)\) is holomorphic on \(V_{\rho_0}\), is real-valued on \([0,1]\), and satisfies \(\sup_{z_0\in\mathcal E_{\rho_0}^{[0,1]}}|g_a(Tz_0)|\le R\), while \(\theta_a\in\mathbb R^d\) and \(\|\theta_a\|_\infty\le\tau\).
It also contains finite real sine--cosine combinations whose time coefficients satisfy the same holomorphy and reality conditions, whose frequencies satisfy the same bound, and whose coefficient-weighted time envelopes sum to at most \(R\).

\begin{proposition}
\label{prop:analytic-signature-rate}
Suppose \((\bm X,\bm h)\in\mathcal U_{\rho_0,\tau,R}^{\mathrm{an}}\).
Then, for every \(1<\rho<\rho_0\), there exists a constant \(C_{\rho}>0\), depending only on \(T,d,q,\bm x_0,R,\rho_0,\rho\), and \(\tau\), such that, for every \(K\ge2\),
\begin{equation}
\label{eq:analytic-signature-upper-bound}
    \inf_{Z\in\mathbb S_K(\bm X)}
    \mathbb E\left|F_{\bm h}(\bm X)-Z\right|^2
    \le
    C_{\rho}R^2\rho^{-2K}.
\end{equation}
\end{proposition}

The strict inequality \(\rho<\rho_0\) allows the subexponential factor generated by localization to be absorbed into the geometric rate; no boundary rate at \(\rho=\rho_0\) is claimed.

\prf
Set \(\bar\rho:=(\rho+\rho_0)/2\) and \(\beta:=3/4\), so that
\(1<\rho<\bar\rho<\rho_0\).
Throughout the proof, \(C,c\), and similarly subscripted constants may
change from line to line and may depend on
\(T,d,q,\bm x_0,R,\rho_0,\rho\), and \(\tau\), but not on \(K\) or on
the particular pair \((\bm X,\bm h)\) in the class.
Set
\begin{equation}
\label{eq:analytic-localization-definitions}
\begin{aligned}
    &m_K:=1+\|\bm x_0\|_\infty+K^\beta,
    \qquad
    X_T^\ast:=\sup_{0\le t\le T}\|\bm X_t\|_\infty,
    \qquad
    A_K:=\{X_T^\ast\le m_K\},
    \\
    &D_K:=[0,T]\times[-2m_K,2m_K]^d,
    \qquad
    \eta_K:=
    \inf\left\{
        t\in[0,T]:
        \|\bm X_t\|_\infty\ge2m_K
    \right\}\wedge T,
\end{aligned}
\end{equation}
where \(\inf\varnothing:=\infty\).
For \(j=1,\ldots,d\), let
\(
    M_t^j
    :=
    \sum_{\ell=1}^q
    \int_0^t
    \sigma_{j\ell}(u,\bm X_u)\,\mathrm dB_u^\ell .
\)
Then
\begin{equation}
    X_T^\ast
    \le
    \|\bm x_0\|_\infty+RT
    +
    \max_{1\le j\le d}\sup_{0\le t\le T}|M_t^j|,
    \qquad
    \langle M^j\rangle_T\le R^2T.
\end{equation}
The exponential martingale inequality and a union bound therefore imply
that there exist constants \(C,c,C_0>0\) such that, for every \(r,s\ge0\),
\begin{equation}
\label{eq:analytic-diffusion-tail-mgf}
    \mathbb P(X_T^\ast>r)
    \le
    C\exp\{-c(r-C_0)_+^2\},
    \qquad
    \mathbb E\exp(sX_T^\ast)
    \le
    C\exp\{C(s+s^2)\}.
\end{equation}
The second bound follows by integrating the first tail bound.
Since \(m_K=O(K^\beta)\) and \(m_K-C_0\ge cK^\beta\) for all
sufficiently large \(K\), enlarging the prefactor to cover the remaining
values gives
\begin{equation}
\label{eq:analytic-AK-tail}
    \mathbb P(A_K^c)
    =
    \mathbb P(X_T^\ast>m_K)
    \le
    C\exp(-cK^{2\beta}).
\end{equation}

For \(a=0,\ldots,d\), define
\(\widetilde h_{a,K}(u,y):=h_a(Tu,2m_Ky)\) on
\(D_0:=[0,1]\times[-1,1]^d\), using the same notation for its
holomorphic extension.
Since \(\bar\rho<\rho_0\),
\(\mathcal E_{\bar\rho}^{[0,1]}\subset V_{\rho_0}\).
Moreover, the maximum absolute imaginary part on
\(\mathcal E_{\bar\rho}\) is
\((\bar\rho-\bar\rho^{-1})/2\), and hence
\(\|\operatorname{Im}(2m_Kz)\|_1
\le d(\bar\rho-\bar\rho^{-1})m_K\) for
\(z\in\mathcal E_{\bar\rho}^d\).
Therefore, \autoref{eq:scale-analytic-class} gives
\begin{equation}
\label{eq:analytic-rescaled-envelope}
    \max_{0\le a\le d}
    \sup_{\substack{z_0\in\mathcal E_{\bar\rho}^{[0,1]}\\
                    z\in\mathcal E_{\bar\rho}^d}}
    |\widetilde h_{a,K}(z_0,z)|
    \le
    R\exp\{d\tau(\bar\rho-\bar\rho^{-1})m_K\}.
\end{equation}
For \(\nu=(\nu_0,\ldots,\nu_d)\in\mathbb N_0^{d+1}\), let
\(T_\nu(u,y):=T_{\nu_0}(2u-1)\prod_{r=1}^dT_{\nu_r}(y_r)\).
Writing \(c_{a,K,\nu}\) for the corresponding tensor-product
Chebyshev coefficient, the product-polyellipse coefficient bound
\citep[Corollary~3.2]{WangZhang2020ACM} and
\autoref{eq:analytic-rescaled-envelope} give
\begin{equation}
\label{eq:analytic-chebyshev-coefficients}
\begin{aligned}
    \widetilde h_{a,K}(u,y)
    &=
    \sum_{\nu\in\mathbb N_0^{d+1}}
    c_{a,K,\nu}T_\nu(u,y),
    \\
    \widetilde p_{a,K}(u,y)
    &:=
    \sum_{|\nu|_1\le K-1}
    c_{a,K,\nu}T_\nu(u,y),
    \\
    |c_{a,K,\nu}|
    &\le
    C R
    \exp\{d\tau(\bar\rho-\bar\rho^{-1})m_K\}
    \bar\rho^{-|\nu|_1}.
\end{aligned}
\end{equation}
The polynomial \(\widetilde p_{a,K}\) has total degree at most
\(K-1\), rather than tensor degree at most \(K-1\).
Since \(|T_j|\le1\) on \([-1,1]\) and the number of multi-indices
satisfying \(|\nu|_1=j\) is \(\binom{j+d}{d}\), summing
\autoref{eq:analytic-chebyshev-coefficients} over \(j\ge K\)
bounds the value tail by
\(CRK^d\exp\{d\tau(\bar\rho-\bar\rho^{-1})m_K\}\bar\rho^{-K}\).
Moreover, \(\|T_{\nu_r}'\|_{L^\infty([-1,1])}=\nu_r^2\), so the
same coefficient bound shows that each differentiated series
converges uniformly absolutely on \(D_0\).
Termwise differentiation is therefore valid and gives the same
tail bound with \(K^d\) replaced by \(K^{d+2}\) for each state
derivative.
Consequently,
\begin{equation}
\label{eq:analytic-rescaled-approx}
\begin{aligned}
    &\max_{0\le a\le d}
    \|\widetilde h_{a,K}-\widetilde p_{a,K}\|_{L^\infty(D_0)}
    +
    \max_{1\le a\le d}
    \|\nabla_y\widetilde h_{a,K}
      -\nabla_y\widetilde p_{a,K}\|_{L^\infty(D_0)}
    \\
    &\qquad\le
    C R K^{d+2}
    \exp\{d\tau(\bar\rho-\bar\rho^{-1})m_K\}
    \bar\rho^{-K}
    \le
    C_{\rho}R\rho^{-K}.
\end{aligned}
\end{equation}
The last inequality follows from \(m_K=O(K^\beta)\),
\(\beta<1\), and \(\rho<\bar\rho\), since the polynomial and
subexponential factors are absorbed by
\((\bar\rho/\rho)^K\).

Returning to the original variables, define
\begin{equation}
\label{eq:analytic-polynomial-original-scale}
    p_{a,K}(t,x)
    :=
    \widetilde p_{a,K}
    \left(
        \frac{t}{T},
        \frac{x}{2m_K}
    \right),
    \qquad
    a=0,\ldots,d.
\end{equation}
The map \((t,x)\mapsto(t/T,x/(2m_K))\) sends \(D_K\) onto
\(D_0\) and does not increase total degree, so each \(p_{a,K}\)
has total degree at most \(K-1\).
Moreover, for \((t,x)\in D_K\), the chain rule gives
\(\nabla_x(h_a-p_{a,K})(t,x)
=(2m_K)^{-1}
(\nabla_y\widetilde h_{a,K}-\nabla_y\widetilde p_{a,K})
(t/T,x/(2m_K))\).
Since \(m_K\ge1\), \autoref{eq:analytic-rescaled-approx} therefore
yields
\begin{equation}
\label{eq:analytic-local-approx-original}
    \max_{0\le a\le d}
    \|h_a-p_{a,K}\|_{L^\infty(D_K)}
    +
    \max_{1\le a\le d}
    \|\nabla_xh_a-\nabla_xp_{a,K}\|_{L^\infty(D_K)}
    \le
    C_{\rho}R\rho^{-K}.
\end{equation}

Let \(\delta_{a,K}:=h_a-p_{a,K}\) and define
\begin{equation}
\label{eq:analytic-YK}
    Y_K
    :=
    \sum_{a=0}^d
    \int_0^T
    p_{a,K}(t,\bm X_t)
    \circ\mathrm d\widehat X_t^a.
\end{equation}
Let \(e_0,\ldots,e_d\) denote the canonical basis of
\(\mathbb R^{d+1}\), and fix \(a\in\{0,\ldots,d\}\).
Since \(p_{a,K}\) has total degree at most \(K-1\) and
\(\bm X_0=\bm x_0\) is deterministic, its evaluation at
\((t,\bm X_t)\) is a polynomial of degree at most \(K-1\) in the
first-level prefix signature coordinates of \(\widehat{\bm X}\).
Indeed,
\begin{equation}
\label{eq:analytic-first-level-prefix}
    t
    =
    \left\langle
        e_0,S(\widehat{\bm X})_{0,t}
    \right\rangle,
    \qquad
    X_t^b
    =
    x_0^b
    +
    \left\langle
        e_b,S(\widehat{\bm X})_{0,t}
    \right\rangle,
    \quad b=1,\ldots,d.
\end{equation}
The deterministic constants \(x_0^b\) are represented by the
level-zero signature coordinate. By the shuffle identity, writing
\(T^{(K-1)}(\mathbb R^{d+1})
:=
\bigoplus_{r=0}^{K-1}(\mathbb R^{d+1})^{\otimes r}\), there exists
a deterministic tensor
\(q_{a,K}\in T^{(K-1)}(\mathbb R^{d+1})\) such that the recursive
definition of the Stratonovich signature gives
\begin{equation}
\label{eq:analytic-prefix-append}
\begin{aligned}
    &p_{a,K}(t,\bm X_t)
    =
    \left\langle
        q_{a,K},
        S(\widehat{\bm X})_{0,t}
    \right\rangle,
    \\
    \mathrm{and} \qquad &\int_0^T
    p_{a,K}(t,\bm X_t)
    \circ\mathrm d\widehat X_t^a
    =
    \left\langle
        q_{a,K}\otimes e_a,
        S(\widehat{\bm X})_{0,T}
    \right\rangle.
\end{aligned}
\end{equation}
Every homogeneous component of \(q_{a,K}\otimes e_a\) has level
at most \(K\). The same identity covers \(a=0\), for which
\(\circ\mathrm d\widehat X_t^0=\mathrm dt\) and the appended letter
is the time letter \(e_0\). Summing over \(a=0,\ldots,d\) therefore
shows that \(Y_K\in\mathbb S_K(\bm X)\), and hence
\begin{equation}
\label{eq:analytic-projection-bound}
    \inf_{Z\in\mathbb S_K(\bm X)}
    \mathbb E\left|
        F_{\bm h}(\bm X)-Z
    \right|^2
    \le
    \mathbb E\left|
        F_{\bm h}(\bm X)-Y_K
    \right|^2.
\end{equation}

We first control the last error on \(A_K\).
Since \(\bm X\) has continuous paths, \(\eta_K\) is a stopping time,
\((t,\bm X_t)\in D_K\) for every \(t\le\eta_K\), including the exit
point, and \(A_K\subseteq\{\eta_K=T\}\). Therefore,
\begin{equation}
\label{eq:analytic-stopped-error-reduction}
\begin{aligned}
    &\mathbb E\left[
        |F_{\bm h}(\bm X)-Y_K|^2\mathbf 1_{A_K}
    \right]
    \\
    &\quad\le
    \mathbb E\left|
        \int_0^{\eta_K}
        \delta_{0,K}(t,\bm X_t)\,\mathrm dt
        +
        \sum_{a=1}^d
        \int_0^{\eta_K}
        \delta_{a,K}(t,\bm X_t)
        \circ\mathrm dX_t^a
    \right|^2.
\end{aligned}
\end{equation}
Joint holomorphy of \(h_a\) and polynomiality of \(p_{a,K}\) ensure
that \(\delta_{a,K}\) is sufficiently smooth for the
Stratonovich--It\^{o} conversion. For \(a=1,\ldots,d\), it gives
\begin{equation}
\label{eq:analytic-stopped-stratonovich}
\begin{aligned}
    \int_0^{\eta_K}
    \delta_{a,K}(t,\bm X_t)\circ\mathrm dX_t^a
    ={}&
    \int_0^{\eta_K}
    \delta_{a,K}(t,\bm X_t)
    b_a(t,\bm X_t)\,\mathrm dt
    \\
    &+
    \sum_{\ell=1}^q
    \int_0^{\eta_K}
    \delta_{a,K}(t,\bm X_t)
    \sigma_{a\ell}(t,\bm X_t)\,\mathrm dB_t^\ell
    \\
    &+
    \frac12
    \int_0^{\eta_K}
    \sum_{r=1}^d
    \partial_{x_r}\delta_{a,K}(t,\bm X_t)
    (\bm\sigma\bm\sigma^\top)_{ra}(t,\bm X_t)\,\mathrm dt.
\end{aligned}
\end{equation}
Only spatial derivatives enter the quadratic-covariation correction,
so no approximation bound for the time derivative is required.
Set
\(\varepsilon_K
:=
\max_{0\le a\le d}\|\delta_{a,K}\|_{L^\infty(D_K)}
+
\max_{1\le a\le d}
\|\nabla_x\delta_{a,K}\|_{L^\infty(D_K)}\).
Then \autoref{eq:analytic-local-approx-original} gives
\(\varepsilon_K\le C_\rho R\rho^{-K}\).
Cauchy--Schwarz for the finite-variation terms and It\^{o} isometry
for the stopped martingale term yield
\begin{equation}
\label{eq:analytic-local-error}
\begin{aligned}
    \mathbb E\left[
        |F_{\bm h}(\bm X)-Y_K|^2\mathbf 1_{A_K}
    \right]
    &\le
    C(1+R^2+R^4)\varepsilon_K^2
    \le
    C_\rho R^2\rho^{-2K}.
\end{aligned}
\end{equation}
Here the fixed factor \(1+R^2+R^4\) is absorbed into \(C_\rho\),
which is allowed to depend on the class bound \(R\).
The stopping argument is needed because \(A_K\) depends on the entire
path and therefore cannot be inserted directly into an It\^{o}
isometry.

It remains to control the complement of \(A_K\).
On the real state space, \autoref{eq:scale-analytic-class} gives
\(\max_{0\le a\le d}|h_a(t,x)|\le R\).
Moreover, fixing a state coordinate \(r\) and applying the Cauchy
integral formula on the unit circle in that coordinate gives
\(\max_{0\le a\le d}|\partial_{x_r}h_a(t,x)|
\le e^\tau R\), uniformly over \((t,x)\in D\).
The Stratonovich--It\^{o} decomposition in
\autoref{eq:analytic-stopped-stratonovich}, now applied on
\([0,T]\) to \(h_a\), together with boundedness of
\(\bm b\) and \(\bm\sigma\) and the Burkholder--Davis--Gundy
inequality, therefore yields
\begin{equation}
\label{eq:analytic-functional-fourth}
\begin{aligned}
    \left\{
        \mathbb E
        |F_{\bm h}(\bm X)|^4
    \right\}^{1/2}
    &\le
    C\left(R+R^2+R^3\right)^2
    \le
    C R^2.
\end{aligned}
\end{equation}
In the last inequality, the fixed factor
\((1+R+R^2)^2\) is absorbed into \(C\), which is allowed to
depend on the class bound \(R\).
Hence, by Cauchy--Schwarz and \autoref{eq:analytic-AK-tail},
\begin{equation}
\label{eq:analytic-functional-tail}
\begin{aligned}
    \mathbb E\left[
        |F_{\bm h}(\bm X)|^2\mathbf 1_{A_K^c}
    \right]
    &\le
    \left\{
        \mathbb E|F_{\bm h}(\bm X)|^4
    \right\}^{1/2}
    \mathbb P(A_K^c)^{1/2}
    \\
    &\le
    C R^2\exp(-cK^{2\beta})
    =
    o\left(R^2\rho^{-2K}\right).
\end{aligned}
\end{equation}
The last equality follows because \(\beta=3/4\), and hence
\(K^{2\beta}=K^{3/2}\) dominates \(K\).

To control \(Y_K\), first note from
\autoref{eq:analytic-local-approx-original} and the real-axis bound on
\(h_a\) that
\(\max_{0\le a\le d}\|p_{a,K}\|_{L^\infty(D_K)}\le CR\).
After transforming the time coordinate by \(z_0=2u-1\), the
tensor-product Chebyshev coefficient formula therefore gives
\(|c_{a,K,\nu}|
\le2^{d+1}\|\widetilde p_{a,K}\|_{L^\infty(D_0)}
\le CR\) whenever \(|\nu|_1\le K-1\).
Moreover, for every real \(y\),
\(|T_j(y)|\le(1+2|y|)^j\) and
\(|T_j'(y)|\le j^2(1+2|y|)^j\).
Since
\(\#\{\nu\in\mathbb N_0^{d+1}:|\nu|_1\le K-1\}
=\binom{K+d}{d+1}\le CK^{d+1}\), while the chain rule contributes
the factor \((2m_K)^{-1}\nu_r^2\le K^2\) to each spatial derivative,
summing over the retained multi-indices yields, uniformly over
\(t\in[0,T]\) and \(x\in\mathbb R^d\),
\begin{equation}
\label{eq:analytic-polynomial-growth}
    \max_{0\le a\le d}|p_{a,K}(t,x)|
    +
    \max_{1\le a\le d}\|\nabla_xp_{a,K}(t,x)\|
    \le
    CRK^{d+3}
    \left(
        1+\frac{\|x\|_\infty}{m_K}
    \right)^K.
\end{equation}

Applying the Stratonovich--It\^{o} conversion to
\autoref{eq:analytic-YK} gives
\begin{equation}
\label{eq:analytic-YK-ito}
\begin{aligned}
    Y_K
    ={}&
    \int_0^T p_{0,K}(t,\bm X_t)\,\mathrm dt
    +
    \sum_{a=1}^d
    \int_0^T
    p_{a,K}(t,\bm X_t)b_a(t,\bm X_t)\,\mathrm dt
    \\
    &+
    \sum_{\ell=1}^q
    \int_0^T
    \left\{
        \sum_{a=1}^d
        p_{a,K}(t,\bm X_t)
        \sigma_{a\ell}(t,\bm X_t)
    \right\}\mathrm dB_t^\ell
    \\
    &+
    \frac12
    \sum_{a,r=1}^d
    \int_0^T
    \partial_{x_r}p_{a,K}(t,\bm X_t)
    (\bm\sigma\bm\sigma^\top)_{ra}(t,\bm X_t)
    \,\mathrm dt.
\end{aligned}
\end{equation}
Consequently, \autoref{eq:analytic-polynomial-growth},
Cauchy--Schwarz for the finite-variation terms, and the
Burkholder--Davis--Gundy inequality imply
\begin{equation}
\label{eq:analytic-YK-fourth}
\begin{aligned}
    &\left\{\mathbb E|Y_K|^4\right\}^{1/2}
    \le
    CR^2(1+R+R^2)^2K^{2d+6}
    \left\{
        \mathbb E
        \left(
            1+\frac{X_T^\ast}{m_K}
        \right)^{4K}
    \right\}^{1/2},
    \\
    &\left\{
        \mathbb E
        \left(
            1+\frac{X_T^\ast}{m_K}
        \right)^{4K}
    \right\}^{1/2}
    \le
    \left\{
        \mathbb E
        \exp\left(
            \frac{4KX_T^\ast}{m_K}
        \right)
    \right\}^{1/2}
    \\
    &\qquad \qquad \qquad \qquad \qquad \,\le
    C\exp\left\{
        C\left(
            \frac{K}{m_K}
            +
            \frac{K^2}{m_K^2}
        \right)
    \right\}
    \le
    C\exp\left(CK^{2-2\beta}\right).
\end{aligned}
\end{equation}
Since the fixed factor \((1+R+R^2)^2\) may be absorbed into the
class constant, it follows that
\begin{equation}
\label{eq:analytic-YK-fourth-final}
    \left\{\mathbb E|Y_K|^4\right\}^{1/2}
    \le
    CR^2K^{2d+6}\exp\left(CK^{2-2\beta}\right).
\end{equation}

Combining \autoref{eq:analytic-YK-fourth-final} with
\autoref{eq:analytic-AK-tail} and applying Cauchy--Schwarz gives
\begin{equation}
\label{eq:analytic-YK-tail}
\begin{aligned}
    \mathbb E\left[
        |Y_K|^2\mathbf 1_{A_K^c}
    \right]
    &\le
    \left\{
        \mathbb E|Y_K|^4
    \right\}^{1/2}
    \mathbb P(A_K^c)^{1/2}
    \\
    &\le
    C R^2K^{2d+6}
    \exp\left\{
        CK^{2-2\beta}-cK^{2\beta}
    \right\}
    \\
    &=
    C R^2K^{2d+6}
    \exp\left\{
        CK^{1/2}-cK^{3/2}
    \right\}
    =
    o\left(R^2\rho^{-2K}\right).
\end{aligned}
\end{equation}
Indeed, after division by \(\rho^{-2K}\), the exponent acquires only
the additional linear term \(2K\log\rho\), which is dominated by
\(-cK^{3/2}\). Finally, using
\(|u-v|^2\le2|u|^2+2|v|^2\), together with
\autoref{eq:analytic-functional-tail} and
\autoref{eq:analytic-YK-tail}, yields
\begin{equation}
\label{eq:analytic-tail-error}
    \mathbb E\left[
        |F_{\bm h}(\bm X)-Y_K|^2
        \mathbf 1_{A_K^c}
    \right]
    \le
    2\mathbb E\left[
        |F_{\bm h}(\bm X)|^2
        \mathbf 1_{A_K^c}
    \right]
    +
    2\mathbb E\left[
        |Y_K|^2
        \mathbf 1_{A_K^c}
    \right]
    =
    o\left(R^2\rho^{-2K}\right).
\end{equation}

Combining \autoref{eq:analytic-projection-bound},
\autoref{eq:analytic-local-error}, and
\autoref{eq:analytic-tail-error}, and enlarging \(C_\rho\) if
necessary, there exists an integer \(K_0\ge2\), depending only on the
fixed class parameters, such that, uniformly over
\((\bm X,\bm h)\in\mathcal U_{\rho_0,\tau,R}^{\mathrm{an}}\),
\begin{equation}
\label{eq:analytic-final-large-K}
\begin{aligned}
    \inf_{Z\in\mathbb S_K(\bm X)}
    \mathbb E\left|
        F_{\bm h}(\bm X)-Z
    \right|^2
    &\le
    \mathbb E\left|
        F_{\bm h}(\bm X)-Y_K
    \right|^2
    \\
    &=
    \mathbb E\left[
        |F_{\bm h}(\bm X)-Y_K|^2\mathbf 1_{A_K}
    \right]
    +
    \mathbb E\left[
        |F_{\bm h}(\bm X)-Y_K|^2\mathbf 1_{A_K^c}
    \right]
    \\
    &\le
    C_\rho R^2\rho^{-2K}
\end{aligned}
\end{equation}
for every \(K\ge K_0\).

For the finitely many integers \(2\le K<K_0\), since
\(0\in\mathbb S_K(\bm X)\), \autoref{eq:analytic-functional-fourth}
gives
\begin{equation}
\label{eq:analytic-final-small-K}
\begin{aligned}
    \inf_{Z\in\mathbb S_K(\bm X)}
    \mathbb E\left|
        F_{\bm h}(\bm X)-Z
    \right|^2
    &\le
    \mathbb E|F_{\bm h}(\bm X)|^2\le
    \left\{
        \mathbb E|F_{\bm h}(\bm X)|^4
    \right\}^{1/2}
    \\
    &\le
    C R^2
    \le
    C R^2\rho^{2(K_0-1)}\rho^{-2K}.
\end{aligned}
\end{equation}
Increasing \(C_\rho\) once more therefore proves
\autoref{eq:analytic-signature-upper-bound} for every \(K\ge2\).
\QED

\section{Theoretical Details for \autoref{sec_Reg}}
\label{Appen_ProofsSec4}

This appendix proves the statistical results in \autoref{sec_Reg}.
Throughout the appendix, we maintain the standing assumption that the underlying path is an It\^{o} diffusion satisfying the conditions in \autoref{sec:signature-approximation}, and that the target variable \(Y\) is generated from the pair \((\bm X,\bm h)\in\mathcal U_{\gamma,R}\).

\subsection{Technical Lemmas}

We first prove several technical lemmas.
The next two lemmas establish moment and sub-Weibull tail bounds for individual signature coordinates. 
We first introduce the \(\psi_\kappa\)-functional used to state the tail result. 
Since the exponent appearing below may be smaller than one, we refer to it as a sub-Weibull functional rather than an Orlicz norm. 
For \(\kappa>0\), define
\begin{equation}
\label{eq:sec4_orlicz}
\|Y\|_{\psi_\kappa}
:=
\inf\left\{
A>0:
\mathbb E\exp\left\{
\left(\frac{|Y|}{A}\right)^\kappa
\right\}
\le 2
\right\}.
\end{equation}
A random variable \(Y\) is called sub-Weibull with tail parameter \(\kappa\) if \(\|Y\|_{\psi_\kappa}<\infty\); see \citep[Section~6.1]{ZhangandChen2020}. Whenever \(0<\|Y\|_{\psi_\kappa}<\infty\), Markov's inequality and \autoref{eq:sec4_orlicz} give
\begin{equation}
\label{eq:sec4_orlicz_tail}
\mathbb P(|Y|>x)
\le
2\exp\left\{
-\left(\frac{x}{\|Y\|_{\psi_\kappa}}\right)^\kappa
\right\},
\qquad x>0.
\end{equation}

Before stating the lemma, we recall an \(L^p\) form of the
Burkholder--Davis--Gundy inequality for continuous martingales; see
\citet{Ren2008BDG}. For \(p\ge2\), write
\(\|Y\|_p=(\mathbb E|Y|^p)^{1/p}\). There exists a universal constant
\(C_{\mathrm{BDG}}>0\) such that, for every \(p\ge2\), \(t\in[0,T]\),
and every \(\mathbb R^q\)-valued predictable process \(\bm H\) for
which the right-hand side is finite,
\begin{equation}
\label{eq:sec4_bdg}
\left\|
\int_0^t \bm H_u^\top\,\mathrm d\bm B_u
\right\|_p
\le
C_{\mathrm{BDG}}\sqrt p
\left\|
\left(
\int_0^t \|\bm H_u\|_2^2\,\mathrm du
\right)^{1/2}
\right\|_p.
\end{equation}

We first establish a weighted moment bound directly from the bounded-coefficient diffusion structure.

\begin{lem}
\label{lem:sigmoment}
Let \(\bm X\) satisfy \autoref{eq:smooth-coefficient-X-SDE}, and suppose that
$\|\bm b(t,\bm x)\|
+
\|\bm\sigma(t,\bm x)\|_F
\le
R$ uniformly in \((t,\bm x)\). 
Define $
C_{\mathrm{sig}}
:=
4\left\{
1
+
(1+R)\sqrt T
+
C_{\mathrm{BDG}}R
+
R
\right\}$.
Then, for every \(p\ge2\), \(k\ge1\),
\(I\in\{0,\ldots,d\}^k\), and \(t\in[0,T]\), we have
\begin{equation}
\left\|
S(\widehat{\bm X})_{0,t}^I
\right\|_p
\le
\frac{
(C_{\mathrm{sig}}\sqrt p)^k t^{k/2}
}{
\sqrt{k!}
},\qquad k=|I|.
\label{eq:signature-factorial-moment}
\end{equation}
\end{lem}

\prf
For a word \(J\), write
\(S_t^J:=S(\widehat{\bm X})_{0,t}^J\), with
\(S_t^\emptyset\equiv1\). For \(p\ge2\), define
\begin{equation}
F_{0,p}(t):=1,
\qquad
F_{k,p}(t)
:=
\max_{J\in\{0,\ldots,d\}^k}
\|S_t^J\|_p,
\qquad k\ge1.
\end{equation}
Let \(I=(i_1,\ldots,i_k)\) and
\(I'=(i_1,\ldots,i_{k-1})\). The recursive definition of the
Stratonovich signature gives $
S_t^I
=
\int_0^t
S_u^{I'}
\circ\mathrm d\widehat X_u^{i_k}$.
If \(i_k=0\), then \(S_t^I=\int_0^t S_u^{I'}\,\mathrm du\). If
\(i_k=a\in\{1,\ldots,d\}\), the Stratonovich--It\^{o} conversion gives
\begin{equation}
\begin{aligned}
S_t^I
=
\int_0^t
S_u^{I'}b_a(u,\bm X_u)\,\mathrm du
+
\sum_{\ell=1}^q
\int_0^t
S_u^{I'}\sigma_{a\ell}(u,\bm X_u)\,\mathrm dB_u^\ell +
\frac12
\left[S^{I'},X^a\right]_t.
\end{aligned}
\end{equation}
When \(k\ge2\) and \(I'=(I'',r)\) for some
\(r\in\{1,\ldots,d\}\), the local martingale part of \(S^{I'}\) is
generated by \(S_u^{I''}\bm\sigma_{r\cdot}(u,\bm X_u)\). Therefore,
\begin{equation}
\mathrm d
\left[S^{I'},X^a\right]_u
=
S_u^{I''}
\sum_{\ell=1}^q
\sigma_{r\ell}(u,\bm X_u)
\sigma_{a\ell}(u,\bm X_u)
\,\mathrm du.
\end{equation}
The quadratic covariation is zero when \(I'=\emptyset\) or the last
letter of \(I'\) is the time channel.

The coefficient bound implies
\begin{equation}
|b_a(t,\bm x)|\le R,
\qquad
\|\bm\sigma_{a\cdot}(t,\bm x)\|_2\le R,
\qquad
\left|
\sum_{\ell=1}^q
\sigma_{r\ell}(t,\bm x)\sigma_{a\ell}(t,\bm x)
\right|
\le R^2.
\end{equation}
Moreover, since \(p/2\ge1\), Minkowski's inequality gives, for every
progressively measurable scalar process \(Z\),
\begin{equation}
\begin{aligned}
\left\|
\left(
\int_0^t |Z_u|^2\,\mathrm du
\right)^{1/2}
\right\|_p^2
&=
\left\|
\int_0^t |Z_u|^2\,\mathrm du
\right\|_{p/2} \\
&\le
\int_0^t
\bigl\||Z_u|^2\bigr\|_{p/2}\,\mathrm du
=
\int_0^t
\|Z_u\|_p^2\,\mathrm du.
\end{aligned}
\end{equation}
Consequently, Minkowski's inequality and
\autoref{eq:sec4_bdg}, together with an enlargement of
\(\max\{1,R\}\) to \(1+R\), imply
\begin{equation}
\begin{aligned}
F_{k,p}(t)
\le{}&
(1+R)
\int_0^t
F_{k-1,p}(u)\,\mathrm du
+
C_{\mathrm{BDG}}R\sqrt p
\left\{
\int_0^t
F_{k-1,p}(u)^2\,\mathrm du
\right\}^{1/2} \\
&+
\frac{R^2}{2}
\int_0^t
F_{k-2,p}(u)\,\mathrm du,
\end{aligned}
\label{eq:signature-moment-recursion}
\end{equation}
where the final term is absent when \(k=1\).
Then, set
\begin{equation}
G_{0,p}(t):=1,
\qquad
G_{k,p}(t)
:=
\frac{
(C_{\mathrm{sig}}\sqrt p)^k t^{k/2}
}{
\sqrt{k!}
},
\qquad k\ge1.
\end{equation}
The assertion is immediate at \(t=0\). Suppose inductively that
\(F_{j,p}(u)\le G_{j,p}(u)\) for every \(j<k\) and \(u\in[0,T]\).
For \(t>0\), direct calculation gives
\begin{equation}
\begin{aligned}
&\int_0^t G_{k-1,p}(u)\,\mathrm du
=
G_{k,p}(t)
\frac{2\sqrt{kt}}
{(k+1)C_{\mathrm{sig}}\sqrt p}, \\
&\left\{
\int_0^t G_{k-1,p}(u)^2\,\mathrm du
\right\}^{1/2}
=
\frac{G_{k,p}(t)}
{C_{\mathrm{sig}}\sqrt p}, \\
&\int_0^t G_{k-2,p}(u)\,\mathrm du
=
G_{k,p}(t)
\frac{2}{C_{\mathrm{sig}}^2p}
\sqrt{\frac{k-1}{k}},
\qquad k\ge2.
\end{aligned}
\end{equation}
Substitution into \autoref{eq:signature-moment-recursion} yields
\begin{equation}
F_{k,p}(t)
\le
G_{k,p}(t)
\left\{
\frac{2(1+R)\sqrt{kt}}
{(k+1)C_{\mathrm{sig}}\sqrt p}
+
\frac{C_{\mathrm{BDG}}R}
{C_{\mathrm{sig}}}
+
\frac{R^2}
{C_{\mathrm{sig}}^2p}
\sqrt{\frac{k-1}{k}}
\right\},
\end{equation}
where the final term is zero when \(k=1\). Since
\(2\sqrt{k}/(k+1)\le1\), \(p\ge2\), and \(t\le T\), the expression
in braces is bounded above by
\begin{equation}
\frac{(1+R)\sqrt T+C_{\mathrm{BDG}}R}
{C_{\mathrm{sig}}}
+
\frac{R^2}{C_{\mathrm{sig}}^2}
\le
\frac14+\frac1{16}
<1.
\end{equation}
Induction therefore gives \(F_{k,p}(t)\le G_{k,p}(t)\) for
every \(k\ge1\), which proves
\autoref{eq:signature-factorial-moment}.
\QED

The factorial moment bound  yields a sub-Weibull bound for each signature coordinate.

\begin{lem}
\label{lem:sigtail}
Under the conditions of \autoref{lem:sigmoment}, define $M_\psi
:=
4\sqrt{e}\,
\max\left\{
1,
C_{\mathrm{sig}}\sqrt T
\right\}$.
Then, for every \(k\ge1\) and
\(I\in\{0,\ldots,d\}^k\),
\begin{equation}
\left\|
S(\widehat{\bm X})_{0,T}^I
-
\mathbb E S(\widehat{\bm X})_{0,T}^I
\right\|_{\psi_{2/k}}
\le
M_\psi^k.
\label{eq:signature-coordinate-psi}
\end{equation}
Consequently, for every \(x>0\),
\begin{equation}
\mathbb P\left(
\left|
S(\widehat{\bm X})_{0,T}^I
-
\mathbb E S(\widehat{\bm X})_{0,T}^I
\right|
>x
\right)
\le
2\exp\left\{
-
\left(
\frac{x}{M_\psi^k}
\right)^{2/k}
\right\}.
\label{eq:signature-coordinate-tail}
\end{equation}
Thus, every centered signature coordinate of word length \(k\) is
sub-Weibull of order \(2/k\).
\end{lem}

\prf
Fix \(k\ge1\) and \(I\in\{0,\ldots,d\}^k\), and write
\begin{equation}
Z_I
:=
S(\widehat{\bm X})_{0,T}^I
-
\mathbb E S(\widehat{\bm X})_{0,T}^I,
\qquad
A
:=
\max\left\{
1,
C_{\mathrm{sig}}\sqrt T
\right\}.
\end{equation}
For every integer \(m\ge1\), the triangle inequality, Jensen's
inequality, and \autoref{eq:signature-factorial-moment} give
\begin{equation}
\begin{aligned}
\|Z_I\|_{2m}
&\le
2
\left\|
S(\widehat{\bm X})_{0,T}^I
\right\|_{2m} \\
&\le
\frac{
2(C_{\mathrm{sig}}\sqrt{2m})^kT^{k/2}
}{
\sqrt{k!}
}
\le
\frac{
2A^k(2m)^{k/2}
}{
\sqrt{k!}
}.
\end{aligned}
\label{eq:centered-signature-even-moment}
\end{equation}
Since \(2m/k\le2m\), Lyapunov's inequality yields
\begin{equation}
\mathbb E|Z_I|^{2m/k}
\le
\left(
\mathbb E|Z_I|^{2m}
\right)^{1/k}
=
\|Z_I\|_{2m}^{2m/k}.
\label{eq:centered-signature-lyapunov}
\end{equation}
Using \(M_\psi=4\sqrt{e}\,A\), monotone convergence and
\autoref{eq:centered-signature-even-moment}--\autoref{eq:centered-signature-lyapunov}
give
\begin{equation}
\mathbb E
\exp\left\{
\left(
\frac{|Z_I|}{M_\psi^k}
\right)^{2/k}
\right\}=
1+
\sum_{m=1}^{\infty}
\frac{
\mathbb E|Z_I|^{2m/k}
}{
m!M_\psi^{2m}
}, 
\end{equation}
and
\begin{equation}
\frac{
\mathbb E|Z_I|^{2m/k}
}{
m!M_\psi^{2m}
}
\le
\frac{
2^{2m/k}(2m)^m
}{
m!(16e)^m(k!)^{m/k}
} 
\le
\frac{1}{m!}
\left(
\frac{m}{2e}
\right)^m
\le
2^{-m}.
\label{eq:signature-exponential-moment-bound}
\end{equation}
Here, the last two inequalities use
\(2^{2m/k}\le4^m\), \(k!\ge1\), and
\(m!\ge(m/e)^m\). It follows that
\begin{equation}
\mathbb E
\exp\left\{
\left(
\frac{|Z_I|}{M_\psi^k}
\right)^{2/k}
\right\}
\le
1+\sum_{m=1}^{\infty}2^{-m}
=
2.
\end{equation}
Therefore, the definition in \autoref{eq:sec4_orlicz} gives
\autoref{eq:signature-coordinate-psi}. Finally, Markov's inequality
implies that, for every \(x>0\),
\begin{equation}
\begin{aligned}
\mathbb P(|Z_I|>x)
&\le
\exp\left\{
-
\left(
\frac{x}{M_\psi^k}
\right)^{2/k}
\right\}
\mathbb E
\exp\left\{
\left(
\frac{|Z_I|}{M_\psi^k}
\right)^{2/k}
\right\} \\
&\le
2\exp\left\{
-
\left(
\frac{x}{M_\psi^k}
\right)^{2/k}
\right\},
\end{aligned}
\end{equation}
which proves \autoref{eq:signature-coordinate-tail}.
\QED

\subsection{Proof of \autoref{thm: OLSCLT}}

\prf
(i) We first establish relative concentration of the sample Gram matrix. Write
\(\bm s_i:=\bm s_K(\widehat{\bm X}_i)\), let
\(\bm u_i:=\bm\Sigma_K^{-1/2}\bm s_i\), and let \(\mathbf U\) be the matrix with rows \(\bm u_i^\top\). Then
\begin{equation}
\widehat{\bm G}_K
:=
\frac{1}{n}\mathbf U^\top\mathbf U
=
\bm\Sigma_K^{-1/2}
\widehat{\bm\Sigma}_K
\bm\Sigma_K^{-1/2},
\qquad
\widehat{\bm\Sigma}_K
:=
\frac{1}{n}\mathbf S^\top\mathbf S,
\label{eq:OLS-standardized-Gram}
\end{equation}
and \(\mathbb E[\bm u_i\bm u_i^\top]=\bm I_{d_K}\). Independence gives
\begin{equation}
\begin{aligned}
\mathbb E\left[
\|\widehat{\bm G}_K-\bm I_{d_K}\|_F^2
\right]
&=
\frac{1}{n}
\mathbb E\left[
\|\bm u_i\bm u_i^\top-\bm I_{d_K}\|_F^2
\right]\\
&=
\frac{1}{n}
\left\{
\mathbb E\|\bm u_i\|^4-d_K
\right\}
=
o(1).
\end{aligned}
\label{eq:OLS-relative-Gram}
\end{equation}
Hence,
\(\|\widehat{\bm G}_K-\bm I_{d_K}\|_{\mathrm{op}}=o_P(1)\). In particular, \(\widehat{\bm G}_K\), and therefore \(\widehat{\bm\Sigma}_K\), is positive definite with probability approaching one. We work on this event and define inverse-based quantities arbitrarily on its complement.

Let \(\bm u_K(\widehat{\bm x}):=\bm\Sigma_K^{-1/2}\bm s_K(\widehat{\bm x})\), write \(v_K:=v_K(\widehat{\bm x})\), and set
\begin{equation}
\widehat v_K^2
:=
\bm s_K(\widehat{\bm x})^\top
\widehat{\bm\Sigma}_K^{-1}
\bm s_K(\widehat{\bm x})
=
\bm u_K(\widehat{\bm x})^\top
\widehat{\bm G}_K^{-1}
\bm u_K(\widehat{\bm x}).
\label{eq:OLS-sample-leverage}
\end{equation}
By \autoref{eq:OLS-relative-Gram} and the inverse perturbation bound,
\begin{equation}
\left|
\frac{\widehat v_K^2}{v_K^2}
-1
\right|
\le
\|\widehat{\bm G}_K^{-1}-\bm I_{d_K}\|_{\mathrm{op}}
=
o_P(1).
\label{eq:OLS-variance-ratio}
\end{equation}
Let \(\bm e_0\) select the retained intercept. Since \(\bm e_0^\top\bm s_K(\widehat{\bm x})=1\) and \(\bm e_0^\top\bm\Sigma_K\bm e_0=1\), the Cauchy--Schwarz inequality gives
\(1\le(\bm e_0^\top\bm\Sigma_K\bm e_0)v_K^2=v_K^2\).
Consequently, \(\widehat v_K/v_K\xrightarrow{P}1\), \(\widehat v_K^2=O_P(d_K)\), and \(\widehat v_K^2\) is bounded away from zero with probability approaching one.

The same relative moment condition controls the maximal sample leverage. On the event \(\|\widehat{\bm G}_K-\bm I_{d_K}\|_{\mathrm{op}}\le1/2\),
\begin{equation}
\bm s_i^\top(\mathbf S^\top\mathbf S)^{-1}\bm s_i
=
\frac{1}{n}\bm u_i^\top\widehat{\bm G}_K^{-1}\bm u_i
\le
\frac{2}{n}\|\bm u_i\|^2.
\end{equation}
Thus, for every \(\delta>0\), the union bound and Markov's inequality give
\begin{equation}
\begin{aligned}
\mathbb P\left(
\max_{1\le i\le n}
\bm s_i^\top(\mathbf S^\top\mathbf S)^{-1}\bm s_i
>
\delta
\right)
&\le
o(1)
+
n\mathbb P\left(
\|\bm u_i\|^2>\frac{\delta n}{2}
\right)\\
&\le
o(1)
+
\frac{4\mathbb E\|\bm u_i\|^4}{\delta^2n}
=
o(1).
\end{aligned}
\label{eq:OLS-maximal-leverage}
\end{equation}

Using \autoref{eq:sample-path-regression-model}, decompose the prediction error as
\begin{equation}
\label{eqn:decomp_fhat}
\begin{aligned}
\widehat f_K(\widehat{\bm x})
-
f_0(\widehat{\bm x})
&=
\bm s_K(\widehat{\bm x})^\top
(\mathbf S^\top\mathbf S)^{-1}
\mathbf S^\top\bm\varepsilon
+
\bm s_K(\widehat{\bm x})^\top
(\mathbf S^\top\mathbf S)^{-1}
\mathbf S^\top\bm\xi_K
-
\xi_K^{\bm x}(f_0)\\
&=:
I_1+I_2+I_3.
\end{aligned}
\end{equation}
By \autoref{assu: PathReg}(ii), \(\mathbb E[n^{-1}\|\bm\xi_K\|^2]=O(K^{-Q})\), and hence \(n^{-1}\|\bm\xi_K\|^2=O_P(K^{-Q})\) by Markov's inequality. Since \(\mathbf S(\mathbf S^\top\mathbf S)^{-1}\mathbf S^\top\) is an orthogonal projector, conditional Chebyshev's inequality and Cauchy--Schwarz give
\begin{equation}
\begin{aligned}
&\mathbb E[I_1^2\mid\mathcal X_n]
=
\frac{\sigma^2}{n}\widehat v_K^2
\qquad \Longrightarrow \qquad
I_1^2
=
O_P\left(\frac{d_K}{n}\right),\\
&I_2^2
\le
\frac{\widehat v_K^2}{n}
\bm\xi_K^\top
\mathbf S(\mathbf S^\top\mathbf S)^{-1}
\mathbf S^\top\bm\xi_K
\le
\frac{\widehat v_K^2}{n}\|\bm\xi_K\|^2
=
O_P\left(\frac{d_K}{K^Q}\right),
\end{aligned}
\end{equation}
and $I_3^2
=
O(K^{-Q})$.
Moreover,
$d_K^2
=
\left\{\mathbb E\|\bm u_i\|^2\right\}^2
\le
\mathbb E\|\bm u_i\|^4,
$ 
so one has \(d_K^2/n\to0\). Since \((I_1+I_2+I_3)^2\le3(I_1^2+I_2^2+I_3^2)\), the preceding bounds and \(d_K/K^Q\to0\) prove part (i).

(ii) Define the \(\mathcal X_n\)-measurable weights
\begin{equation}
\begin{aligned}
&w_{ni}
:=
\frac{
\bm s_K(\widehat{\bm x})^\top
\widehat{\bm\Sigma}_K^{-1}\bm s_i
}{
\sqrt n\,\widehat v_K
},
\quad 
\frac{\sqrt n\,I_1}{\widehat v_K}
=
\sum_{i=1}^n w_{ni}\varepsilon_i,
\quad
\sum_{i=1}^n w_{ni}^2
=
1,\\
&\max_{1\le i\le n}w_{ni}^2
\le
\max_{1\le i\le n}
\bm s_i^\top
(\mathbf S^\top\mathbf S)^{-1}
\bm s_i
=
o_P(1).
\end{aligned}
\end{equation}
Conditional on \(\mathcal X_n\), the weights are fixed and the errors are independent, centered, and have variance \(\sigma^2\). Hence, their conditional Lyapunov ratio satisfies
\begin{equation}
\frac{1}{\sigma^4}
\sum_{i=1}^n
\mathbb E\left[
|w_{ni}\varepsilon_i|^4
\mid
\mathcal X_n
\right]
\le
\frac{C_\varepsilon}{\sigma^4}
\left(
\max_{1\le i\le n}w_{ni}^2
\right)
\sum_{i=1}^n w_{ni}^2
=
o_P(1).
\end{equation}
Every subsequence therefore contains a further subsequence along which this Lyapunov ratio converges to zero almost surely. Conditional Lindeberg--Feller then gives almost-sure convergence of the conditional characteristic function along this further subsequence to \(\exp(-\sigma^2t^2/2)\). Hence, the conditional characteristic functions converge to this limit in probability along the original sequence and, since they are uniformly bounded, also in \(L^1\). Taking expectations and applying L\'evy's continuity theorem yields \(\sqrt n\,I_1/\widehat v_K\xrightarrow{d}\mathcal N(0,\sigma^2)\).

Finally, the projection bound above and the fact that \(\widehat v_K^2\) is bounded away from zero with probability approaching one give
\begin{equation}
\frac{nI_2^2}{\widehat v_K^2}
\le
\|\bm\xi_K\|^2
=
O_P\left(\frac{n}{K^Q}\right),
\qquad
\frac{nI_3^2}{\widehat v_K^2}
=
O_P\left(\frac{n}{K^Q}\right).
\end{equation}
Hence, \(n/K^Q\to0\) implies \(\sqrt n\,I_2/\widehat v_K=o_P(1)\) and \(\sqrt n\,I_3/\widehat v_K=o_P(1)\). Combining these bounds with \autoref{eqn:decomp_fhat}, the central limit theorem above, and \(\widehat v_K/v_K\xrightarrow{P}1\), we obtain
\begin{equation}
\begin{aligned}
\frac{
\sqrt n\{\widehat f_K(\widehat{\bm x})-f_0(\widehat{\bm x})\}
}{
v_K
}
&=
\frac{\widehat v_K}{v_K}
\left(
\frac{\sqrt n\,I_1}{\widehat v_K}
+
\frac{\sqrt n\,I_2}{\widehat v_K}
+
\frac{\sqrt n\,I_3}{\widehat v_K}
\right)
\xrightarrow{d}
\mathcal N(0,\sigma^2),
\end{aligned}
\end{equation}
which proves part (ii).
\QED

\subsection{Proof of \autoref{thm:lasso_finite_sample}}

\prf
Let
\(\bm D_K:=\operatorname{diag}(w_{K1},\ldots,w_{Kd_K})\), where \(T>0\) and
\(w_{Kj}=T^{k_j/2}/\sqrt{k_j!}\), and write
\(\overline{\bm s}_K:=\bm D_K^{-1}\bm s_K\),
\(\overline{\mathbf S}:=\mathbf S\bm D_K^{-1}\),
\(\bm\Gamma_K:=
\mathbb E[\overline{\bm s}_K(\widehat{\bm X})
\overline{\bm s}_K(\widehat{\bm X})^\top]\),
\(\bm\beta_{0K}:=\bm D_K\bm L_{0K}\), and
\(\widehat{\bm\beta}_K:=\bm D_K\widehat{\bm L}_K\).
Since \(\bm D_K\) is nonsingular, the weighted LASSO problem is equivalent to
\begin{equation}
\label{eq:weighted-lasso-equivalent}
\widehat{\bm\beta}_K
\in
\arg\min_{\bm\beta\in\mathbb R^{d_K}}
\left\{
\frac{1}{2n}
\|\mathbf Z-\overline{\mathbf S}\bm\beta\|^2
+\lambda_n\|\bm\beta\|_1
\right\}.
\end{equation}
This transformation preserves supports and signs. It also preserves the population normal equation, so
\(\mathbb E[\overline{\bm s}_K(\widehat{\bm X})
\xi_K^{\bm X}(f_0)]=\bm0\).

For \(j=1,\ldots,d_K\), let \(\overline s_{ij}\) denote the \(j\)-th coordinate of
\(\overline{\bm s}_K(\widehat{\bm X}_i)\), and let
\(k_j\in\{0,\ldots,K\}\) be the length of the corresponding signature word.
For \(k_j=0\), \(\overline s_{ij}=1\). For \(k_j\ge1\),
\autoref{lem:sigmoment} and the definition of \(w_{Kj}\) give
\begin{equation}
\label{eq:weighted-signature-moment}
\|\overline s_{ij}\|_p
\le
(C_{\mathrm{sig}}\sqrt p)^{k_j},
\qquad p\ge2.
\end{equation}
Fix coordinates \(j,m\), and set \(r:=k_j+k_m\). When \(r=0\), the centered product is identically zero. For \(r\ge1\), H\"older's inequality and \autoref{eq:weighted-signature-moment} give
\begin{equation}
\label{eq:weighted-lasso-product-moment}
\left\|
\overline s_{ij}\overline s_{im}
-\mathbb E(\overline s_{ij}\overline s_{im})
\right\|_p
\le
C^r p^{r/2},
\qquad p\ge2,
\end{equation}
where \(C>0\) depends only on \(R\) and \(T\).

The moment-to-\(\psi\) argument used in the proof of \autoref{lem:sigtail}
shows that \(\|W\|_p\le B_rp^{r/2}\) implies
\(\|W\|_{\psi_{2/r}}\le D^rB_r\) for a universal constant \(D\).
Thus, the centered product in \autoref{eq:weighted-lasso-product-moment}
has \(\psi_{2/r}\)-functional bounded by \(C^r\).
For \(r=1,2\), the usual sub-Gaussian and sub-exponential inequalities apply. Since the observations are i.i.d., for \(r\ge3\), applying \citep[Corollary~6.4]{ZhangandChen2020} with tail parameter \(2/r\) and equal weights \(1/n\) gives, after enlarging \(C\),
\begin{equation}
\label{eq:weighted-lasso-product-concentration}
\mathbb P\left(
\left|
\frac1n\sum_{i=1}^n
\left\{
\overline s_{ij}\overline s_{im}
-\mathbb E(\overline s_{ij}\overline s_{im})
\right\}
\right|
>
(Cr)^{r/2}
\left\{
\sqrt{\frac{t}{n}}
+\frac{t^{r/2}}{n}
\right\}
\right)
\le2e^{-t}
\end{equation}
for every \(t\ge1\). Indeed, the concentration constant in that result satisfies
\(C(2/r)4^{r/2}\le(Cr)^{r/2}\); the same upper bound also absorbs the coefficient of the square-root term.

Define $
\Delta_{nK}
:=
(CK)^K
\left\{
\sqrt{\frac{\log(e d_K)}{n}}
+\frac{\{\log(e d_K)\}^K}{n}
\right\}$,
where \(C>0\) is a sufficiently large fixed constant. The rate and penalty conditions in \autoref{thm:lasso_finite_sample} are equivalently
\(q_K\Delta_{nK}\to0\) and
\(\lambda_n\asymp\sqrt{q_K}\Delta_{nK}\), with a sufficiently large leading constant in the latter relation. Write \(\|\cdot\|_{\max}\) for the entrywise maximum norm. Taking
\(t=A\log(e d_K)\) in \autoref{eq:weighted-lasso-product-concentration} for a sufficiently large fixed \(A\), using \(r\le2K\), and applying a union bound over all pairs of coordinates yield
\begin{equation}
\label{eq:weighted-lasso-gram-concentration}
\left\|
\widehat{\bm\Gamma}_K-\bm\Gamma_K
\right\|_{\max}
=O_P(\Delta_{nK}),
\qquad
\widehat{\bm\Gamma}_K
:=
\frac1n\overline{\mathbf S}^{\top}\overline{\mathbf S}.
\end{equation}
Since \(d_K\to\infty\), the constants can be chosen so that the event on which the left-hand side is bounded by a fixed multiple of \(\Delta_{nK}\) has probability approaching one.

We next transfer the population design conditions to the sample. Write
\(\bm A:=(\bm\Gamma_K)_{A_0A_0}\),
\(\widehat{\bm A}:=(\widehat{\bm\Gamma}_K)_{A_0A_0}\), and
\(\bm B:=(\bm\Gamma_K)_{A_0^cA_0}\), with
\(\widehat{\bm B}\) defined analogously. By \autoref{assu:lasso}(i), there exist proof constants \(C_\Gamma<\infty\) and \(\eta\in(0,1]\) such that
\(\|\bm A^{-1}\|_\infty\le C_\Gamma\) and
\(\|\bm B\bm A^{-1}\|_\infty\le1-\eta\) uniformly in \(K\). By
\autoref{eq:weighted-lasso-gram-concentration},
\begin{equation}
\left\|
(\widehat{\bm A}-\bm A)\bm A^{-1}
\right\|_\infty
\le
q_K
\left\|
\widehat{\bm\Gamma}_K-\bm\Gamma_K
\right\|_{\max}
\|\bm A^{-1}\|_\infty
=o_P(1).
\end{equation}
The inverse perturbation identity and \autoref{assu:lasso}(i) therefore imply that \(\widehat{\bm A}\) is nonsingular and
\(\|\widehat{\bm A}^{-1}\|_\infty\le2C_\Gamma\) on an event whose probability tends to one.
Since \(\widehat{\bm A}\) is a sample Gram matrix and hence positive semidefinite, it is positive definite on this event. Furthermore,
\begin{equation}
\label{eq:weighted-lasso-irc-perturbation}
\begin{aligned}
&\,\,\,\widehat{\bm B}\widehat{\bm A}^{-1}
-\bm B\bm A^{-1}
=
\left\{
(\widehat{\bm B}-\bm B)
-\bm B\bm A^{-1}(\widehat{\bm A}-\bm A)
\right\}
\widehat{\bm A}^{-1},\\
&\left\|
\widehat{\bm B}\widehat{\bm A}^{-1}
-\bm B\bm A^{-1}
\right\|_\infty
\le
\left\{
1+\|\bm B\bm A^{-1}\|_\infty
\right\}
q_K
\left\|
\widehat{\bm\Gamma}_K-\bm\Gamma_K
\right\|_{\max}
\|\widehat{\bm A}^{-1}\|_\infty
=o_P(1).
\end{aligned}
\end{equation}
Consequently, with probability approaching one,
\begin{equation}
\label{eq:weighted-lasso-sample-design}
\left\|
(\widehat{\bm\Gamma}_K)_{A_0A_0}^{-1}
\right\|_\infty
\le2C_\Gamma,
\qquad
\left\|
(\widehat{\bm\Gamma}_K)_{A_0^cA_0}
(\widehat{\bm\Gamma}_K)_{A_0A_0}^{-1}
\right\|_\infty
\le1-\frac{\eta}{2}.
\end{equation}

We now control the effective score. Since
\(\overline{\mathbf S}\) is \(\mathcal X_n\)-measurable, conditional independence and \autoref{assu:lasso}(iii) imply, for every \(j\) and \(t\in\mathbb R\),
\begin{equation}
\mathbb E\left[
\exp\left\{
\frac{t}{n}
\sum_{i=1}^n\overline s_{ij}\varepsilon_i
\right\}
\biggm|\mathcal X_n
\right]
\le
\exp\left\{
\frac{\sigma_0^2t^2}{2n^2}
\sum_{i=1}^n\overline s_{ij}^2
\right\}.
\end{equation}
By \autoref{eq:weighted-signature-moment} and
\autoref{eq:weighted-lasso-gram-concentration},
\(\max_j n^{-1}\sum_i\overline s_{ij}^2=O_P(C^{2K})\).
The preceding conditional moment-generating-function bound and a union bound over \(j\) therefore give
\begin{equation}
\label{eq:weighted-lasso-noise-score}
\frac1n
\left\|
\overline{\mathbf S}^{\top}\bm\varepsilon
\right\|_\infty
=
O_P\left(
C^K\sqrt{\frac{\log(e d_K)}{n}}
\right)
=O_P(\Delta_{nK}).
\end{equation}

The bound
\(\|\bm A^{-1}\|_\infty\le C_\Gamma\) implies
\(\lambda_{\min}(\bm A)\ge C_\Gamma^{-1}\). Let
\(Y:=f_0(\widehat{\bm X})\), with \(Y_i:=f_0(\widehat{\bm X}_i)\).
The population normal equation and
\(\mathbb E(Y^2)\le C\) give
\begin{equation}
\label{eq:weighted-lasso-coefficient-bound}
\begin{aligned}
C_\Gamma^{-1}
\|\bm\beta_{0K}\|_2^2
&\le
\bm\beta_{0K,A_0}^{\top}
\bm A
\bm\beta_{0K,A_0}=\mathbb E\left[
Y\,
\overline{\bm s}_{K,A_0}(\widehat{\bm X})^{\top}
\bm\beta_{0K,A_0}
\right]\\
&\le
\{\mathbb E(Y^2)\}^{1/2}
\left\{
\bm\beta_{0K,A_0}^{\top}
\bm A
\bm\beta_{0K,A_0}
\right\}^{1/2}.
\end{aligned}
\end{equation}
Hence,
\(\|\bm\beta_{0K}\|_2=O(1)\) and
\(\|\bm\beta_{0K}\|_1=O(\sqrt{q_K})\).

The Stratonovich--It\^{o} conversion expresses \(Y\) as bounded finite-variation terms plus an It\^{o} integral with a uniformly bounded integrand; the Burkholder--Davis--Gundy inequality therefore gives
\(\|Y\|_p\le C\sqrt p\) for every \(p\ge2\). Thus, by
\autoref{eq:weighted-signature-moment} and H\"older's inequality,
$
\left\|
\overline s_{ij}Y_i
-\mathbb E(\overline s_{ij}Y_i)
\right\|_p
\le
C^{k_j+1}p^{(k_j+1)/2}$.
Since \(k_j+1\le K+1\le2K\), the same sub-Weibull concentration argument as above, followed by a union bound over \(j\), yields
\begin{equation}
\label{eq:weighted-lasso-response-concentration}
\max_{1\le j\le d_K}
\left|
\frac1n\sum_{i=1}^n
\left\{
\overline s_{ij}Y_i
-\mathbb E(\overline s_{ij}Y_i)
\right\}
\right|
=O_P(\Delta_{nK}),
\end{equation}
with a corresponding fixed-multiple bound on an event whose probability approaches one.

Since
\(\xi_K^{\bm X_i}(f_0)
=Y_i-\overline{\bm s}_K(\widehat{\bm X}_i)^{\top}\bm\beta_{0K}\), the population normal equation gives, for every \(j\),
\begin{equation}
\begin{aligned}
\frac1n\sum_{i=1}^n
\overline s_{ij}\xi_K^{\bm X_i}(f_0)
&=
\frac1n\sum_{i=1}^n
\left\{
\overline s_{ij}Y_i
-\mathbb E(\overline s_{ij}Y_i)
\right\}\\
&\quad-
\sum_{m\in A_0}\beta_{0K,m}
\frac1n\sum_{i=1}^n
\left\{
\overline s_{ij}\overline s_{im}
-\mathbb E(\overline s_{ij}\overline s_{im})
\right\}.
\end{aligned}
\end{equation}
Combining \autoref{eq:weighted-lasso-gram-concentration},
\autoref{eq:weighted-lasso-response-concentration}, and
\(\|\bm\beta_{0K}\|_1=O(\sqrt{q_K})\) yields
\begin{equation}
\label{eq:weighted-lasso-approximation-score}
\frac1n
\left\|
\overline{\mathbf S}^{\top}\bm\xi_K
\right\|_\infty
=O_P(\sqrt{q_K}\Delta_{nK}).
\end{equation}
The explicit exponential bounds used above allow all constants to be chosen uniformly. Therefore, for
\(\bm u:=\bm\xi_K+\bm\varepsilon\), the stated penalty choice gives
\begin{equation}
\label{eq:weighted-lasso-score}
\frac1n
\left\|
\overline{\mathbf S}^{\top}\bm u
\right\|_\infty
\le
\frac{\eta\lambda_n}{4}
\end{equation}
on an event whose probability approaches one.

Work on the event on which
\autoref{eq:weighted-lasso-sample-design} and
\autoref{eq:weighted-lasso-score} hold.

\noindent\textit{(i) Support containment and uniqueness.}
Let
\(\widetilde{\bm\beta}_K\) solve the restricted problem on \(A_0\), with
\(\widetilde{\bm\beta}_{K,A_0^c}=\bm0\). It is unique because
\((\widehat{\bm\Gamma}_K)_{A_0A_0}\) is positive definite. Its KKT condition is
\begin{equation}
\label{eq:weighted-lasso-restricted-kkt}
(\widehat{\bm\Gamma}_K)_{A_0A_0}
\left(
\widetilde{\bm\beta}_{K,A_0}
-\bm\beta_{0K,A_0}
\right)
=
\frac1n
\overline{\mathbf S}_{A_0}^{\top}\bm u
-\lambda_n\bm z_{A_0},
\qquad
\|\bm z_{A_0}\|_\infty\le1.
\end{equation}
Define the inactive dual vector by
\begin{equation}
\bm z_{A_0^c}
:=
\frac{
\overline{\mathbf S}_{A_0^c}^{\top}
\{\mathbf Z-
\overline{\mathbf S}_{A_0}\widetilde{\bm\beta}_{K,A_0}\}
}{n\lambda_n}.
\end{equation}
Substituting \autoref{eq:weighted-lasso-restricted-kkt} and using
\autoref{eq:weighted-lasso-sample-design}--\autoref{eq:weighted-lasso-score} give
\begin{equation}
\label{eq:weighted-lasso-strict-dual}
\begin{aligned}
\|\bm z_{A_0^c}\|_\infty
&\le
\frac{
\|\overline{\mathbf S}_{A_0^c}^{\top}\bm u\|_\infty
}{n\lambda_n}+
\left\|
(\widehat{\bm\Gamma}_K)_{A_0^cA_0}
(\widehat{\bm\Gamma}_K)_{A_0A_0}^{-1}
\right\|_\infty
\left\{
\frac{
\|\overline{\mathbf S}_{A_0}^{\top}\bm u\|_\infty
}{n\lambda_n}
+\|\bm z_{A_0}\|_\infty
\right\}\\
&\le
\frac{\eta}{4}
+\left(1-\frac{\eta}{2}\right)
\left(1+\frac{\eta}{4}\right)
=1-\frac{\eta^2}{8}
<1.
\end{aligned}
\end{equation}
Thus, \((\widetilde{\bm\beta}_K,\bm z)\) satisfies the KKT conditions for the full problem, with strict dual feasibility on \(A_0^c\).

Suppose that \(\bm\beta^*\) is another minimizer. Convexity of the squared loss in the fitted value implies
\(\overline{\mathbf S}\bm\beta^*
=\overline{\mathbf S}\widetilde{\bm\beta}_K\). The two minimizers have the same \(\ell_1\)-norm, and the KKT equation gives
\(\bm z^{\top}(\bm\beta^*-\widetilde{\bm\beta}_K)=0\). Equality in the subgradient inequality, together with
\autoref{eq:weighted-lasso-strict-dual}, forces
\(\bm\beta^*_{A_0^c}=\bm0\). Positive definiteness of
\((\widehat{\bm\Gamma}_K)_{A_0A_0}\) then gives
\(\bm\beta^*=\widetilde{\bm\beta}_K\). Hence, the full solution is unique with probability approaching one and
\(\operatorname{supp}(\widehat{\bm\beta}_K)\subseteq A_0\).
Because \(\bm D_K\) is diagonal with positive entries, this is equivalent to
\(\operatorname{supp}(\widehat{\bm L}_K)\subseteq A_0\), which proves (i).

\noindent\textit{(ii) Weighted coefficient rates.}
It follows from
\autoref{eq:weighted-lasso-restricted-kkt},
\autoref{eq:weighted-lasso-sample-design}, and
\autoref{eq:weighted-lasso-score} that, for some fixed \(C_0>0\),
\begin{equation}
\label{eq:weighted-lasso-restricted-error}
\left\|
\widetilde{\bm\beta}_{K,A_0}
-\bm\beta_{0K,A_0}
\right\|_\infty
\le C_0\lambda_n,
\qquad
\left\|
\widetilde{\bm\beta}_{K,A_0}
-\bm\beta_{0K,A_0}
\right\|_2
\le C_0\sqrt{q_K}\lambda_n.
\end{equation}
Since \(\widehat{\bm\beta}_K=\widetilde{\bm\beta}_K\) and both
\(\widehat{\bm\beta}_{K,A_0^c}\) and
\(\bm\beta_{0K,A_0^c}\) vanish,
\begin{equation}
\label{eq:weighted-lasso-support-error}
\begin{aligned}
\|\widehat{\bm\beta}_K-\bm\beta_{0K}\|_\infty
&=O_P(\lambda_n),
&
\|\widehat{\bm\beta}_K-\bm\beta_{0K}\|_2
&=O_P(\sqrt{q_K}\lambda_n).
\end{aligned}
\end{equation}
Since
\(\widehat{\bm\beta}_K-\bm\beta_{0K}
=\bm D_K(\widehat{\bm L}_K-\bm L_{0K})\), its \(\ell_\infty\)-norm is
\(\max_j w_{Kj}|\widehat L_{K,j}-L_{0K,j}|\), while its squared
\(\ell_2\)-norm is
\(\sum_jw_{Kj}^2(\widehat L_{K,j}-L_{0K,j})^2\).
This proves (ii).

\noindent\textit{(iii) Pointwise prediction.}
Write \(\bm D_{K,A_0}:=(\bm D_K)_{A_0A_0}\). Since
\begin{equation}
\begin{aligned}
\bm s_{K,A_0}(\widehat{\bm x})
=\bm D_{K,A_0}\overline{\bm s}_{K,A_0}(\widehat{\bm x}),\quad
(\bm\Sigma_K)_{A_0A_0}
=\bm D_{K,A_0}(\bm\Gamma_K)_{A_0A_0}\bm D_{K,A_0},
\end{aligned}
\end{equation}
the raw-Gram leverage condition in the main text is exactly equivalent to
\begin{equation}
\begin{aligned}
&\bm s_{K,A_0}(\widehat{\bm x})^\top
(\bm\Sigma_K)_{A_0A_0}^{-1}
\bm s_{K,A_0}(\widehat{\bm x})=
\overline{\bm s}_{K,A_0}(\widehat{\bm x})^\top
(\bm\Gamma_K)_{A_0A_0}^{-1}
\overline{\bm s}_{K,A_0}(\widehat{\bm x})
=O(1).
\end{aligned}
\end{equation}
The Gram concentration bound and
\(q_K\Delta_{nK}\to0\) imply
$(\widehat{\bm\Gamma}_K)_{A_0A_0}
\succeq
\frac12
(\bm\Gamma_K)_{A_0A_0}$
with probability approaching one. Indeed,
\begin{equation}
\|\bm A^{-1/2}(\widehat{\bm A}-\bm A)\bm A^{-1/2}\|_{\mathrm{op}}
\le C_\Gamma q_K
\|\widehat{\bm\Gamma}_K-\bm\Gamma_K\|_{\max}
=o_P(1).
\end{equation}
Hence, the test-path leverage condition gives
$
\overline{\bm s}_{K,A_0}(\widehat{\bm x})^{\top}
(\widehat{\bm\Gamma}_K)_{A_0A_0}^{-1}
\overline{\bm s}_{K,A_0}(\widehat{\bm x})
=O_P(1).
$
Moreover, \autoref{eq:weighted-lasso-restricted-kkt},
\autoref{eq:weighted-lasso-sample-design}, and
\autoref{eq:weighted-lasso-score} imply
\begin{equation}
\label{eq:weighted-lasso-active-prediction}
\begin{aligned}
&\left(
\widehat{\bm\beta}_{K,A_0}-\bm\beta_{0K,A_0}
\right)^{\top}
(\widehat{\bm\Gamma}_K)_{A_0A_0}
\left(
\widehat{\bm\beta}_{K,A_0}-\bm\beta_{0K,A_0}
\right)\\
&\quad=
\left\{
\frac1n\overline{\mathbf S}_{A_0}^{\top}\bm u
-\lambda_n\bm z_{A_0}
\right\}^{\top}
(\widehat{\bm\Gamma}_K)_{A_0A_0}^{-1}
\left\{
\frac1n\overline{\mathbf S}_{A_0}^{\top}\bm u
-\lambda_n\bm z_{A_0}
\right\}
=O_P(q_K\lambda_n^2).
\end{aligned}
\end{equation}
Cauchy--Schwarz in the
\((\widehat{\bm\Gamma}_K)_{A_0A_0}\)-geometry therefore gives
\begin{equation}
\label{eq:weighted-lasso-pointwise-prediction}
\begin{aligned}
\left|
\widehat f_K(\widehat{\bm x})
-f_0(\widehat{\bm x})
\right|^2
&\le
2\left|
\overline{\bm s}_{K,A_0}(\widehat{\bm x})^{\top}
\left(
\widehat{\bm\beta}_{K,A_0}
-\bm\beta_{0K,A_0}
\right)
\right|^2
+2|\xi_K^{\bm x}(f_0)|^2\\
&=
O_P\left(q_K\lambda_n^2+K^{-Q}\right).
\end{aligned}
\end{equation}
Finally,
\(q_K\lambda_n^2
\asymp(q_K\Delta_{nK})^2=o(1)\), while
\(K^{-Q}\to0\). This proves (iii).

\noindent\textit{(iv) Sign consistency.}
On the same event, \autoref{eq:weighted-lasso-restricted-error} gives
$\max_{j\in A_0}
w_{Kj}|\widehat L_{K,j}-L_{0K,j}|
\le C_0\lambda_n$.
Thus, if
\(\min_{j\in A_0}w_{Kj}|L_{0K,j}|>C_\beta\lambda_n\)
for a sufficiently large fixed constant \(C_\beta>C_0\), every active coefficient retains its sign. By (i), all inactive estimated coefficients vanish. Hence,
\(\mathbb P\{\operatorname{sign}(\widehat{\bm L}_K)
=\operatorname{sign}(\bm L_{0K})\}\to1\), proving (iv).
\QED

\subsection{Proof of \autoref{thm:classification-consistency}}

\prf
Throughout the proof, write
\(\bm s_i:=\bm s_K(\widehat{\bm X}_i)\) and
\(\bm u_i:=\bm\Sigma_K^{-1/2}\bm s_i\), with
\(\bm u:=\bm\Sigma_K^{-1/2}\bm s_K(\widehat{\bm X})\) denoting the
population counterpart. Then
\(\mathbb E[\bm u_i\bm u_i^\top]=\bm I_{d_K}\) and
\(\mathbb E\|\bm u_i\|^2=d_K\). Independence gives
\begin{equation}
\label{eq:classification-relative-Gram}
\begin{aligned}
\mathbb E\left[
\left\|
\bm\Sigma_K^{-1/2}
\left(\frac{1}{n}\mathbf S^\top\mathbf S\right)
\bm\Sigma_K^{-1/2}
-\bm I_{d_K}
\right\|_F^2
\right]
&=
\frac{1}{n}
\left\{
\mathbb E\|\bm u_i\|^4-d_K
\right\}
=o(1),
\end{aligned}
\end{equation}
where the last equality follows from
\autoref{assu:classification}(i). Hence,
the operator norm of the matrix inside the norm is \(o_P(1)\), so
\(\mathbf S\) has full column rank with probability approaching one.

We next work in the population Gram geometry. For
\(\bm\theta\in\mathbb R^{d_K}\), define
\begin{equation}
\begin{aligned}
&\mathcal R_{nK}(\bm\theta)
:=
\frac{1}{n}\sum_{i=1}^n
\left[
\log\left\{1+\exp\left(
\bm s_i^\top\bm L_{0K}+\bm u_i^\top\bm\theta
\right)\right\}
-Z_i\left(
\bm s_i^\top\bm L_{0K}+\bm u_i^\top\bm\theta
\right)
\right],\\
&\mathcal R_K(\bm\theta)
:=
\mathbb E\left[
\log\left\{1+\exp\left(
\bm s_K(\widehat{\bm X})^\top\bm L_{0K}
+\bm u^\top\bm\theta
\right)\right\}
-Z\left\{
\bm s_K(\widehat{\bm X})^\top\bm L_{0K}
+\bm u^\top\bm\theta
\right\}
\right].
\end{aligned}
\label{eq:classification-standardized-risks}
\end{equation}
The population score at \(\bm0\) is
\begin{equation}
\bm g_K
:=
\nabla\mathcal R_K(\bm0)
=
\mathbb E\left[
\bm u
\left[
\Lambda\{\bm s_K(\widehat{\bm X})^\top\bm L_{0K}\}
-
\Lambda\{f_0(\widehat{\bm X})\}
\right]
\right].
\end{equation}
Since \(\Lambda\) is \(1/4\)-Lipschitz, Cauchy--Schwarz and
\autoref{assu:classification}(iii) imply
\begin{equation}
\label{eq:classification-population-score}
\begin{aligned}
\|\bm g_K\|
&=
\sup_{\|\bm a\|\le1}
\left|
\mathbb E\left[
(\bm a^\top\bm u)
\left[
\Lambda\{\bm s_K(\widehat{\bm X})^\top\bm L_{0K}\}
-
\Lambda\{f_0(\widehat{\bm X})\}
\right]
\right]
\right|\\
&\le
\frac{1}{4}
\sup_{\|\bm a\|\le1}
\{\mathbb E(\bm a^\top\bm u)^2\}^{1/2}
\{\mathbb E[\{\xi_K^{\bm X}(f_0)\}^2]\}^{1/2}
=O(K^{-Q/2}).
\end{aligned}
\end{equation}
Moreover,
\begin{equation}
\nabla^2\mathcal R_K(\bm\theta)
=
\bm\Sigma_K^{-1/2}
\bm H_K
\left(
\bm L_{0K}+\bm\Sigma_K^{-1/2}\bm\theta
\right)
\bm\Sigma_K^{-1/2}.
\end{equation}
Thus, \autoref{assu:classification}(ii) and Taylor's theorem give,
uniformly over \(\|\bm\theta\|\le r\),
\begin{equation}
\label{eq:classification-population-curvature}
\mathcal R_K(\bm\theta)-\mathcal R_K(\bm0)
\ge
-\|\bm g_K\|\|\bm\theta\|
+
\frac{c}{2}\|\bm\theta\|^2.
\end{equation}

It remains to control the empirical loss increments. For
\(0<\rho\le r\), let
\begin{equation}
\label{eq:classification-loss-increments}
\begin{aligned}
\Delta_{nK}(\rho)
:=
\sup_{\|\bm\theta\|\le\rho}
\Big|
&\{\mathcal R_{nK}(\bm\theta)-\mathcal R_{nK}(\bm0)\}-
\{\mathcal R_K(\bm\theta)-\mathcal R_K(\bm0)\}
\Big|.
\end{aligned}
\end{equation}
Conditional on the observations, the \(i\)th loss increment, viewed
as a function of \(\bm u_i^\top\bm\theta\), vanishes at zero and is
1-Lipschitz. Let \(\zeta_1,\ldots,\zeta_n\) be i.i.d. Rademacher
variables independent of the observations. Symmetrization,
the contraction inequality, and Jensen's inequality yield
\begin{equation}
\label{eq:classification-increment-bound}
\begin{aligned}
\mathbb E[\Delta_{nK}(\rho)]
&\le
4\,
\mathbb E\left[
\sup_{\|\bm\theta\|\le\rho}
\left|
\frac{1}{n}
\sum_{i=1}^n
\zeta_i\bm u_i^\top\bm\theta
\right|
\right]=
4\rho\,
\mathbb E\left\|
\frac{1}{n}
\sum_{i=1}^n\zeta_i\bm u_i
\right\|\\
&\le
4\rho
\left\{
\mathbb E\left\|
\frac{1}{n}
\sum_{i=1}^n\zeta_i\bm u_i
\right\|^2
\right\}^{1/2}
=
4\rho\sqrt{\frac{d_K}{n}}.
\end{aligned}
\end{equation}

We first establish existence and uniqueness. Since
\(\mathbb E\|\bm u_i\|^2=d_K\), Jensen's inequality and
\autoref{assu:classification}(i) give
\begin{equation}
\frac{d_K^2}{n}
\le
\frac{1}{n}\mathbb E\|\bm u_i\|^4
\to0.
\label{eq:classification-dimension-growth}
\end{equation}
Together with \(K\to\infty\),
\autoref{eq:classification-population-score} and
\autoref{eq:classification-increment-bound} imply
\(\|\bm g_K\|=o(1)\) and \(\Delta_{nK}(r)=o_P(1)\).
Consequently, \autoref{eq:classification-population-curvature} gives,
with probability approaching one,
\begin{equation}
\begin{aligned}
\inf_{\|\bm\theta\|=r}
\left\{
\mathcal R_{nK}(\bm\theta)-\mathcal R_{nK}(\bm0)
\right\}
&\ge
-r\|\bm g_K\|
+
\frac{cr^2}{2}
-
\Delta_{nK}(r)
>0.
\end{aligned}
\label{eq:classification-boundary-separation}
\end{equation}
By continuity, \(\mathcal R_{nK}\) attains its minimum over the closed
ball \(\{\bm\theta:\|\bm\theta\|\le r\}\), and the preceding
separation places this minimum in the interior. Its gradient therefore
vanishes, so convexity makes it a global minimizer. On the event that
\(\mathbf S\) has full column rank, the empirical Hessian
satisfies, for every finite \(\bm\theta\),
\begin{equation}
\nabla^2\mathcal R_{nK}(\bm\theta)
=
\frac{1}{n}
\sum_{i=1}^n
\Lambda(\bm s_i^\top\bm L_{0K}+\bm u_i^\top\bm\theta)
\{1-\Lambda(\bm s_i^\top\bm L_{0K}+\bm u_i^\top\bm\theta)\}
\bm u_i\bm u_i^\top
\succ\bm0.
\label{eq:classification-strict-convexity}
\end{equation}
Hence, \(\mathcal R_{nK}\), and equivalently \(\mathcal L_n\), is
strictly convex. The global minimizer is therefore finite and unique
with probability approaching one. We define \(\widehat{\bm L}_K\)
arbitrarily on the complement of this event.

We finally derive the rate. Let
\(b_{nK}:=\sqrt{d_K/n}+K^{-Q/2}\). By
\autoref{eq:classification-dimension-growth} and \(K\to\infty\),
\(b_{nK}\to0\). For any fixed \(M>0\) and all sufficiently large \(n\),
\(Mb_{nK}\le r\). From
\autoref{eq:classification-population-curvature},
\begin{equation}
\begin{aligned}
\inf_{\|\bm\theta\|=Mb_{nK}}
\left\{
\mathcal R_{nK}(\bm\theta)-\mathcal R_{nK}(\bm0)
\right\}
\ge{}&
-Mb_{nK}\|\bm g_K\|
+
\frac{cM^2b_{nK}^2}{2}-
\Delta_{nK}(Mb_{nK}).
\end{aligned}
\label{eq:classification-rate-boundary}
\end{equation}
By \autoref{eq:classification-population-score},
\(\|\bm g_K\|\le Cb_{nK}\) for some constant \(C>0\) and all
sufficiently large \(n\).
Moreover, \autoref{eq:classification-increment-bound} and Markov's
inequality give
\begin{equation}
\begin{aligned}
\mathbb P\left\{
\Delta_{nK}(Mb_{nK})
>
\frac{cM^2b_{nK}^2}{4}
\right\}
&\le
\frac{16}{cM}
\frac{\sqrt{d_K/n}}{b_{nK}}
\le
\frac{16}{cM}.
\end{aligned}
\label{eq:classification-rate-increment}
\end{equation}
Choose \(M>4C/c\). On the complement of the event in
\autoref{eq:classification-rate-increment}, the right-hand side of
\autoref{eq:classification-rate-boundary} is positive. On the
preceding existence event, convexity and the minimizing property of
\(\widehat{\bm L}_K\) then imply
\(\|\bm\Sigma_K^{1/2}(\widehat{\bm L}_K-\bm L_{0K})\|
\le Mb_{nK}\). Since \(M\) can be made arbitrarily large, we obtain
\begin{equation}
\label{eq:classification-coef-rate}
\left\|
\bm\Sigma_K^{1/2}
(\widehat{\bm L}_K-\bm L_{0K})
\right\|
=
O_P\left(
\sqrt{\frac{d_K}{n}}
+
K^{-Q/2}
\right).
\end{equation}

For the fixed test path \(\widehat{\bm x}\), Cauchy--Schwarz gives
\begin{equation}
\begin{aligned}
\left|
\widehat f_K(\widehat{\bm x})-f_0(\widehat{\bm x})
\right|
\le{}&
\left\{
\bm s_K(\widehat{\bm x})^\top
\bm\Sigma_K^{-1}
\bm s_K(\widehat{\bm x})
\right\}^{1/2}
\left\|
\bm\Sigma_K^{1/2}
(\widehat{\bm L}_K-\bm L_{0K})
\right\|+
|\xi_K^{\bm x}(f_0)|.
\end{aligned}
\end{equation}
The test-path conditions and
\autoref{eq:classification-coef-rate} therefore yield
\begin{equation}
\left|
\widehat f_K(\widehat{\bm x})-f_0(\widehat{\bm x})
\right|
=
O_P\left(
\frac{d_K}{\sqrt n}
+
\frac{\sqrt{d_K}}{K^{Q/2}}
\right).
\end{equation}
The first term converges to zero by
\autoref{eq:classification-dimension-growth}, and the second converges
to zero because \(d_K/K^Q\to0\). This proves the theorem.
\QED
\section{Additional Experimental Details}
\label{app:additional-experiments}

This appendix provides supplementary details for the simulation study in \autoref{sec_simulation} and the three real-data applications in \autoref{sec:real-data}.
The purpose is to document data-generating processes, implementation choices, evaluation metrics, robustness diagnostics, and additional visual evidence that support the main-text findings.

\subsection{Simulation Details}
\label{app:simulation-details}

This subsection records the data-generating processes, implementation choices, and evaluation metrics for \autoref{sec_simulation}.
All reported errors were evaluated on an independent test sample and then averaged over independent Monte Carlo replications.

\subsubsection{Population approximation.}
The first experiment isolates finite-level signature approximation error.
For the Brownian benchmark, we simulated Brownian paths \(B\) on \([0,1]\) and defined $F_{\gamma}^{\mathrm B}(B)
    =
    \int_0^1 g_{\gamma}(t)\,\mathrm dB_t$ where
\begin{equation}
    g_{\gamma}(t)
    =
    \sum_{j=1}^{J} j^{-(\gamma+1/2)}\sin(j\pi t),
    \,\,
    \gamma\in\{1,2\}.
    \label{eq:app-sim-brownian-functional}
\end{equation}
For the fixed finite series used in the simulation, the target belongs to the corresponding smooth coefficient class \(\mathcal U_{\gamma,R}\) for some finite \(R\).
Specifically, define
\begin{equation}
R
=
\sum_{\ell=0}^{\gamma}
\sup_{t\in[0,1]}
\left|
\pi^\ell
\sum_{j=1}^{J}
j^{\ell-\gamma-1/2}
\sin\left(j\pi t+\frac{\ell\pi}{2}\right)
\right|.
\label{eq:app-sim-brownian-radius}
\end{equation}
Then
\((B,h)\in
\mathcal U_{\gamma,R}
\setminus
\mathcal U_{\gamma+1,R}\).

For the diffusion benchmark, we simulated the Ornstein--Uhlenbeck process
\begin{equation}
    \mathrm dX_t=-\theta X_t\,\mathrm dt+\sigma\,\mathrm dB_t,
    \qquad X_0=x_0,
    \label{eq:app-sim-approx-ou-process}
\end{equation}
and used the smooth coefficient functional
\begin{equation}
    F_{\gamma}^{\mathrm{OU}}(X)
    =
    \int_0^1 g_\gamma(t)\circ\mathrm dX_t
    +
    \rho\int_0^1 X_t\,\mathrm dt .
    \label{eq:app-sim-ou-functional}
\end{equation}
This has the same smooth-coefficient form as \autoref{eq:smooth-coefficient-functional}, with \(h_1(t,x)=g_\gamma(t)\) and \(h_0(t,x)=\rho x\), and serves as a non-Brownian diffusion benchmark for the same approximation mechanism.
To match the bounded-coefficient assumption in \(\mathcal U_{\gamma,R}\) literally, the linear term \(x\) can be replaced by a bounded smooth localization that agrees with \(x\) on the simulated state range; the reported computation uses the expression in \autoref{eq:app-sim-ou-functional}.
For each target, each \(\gamma\in\{1,2\}\), and \(K=2,\ldots,9\), we fitted the least-squares projection of the target on \(\bm s_K(\widehat X)\) using simulated training paths and reported the test projection error
\begin{equation}
    \operatorname{MSE}_K
    =
    \frac{1}{N_{\mathrm{test}}}
    \sum_{i=1}^{N_{\mathrm{test}}}
    \left\{
        \widehat f_K(\widehat X_i^{\,\mathrm{test}})
        -
        F_\gamma(X_i^{\mathrm{test}})
    \right\}^2 .
    \label{eq:app-sim-projection-mse}
\end{equation}
The figures use \(12{,}000\) simulated training paths and an independent test pool of \(4{,}000\) paths, with the same numerical discretization and scalar parameters held fixed across \(K\) and \(\gamma\).
The theoretical benchmark is the smoothness-dependent ordering \(\operatorname{MSE}_K=O(K^{-2\gamma})\).

\subsubsection{OLS and logistic regression.}
The second experiment uses the Brownian target \(f(B)=F_2^{\mathrm B}(B)\) from \autoref{eq:app-sim-brownian-functional} as the latent signal.
For OLS, the observed response is
\begin{equation}
    Z_i
    =
    f(B_i)+\varepsilon_i,
    \,\,\,
    \varepsilon_i\sim N(0,\sigma_\varepsilon^2),
    \,\,\,
    \sigma_\varepsilon=0.25 .
    \label{eq:app-sim-ols-response}
\end{equation}
For logistic regression, the response is binary with
\begin{equation}
    Z_i\mid B_i
    \sim
    \operatorname{Bernoulli}\{\Lambda(\alpha f(B_i))\},
    \,\,\,
    \alpha=1.25 .
    \label{eq:app-sim-logit-response}
\end{equation}
Thus the OLS conditional mean and the logistic log-odds were generated by the same smooth path functional.
We fitted OLS and logistic models on \(\bm s_K(\widehat B_i)\) with \(K\in\{2,4,6,8\}\), corresponding to \(d_K\in\{6,30,126,510\}\).
OLS performance was evaluated against the noise-free conditional mean:
\begin{equation}
    \operatorname{RMSE}_{K,n}
    =
    \sqrt{
    \frac{1}{N_{\mathrm{test}}}
    \sum_{i=1}^{N_{\mathrm{test}}}
    \left\{
        \widehat f_{K,n}(\widehat B_i^{\,\mathrm{test}})
        -
        f(B_i^{\mathrm{test}})
    \right\}^2}
    \label{eq:app-sim-ols-rmse}
\end{equation}
For logistic regression, let \(p(B)=\Lambda(\alpha f(B))\), let \(A^*=\mathbb E[\max\{p(B),1-p(B)\}]\) denote the Bayes accuracy, and let \(\widehat h_{K,n}(\widehat B)=\mathbf 1\{\widehat p_{K,n}(\widehat B)\ge 1/2\}\).
We reported the excess classification error
\begin{equation}
    \operatorname{ExErr}_{K,n}
    =
    A^*
    -
    \frac{1}{N_{\mathrm{test}}}
    \sum_{i=1}^{N_{\mathrm{test}}}
    \mathbf 1\{\widehat h_{K,n}(\widehat B_i^{\,\mathrm{test}})=Z_i^{\mathrm{test}}\} .
    \label{eq:app-sim-logit-excess-error}
\end{equation}
The OLS and logistic curves were averaged over independent Monte Carlo replications, using the sample sizes shown in \autoref{fig:sim-summary-regression}.

\subsubsection{Signature-LASSO.}
The third experiment is motivated by the sparse OU representation in \autoref{exmp: OU}, but the observed path used for signatures is an It\^{o} diffusion rather than Brownian motion itself.
We simulated
\begin{equation}
    \mathrm dX_t=-\kappa X_t\,\mathrm dt+\eta\,\mathrm dB_t,
    \qquad
    \mathrm dY_t=-\theta Y_t\,\mathrm dt+\sigma\,\mathrm dX_t,
    \label{eq:app-sim-lasso-dynamics}
\end{equation}
with \(X_0=x_0\) and \(Y_0=y_0\), and used the terminal target
\begin{equation}
    Y_T
    =
    y_0e^{-\theta T}
    +
    \sigma e^{-\theta T}
    \int_0^T e^{\theta s}\,\mathrm dX_s .
    \label{eq:app-sim-lasso-terminal}
\end{equation}
Since the integrand is deterministic and has zero quadratic covariation with \(X\), the It\^{o} and Stratonovich versions of \autoref{eq:app-sim-lasso-terminal} coincide.

We generated \(Y_T\) from \autoref{eq:app-sim-lasso-terminal} and observed \(Z=Y_T+\varepsilon\), with \(\varepsilon\sim N(0,\sigma_\varepsilon^2)\).
We fitted LASSO on scaled non-deterministic signature coordinates for \(K\in\{3,4,5,6\}\), removing pure-time coordinates because \(T\) is fixed.
The main-text figure uses \(T=1\), \(x_0=y_0=0\), \(\eta=\sigma=1\), \(\theta=6\), \(\kappa=2\), and \(\sigma_\varepsilon=1.2\), across the training sample sizes shown in \autoref{fig:sim-summary-lasso}.
Results were averaged over \(100\) independent replications.

Let \(\widehat A\) denote the selected support.
Prediction was evaluated by independent test MSE.
The selection metric in the main text is precision for the OU-relevant family,
\begin{equation}
    \operatorname{Precision}
    =
    \frac{|\widehat A\cap A_K^{\mathrm{OU}}|}{|\widehat A|},
    \label{eq:app-sim-lasso-precision}
\end{equation}
with the convention that precision is zero when \(\widehat A=\varnothing\).
For the selection diagnostic, the LASSO penalty was calibrated on a validation set by maximizing $\mathrm{F}_1$ relative to \(A_K^{\mathrm{OU}}\), where
\begin{equation}
    \operatorname{Recall}
    =
    \frac{|\widehat A\cap A_K^{\mathrm{OU}}|}{|A_K^{\mathrm{OU}}|},
    \qquad
    \operatorname{F}_1
    =
    \frac{2\operatorname{Precision}\operatorname{Recall}}
    {\operatorname{Precision}+\operatorname{Recall}} .
    \label{eq:app-sim-lasso-recall-f1}
\end{equation}
This oracle calibration was used only to study selection behavior in the controlled simulation.
Since \(Y_T\) was generated from the infinite integral target \autoref{eq:app-sim-lasso-terminal}, the support metric should be interpreted as recovery of the finite-level OU-relevant family \(A_K^{\mathrm{OU}}\), not as exact support recovery of a finite sparse linear model.

\subsection{Foreign Exchange Case Study}
\label{Appen_exchangerates}

This subsection provides additional details on the signature-based models estimated for the foreign exchange realized volatility application discussed in \autoref{subsec:fx_rv_application}.
The 12 modeled U.S. dollar pairs are USD/EUR, USD/GBP, USD/NZD, USD/CAD,
USD/CHF, USD/CZK, USD/DKK, USD/HUF, USD/MXN, USD/PLN, USD/SEK, and USD/SGD.

\subsubsection{Data construction and forecasting target.}
The HistData files contained one-minute bid and ask exchange rates over 2013--2023.
We used the mid-price and formed intraday log returns \(r_{t,j}\) for each trading day
\(t\). Daily realized volatility and the one-step-ahead response are
\begin{equation}
    \mathrm{RV}_t
    =
    \left(\sum_{j=1}^{M_t} r_{t,j}^2\right)^{1/2},
    \qquad
    Z_{t+1}=\log \mathrm{RV}_{t+1},
    \label{eq:fx-response}
\end{equation}
where \(M_t\) is the number of available one-minute returns on day \(t\). The log
transformation reduces the pronounced right skewness of realized volatility
\citep{andersen2001distribution,andersen2003modeling,BarndorffNielsenShephard2002JRSSB}.
For each currency pair, observations were ordered chronologically; the first 80\%
formed the training sample and the final 20\% formed the test sample. All feature
selection, model fitting, and tuning were performed on the training period, and the
final observations were used only for out-of-sample evaluation.

\subsubsection{Path representation and estimation.}
Within each day, we cumulated the one-minute returns and applied a time-augmented
lead--lag transformation. Denote the resulting three-dimensional path by
\(\widetilde{\bm X}_t\). The lead--lag coordinates encode successive return
increments and their quadratic variation, while the time coordinate preserves the
location of those movements within the trading day. We computed
\(S(\widetilde{\bm X}_t)^{\leq 3}\), which contains
\(3+3^2+3^3=39\) non-constant coordinates, and concatenated the signatures for
days \(t-21,\ldots,t\). The resulting 22-day predictor contained \(22\times39=858\)
coordinates. This window parallels the multi-horizon motivation of HAR models but
retains the ordered geometry within each day.

The training samples contained approximately 2,000 daily observations. With effective
path dimension three, the logarithmic truncation scale in \autoref{sec_Reg} is
\(\log(n)/\{2\log(3)\}\approx3.45\), which motivated the main choice \(K=3\).
We first screened the 858-coordinate dictionary with Signature-LASSO,
choosing its penalty by five-fold cross-validation within the training sample. We
then refit OLS on the selected coordinates, giving the reported Signature-OLS
forecast. The post-selection refit reduces shrinkage bias, but conventional refit $p$-values do not account for the preceding selection step and were therefore treated
as descriptive only.

The non-signature comparisons were a last-observation naive forecast, AR(1), and four
HAR-family specifications. HAR-RV uses daily, weekly, and monthly realized-volatility
components; HAR-RS replaces or augments these with realized semivariances; HAR-RV-J
adds a jump component; and HAR-CJ separates continuous and jump variation. These
benchmarks follow the constructions in
\citet{corsi2009simple,PattonSheppard2015REStat,AndersenBollerslevDiebold2007REStat,BarndorffNielsenShephard2004JFEC}.

\subsubsection{Evaluation and selected coordinates.}
We reported MAE and RMSE for the log-volatility response, out-of-sample
\(R^2_{\mathrm{oos}}\), and QLIKE on the variance scale \citep{Patton2011JoE}.
For observed variance
\(V_{t+1}=\mathrm{RV}_{t+1}^2\) and forecast \(\widehat V_{t+1}>0\), the QLIKE
contribution is \(\log \widehat V_{t+1}+V_{t+1}/\widehat V_{t+1}\), up to an
additive constant common to all forecasts. Lower MAE, RMSE, and QLIKE and higher
\(R^2_{\mathrm{oos}}\) indicate better performance. The main-text figure reports
pair-level RMSE comparisons, and \autoref{tab:fx_performance_summary} reports
cross-pair averages and win counts for all four metrics. The remainder of this
subsection examines which signature coordinates drive the post-selection forecasts.

\autoref{tab:summary_features} summarizes the signature features selected in at least \(40\%\) of the currency pairs.
Panel A groups the selected terms by signature channel, and Panel B groups them by the number of days ago from which the corresponding daily signature feature is taken.

\begin{table}[h!]
\centering
\caption{Summary of signature features selected in at least \(40\%\) of currency pairs.}
\label{tab:summary_features}
\small
\setlength{\tabcolsep}{18pt}
\renewcommand{\arraystretch}{0.6}
\begin{tabular}{lccc}
\toprule
\textbf{Channel/Lag} & \textbf{Total selections} & \textbf{Pairs covered} & \textbf{Pair frequency} \\
\midrule
\rowcolor{gray!20}
\multicolumn{4}{l}{\textbf{Panel A: By signature channel}}\\
\((\mathrm{lead},\mathrm{lag},t)\)   & 51 & 12 & 1.000 \\
\((\mathrm{lead},\mathrm{lag})\)     & 26 & 10 & 0.833 \\
\((t,\mathrm{lead},\mathrm{lag})\)   &  9 &  7 & 0.583 \\
\((\mathrm{lag},\mathrm{lead},t)\)   & 14 &  6 & 0.500 \\[1ex]

\rowcolor{gray!20}
\multicolumn{4}{l}{\textbf{Panel B: By lag}} \\
6    & 19 & 12 & 1.000 \\
1    & 18 & 12 & 1.000 \\
18   & 16 & 12 & 1.000 \\
12   & 15 & 11 & 0.917 \\
5    & 12 &  9 & 0.750 \\
11   &  8 &  5 & 0.417 \\
\bottomrule
\end{tabular}
\end{table}

The selected features show two clear patterns.
First, the selections are dominated by lead--lag signature terms, especially \((\mathrm{lead},\mathrm{lag},t)\), which appears in all 12 currency pairs.
Other frequently selected channels, such as \((\mathrm{lead},\mathrm{lag})\), \((t,\mathrm{lead},\mathrm{lag})\), and \((\mathrm{lag},\mathrm{lead},t)\), also encode ordered interactions between intraday price movements and time.
These terms go beyond scalar realized volatility summaries by retaining the temporal order of intraday fluctuations.

Second, the selected lags display a clear multi-horizon pattern.
Features from 1, 6, 12, and 18 trading days ago appear in nearly all currency pairs, with additional selections around one week, two weeks, and one month.
This partly echoes the daily, weekly, and monthly structure of HAR-type models.
However, HAR-type models summarize each horizon through realized volatility, semivariance, or jump-related scalar measures, whereas the selected signature features retain ordered lead--lag interactions within each intraday path.
Thus the feature-selection results suggest that the gains from Signature-OLS come not only from using multiple horizons, but also from capturing pathwise information that conventional HAR summaries miss.

\autoref{tab:fx_selected_features_long} reports the selected signature features and the post-selection OLS coefficients for each currency pair.
Within each currency pair, features are ordered by absolute coefficient size.
The reported p-values were computed from the OLS refit after feature selection and should be interpreted as descriptive evidence rather than formal post-selection inference.
\begingroup
\footnotesize
\setlength{\tabcolsep}{24pt}
\renewcommand{\arraystretch}{0.5}
\begin{longtable}{lccrr}
\caption{Selected signature features and post-selection OLS coefficients by currency pair.}
\label{tab:fx_selected_features_long}\\
\toprule
\textbf{Channel} & \textbf{Order} & \textbf{Days ago} & \textbf{Coef.} & \textbf{p-value} \\
\midrule
\endfirsthead

\multicolumn{5}{c}{\tablename~\thetable{} -- continued from previous page} \\
\toprule
\textbf{Channel} & \textbf{Order} & \textbf{Days ago} & \textbf{Coef.} & \textbf{p-value} \\
\midrule
\endhead

\midrule
\multicolumn{5}{r}{\emph{Continued on next page}} \\
\endfoot

\bottomrule
\endlastfoot

\rowcolor{gray!20}
\multicolumn{5}{l}{\textbf{\(\bullet\) Panel A: USD/CAD}} \\
\((\mathrm{lead},\mathrm{lag},t)\) & 3 & 6  & 0.3076  & \(7.59\times10^{-31}\) \\
\((\mathrm{lead},\mathrm{lag},t)\) & 3 & 12 & 0.2941  & \(8.97\times10^{-25}\) \\
\((\mathrm{lead},\mathrm{lag},t)\) & 3 & 18 & 0.2521  & \(1.14\times10^{-19}\) \\
\((\mathrm{lag},\mathrm{lead},t)\) & 3 & 12 & -0.1805 & \(3.49\times10^{-14}\) \\
\((\mathrm{lag},\mathrm{lead},t)\) & 3 & 18 & -0.1588 & \(7.11\times10^{-12}\) \\
\((\mathrm{lag},\mathrm{lead},t)\) & 3 & 6  & -0.1536 & \(1.15\times10^{-11}\) \\
\((t,\mathrm{lead},\mathrm{lag})\) & 3 & 1  & 0.0871  & \(2.38\times10^{-6}\) \\
\addlinespace

\rowcolor{gray!20}
\multicolumn{5}{l}{\textbf{\(\bullet\) Panel B: USD/CHF}} \\
\((\mathrm{lead},\mathrm{lag},t)\) & 3 & 6  & 0.3802  & \(4.68\times10^{-51}\) \\
\((\mathrm{lead},\mathrm{lag},t)\) & 3 & 18 & 0.3446  & \(8.43\times10^{-44}\) \\
\((\mathrm{lag},\mathrm{lead},t)\) & 3 & 6  & -0.2094 & \(1.80\times10^{-20}\) \\
\((\mathrm{lag},\mathrm{lead},t)\) & 3 & 18 & -0.2005 & \(9.59\times10^{-19}\) \\
\((t,\mathrm{lead},\mathrm{lag})\) & 3 & 1  & 0.1356  & \(1.15\times10^{-13}\) \\
\((\mathrm{lead},\mathrm{lag})\)   & 2 & 5  & 0.0868  & \(0.0637\) \\
\((\mathrm{lead},\mathrm{lag})\)   & 2 & 4  & 0.0472  & \(0.0103\) \\
\((\mathrm{lead},\mathrm{lag})\)   & 2 & 17 & 0.0327  & \(0.0741\) \\
\((\mathrm{lead},\mathrm{lag},t)\) & 3 & 5  & 0.0056  & \(0.9046\) \\
\addlinespace

\rowcolor{gray!20}
\multicolumn{5}{l}{\textbf{\(\bullet\) Panel C: USD/CZK}} \\
\((\mathrm{lead},\mathrm{lag})\)   & 2 & 6  & 0.1929  & \(6.31\times10^{-19}\) \\
\((t,\mathrm{lead},\mathrm{lag})\) & 3 & 1  & 0.1813  & \(2.62\times10^{-17}\) \\
\((\mathrm{lead},\mathrm{lag},t)\) & 3 & 12 & 0.1746  & \(1.53\times10^{-15}\) \\
\((\mathrm{lead},\mathrm{lag},t)\) & 3 & 18 & 0.1305  & \(2.53\times10^{-9}\) \\
\((t,\mathrm{lead},\mathrm{lag})\) & 3 & 5  & 0.1031  & \(1.62\times10^{-6}\) \\
\((\mathrm{lead},\mathrm{lag})\)   & 2 & 22 & 0.0850  & \(0.1244\) \\
\((t,\mathrm{lead},\mathrm{lag})\) & 3 & 11 & 0.0540  & \(0.2525\) \\
\((\mathrm{lead},\mathrm{lag})\)   & 2 & 11 & 0.0496  & \(0.2947\) \\
\((\mathrm{lead},\mathrm{lag},t)\) & 3 & 22 & 0.0058  & \(0.9166\) \\
\addlinespace

\rowcolor{gray!20}
\multicolumn{5}{l}{\textbf{\(\bullet\) Panel D: USD/DKK}} \\
\((\mathrm{lead},\mathrm{lag},t)\) & 3 & 6  & 0.3858  & \(8.99\times10^{-35}\) \\
\((\mathrm{lead},\mathrm{lag},t)\) & 3 & 12 & 0.3243  & \(8.11\times10^{-26}\) \\
\((\mathrm{lag},\mathrm{lead},t)\) & 3 & 6  & -0.2391 & \(1.09\times10^{-16}\) \\
\((\mathrm{lag},\mathrm{lead},t)\) & 3 & 12 & -0.2061 & \(1.44\times10^{-13}\) \\
\((\mathrm{lead},\mathrm{lag},t)\) & 3 & 18 & 0.1181  & \(7.73\times10^{-8}\) \\
\((t,\mathrm{lead},\mathrm{lag})\) & 3 & 1  & 0.0957  & \(0.0175\) \\
\((\mathrm{lead},\mathrm{lag})\)   & 2 & 1  & 0.0693  & \(0.0866\) \\
\((\mathrm{lead},\mathrm{lag})\)   & 2 & 11 & 0.0664  & \(0.1497\) \\
\((\mathrm{lead},\mathrm{lag})\)   & 2 & 5  & 0.0594  & \(0.2244\) \\
\((\mathrm{lead},\mathrm{lag})\)   & 2 & 4  & 0.0564  & \(0.0078\) \\
\((\mathrm{lead},\mathrm{lag})\)   & 2 & 22 & 0.0531  & \(0.0102\) \\
\((\mathrm{lead},\mathrm{lag},t)\) & 3 & 5  & 0.0454  & \(0.3502\) \\
\((\mathrm{lead},\mathrm{lag})\)   & 2 & 17 & 0.0249  & \(0.2399\) \\
\((\mathrm{lead},\mathrm{lag},t)\) & 3 & 11 & -0.0157 & \(0.7391\) \\
\addlinespace

\rowcolor{gray!20}
\multicolumn{5}{l}{\textbf{\(\bullet\) Panel E: USD/EUR}} \\
\((\mathrm{lead},\mathrm{lag},t)\) & 3 & 6  & 0.3533  & \(7.08\times10^{-29}\) \\
\((\mathrm{lead},\mathrm{lag},t)\) & 3 & 12 & 0.2942  & \(6.05\times10^{-21}\) \\
\((\mathrm{lead},\mathrm{lag},t)\) & 3 & 18 & 0.2369  & \(2.54\times10^{-15}\) \\
\((\mathrm{lag},\mathrm{lead},t)\) & 3 & 6  & -0.2100 & \(5.17\times10^{-13}\) \\
\((\mathrm{lag},\mathrm{lead},t)\) & 3 & 12 & -0.1828 & \(8.44\times10^{-11}\) \\
\((\mathrm{lag},\mathrm{lead},t)\) & 3 & 18 & -0.1572 & \(8.25\times10^{-9}\) \\
\((t,\mathrm{lead},\mathrm{lag})\) & 3 & 1  & 0.1095  & \(0.0053\) \\
\((\mathrm{lead},\mathrm{lag})\)   & 2 & 11 & 0.0721  & \(0.1100\) \\
\((\mathrm{lead},\mathrm{lag})\)   & 2 & 5  & 0.0611  & \(0.2002\) \\
\((\mathrm{lead},\mathrm{lag})\)   & 2 & 4  & 0.0543  & \(0.0093\) \\
\((\mathrm{lead},\mathrm{lag})\)   & 2 & 1  & 0.0540  & \(0.1713\) \\
\((\mathrm{lead},\mathrm{lag},t)\) & 3 & 5  & 0.0467  & \(0.3243\) \\
\((\mathrm{lead},\mathrm{lag})\)   & 2 & 22 & 0.0396  & \(0.0533\) \\
\((\mathrm{lead},\mathrm{lag},t)\) & 3 & 11 & -0.0224 & \(0.6269\) \\
\((\mathrm{lead},\mathrm{lag})\)   & 2 & 17 & 0.0209  & \(0.3163\) \\
\addlinespace

\rowcolor{gray!20}
\multicolumn{5}{l}{\textbf{\(\bullet\) Panel F: USD/GBP}} \\
\((\mathrm{lead},\mathrm{lag},t)\) & 3 & 6  & 0.2411  & \(2.28\times10^{-7}\) \\
\((\mathrm{lead},\mathrm{lag},t)\) & 3 & 1  & 0.1824  & \(1.37\times10^{-16}\) \\
\((\mathrm{lead},\mathrm{lag},t)\) & 3 & 18 & 0.1810  & \(5.65\times10^{-16}\) \\
\((\mathrm{lead},\mathrm{lag},t)\) & 3 & 12 & 0.1787  & \(3.38\times10^{-15}\) \\
\((\mathrm{lag},\mathrm{lag},t)\)  & 3 & 6  & 0.0361  & \(0.4229\) \\
\addlinespace

\rowcolor{gray!20}
\multicolumn{5}{l}{\textbf{\(\bullet\) Panel G: USD/HUF}} \\
\((\mathrm{lead},\mathrm{lag},t)\) & 3 & 6  & 0.3934  & \(3.99\times10^{-44}\) \\
\((\mathrm{lead},\mathrm{lag},t)\) & 3 & 12 & 0.2490  & \(2.30\times10^{-17}\) \\
\((\mathrm{lag},\mathrm{lead},t)\) & 3 & 6  & -0.2210 & \(7.74\times10^{-19}\) \\
\((\mathrm{lead},\mathrm{lag},t)\) & 3 & 18 & 0.2010  & \(6.10\times10^{-13}\) \\
\((\mathrm{lag},\mathrm{lead},t)\) & 3 & 12 & -0.1457 & \(1.51\times10^{-8}\) \\
\((\mathrm{lag},\mathrm{lead},t)\) & 3 & 18 & -0.1401 & \(1.37\times10^{-8}\) \\
\((\mathrm{lead},\mathrm{lag})\)   & 2 & 1  & 0.1054  & \(7.62\times10^{-8}\) \\
\((\mathrm{lead},\mathrm{lag},t)\) & 3 & 5  & 0.0491  & \(0.0142\) \\
\addlinespace

\rowcolor{gray!20}
\multicolumn{5}{l}{\textbf{\(\bullet\) Panel H: USD/MXN}} \\
\((\mathrm{lead},\mathrm{lag},t)\) & 3 & 6  & 0.6665  & \(2.00\times10^{-48}\) \\
\((\mathrm{lag},\mathrm{lead},t)\) & 3 & 6  & -0.3008 & \(4.15\times10^{-22}\) \\
\((\mathrm{lag},\mathrm{lag},\mathrm{lead})\) & 3 & 7  & -0.2678 & \(6.26\times10^{-17}\) \\
\((\mathrm{lag})\)                 & 1 & 1  & 0.1457  & \(2.07\times10^{-8}\) \\
\((\mathrm{lead},\mathrm{lag},t)\) & 3 & 12 & 0.0925  & \(0.0013\) \\
\((\mathrm{lead},\mathrm{lag},t)\) & 3 & 18 & 0.0882  & \(0.0017\) \\
\((\mathrm{lead},\mathrm{lag})\)   & 2 & 5  & 0.0728  & \(0.0036\) \\
\((\mathrm{lead},\mathrm{lag},t)\) & 3 & 2  & 0.0652  & \(0.0278\) \\
\((\mathrm{lead},\mathrm{lag})\)   & 2 & 1  & 0.0550  & \(0.1066\) \\
\((\mathrm{lead},\mathrm{lag},t)\) & 3 & 1  & 0.0430  & \(0.2638\) \\
\addlinespace

\rowcolor{gray!20}
\multicolumn{5}{l}{\textbf{\(\bullet\) Panel I: USD/NZD}} \\
\((\mathrm{lead},\mathrm{lag},t)\) & 3 & 6  & 0.2455  & \(8.88\times10^{-42}\) \\
\((\mathrm{lead},\mathrm{lag})\)   & 2 & 1  & 0.1203  & \(0.0005\) \\
\((\mathrm{lead},\mathrm{lag},t)\) & 3 & 18 & 0.1130  & \(4.16\times10^{-11}\) \\
\((\mathrm{lead},\mathrm{lag},t)\) & 3 & 5  & 0.1099  & \(3.85\times10^{-10}\) \\
\((t,\mathrm{lead},\mathrm{lag})\) & 3 & 1  & 0.0975  & \(0.0041\) \\
\((\mathrm{lead},\mathrm{lag},t)\) & 3 & 12 & 0.0811  & \(2.82\times10^{-5}\) \\
\((\mathrm{lead},\mathrm{lag},t)\) & 3 & 11 & 0.0209  & \(0.2677\) \\
\addlinespace

\rowcolor{gray!20}
\multicolumn{5}{l}{\textbf{\(\bullet\) Panel J: USD/PLN}} \\
\((\mathrm{lead},\mathrm{lag},t)\) & 3 & 6  & 0.2626  & \(1.55\times10^{-32}\) \\
\((\mathrm{lead},\mathrm{lag},t)\) & 3 & 18 & 0.1824  & \(3.43\times10^{-17}\) \\
\((\mathrm{lead},\mathrm{lag},t)\) & 3 & 12 & 0.1776  & \(7.98\times10^{-16}\) \\
\((\mathrm{lead},\mathrm{lag})\)   & 2 & 1  & 0.0915  & \(0.0013\) \\
\((\mathrm{lead},\mathrm{lag},t)\) & 3 & 1  & 0.0895  & \(0.0017\) \\
\((\mathrm{lead},\mathrm{lag},t)\) & 3 & 5  & 0.0879  & \(1.79\times10^{-5}\) \\
\addlinespace

\rowcolor{gray!20}
\multicolumn{5}{l}{\textbf{\(\bullet\) Panel K: USD/SEK}} \\
\((\mathrm{lead},\mathrm{lag},t)\) & 3 & 6  & 0.2572  & \(1.07\times10^{-36}\) \\
\((\mathrm{lead},\mathrm{lag},t)\) & 3 & 12 & 0.1916  & \(1.99\times10^{-21}\) \\
\((\mathrm{lead},\mathrm{lag},t)\) & 3 & 18 & 0.1436  & \(7.29\times10^{-13}\) \\
\((\mathrm{lead},\mathrm{lag})\)   & 2 & 1  & 0.1240  & \(6.20\times10^{-11}\) \\
\addlinespace

\rowcolor{gray!20}
\multicolumn{5}{l}{\textbf{\(\bullet\) Panel L: USD/SGD}} \\
\((\mathrm{lead},\mathrm{lag},t)\) & 3 & 6  & 0.1928  & \(1.93\times10^{-28}\) \\
\((\mathrm{lead},\mathrm{lag})\)   & 2 & 1  & 0.1554  & \(2.05\times10^{-21}\) \\
\((\mathrm{lead},\mathrm{lag},t)\) & 3 & 12 & 0.1318  & \(8.01\times10^{-13}\) \\
\((\mathrm{lead},\mathrm{lag},t)\) & 3 & 18 & 0.1112  & \(2.87\times10^{-10}\) \\
\((\mathrm{lead},\mathrm{lag},t)\) & 3 & 5  & 0.0996  & \(7.71\times10^{-9}\) \\
\((\mathrm{lead},\mathrm{lag})\)   & 2 & 4  & 0.0927  & \(6.06\times10^{-9}\) \\
\((t,\mathrm{lead},\mathrm{lag})\) & 3 & 13 & 0.0726  & \(1.38\times10^{-5}\) \\
\((\mathrm{lead},\mathrm{lag},t)\) & 3 & 17 & 0.0692  & \(6.47\times10^{-5}\) \\
\((\mathrm{lead},\mathrm{lag},t)\) & 3 & 11 & 0.0549  & \(0.0015\) \\

\end{longtable}
\endgroup

\subsection{Battery End-of-Life Case Study}
\label{Appen_battery}

This subsection provides supplementary details for the battery EOL case study in
\autoref{sec:battery}. The retained sample contained 241 cells, with 166 cells in the
training set and 75 cells in the test set. \autoref{tab:battery_split_summary} summarizes
the empirical distribution of \(Z_j=\mathrm{EOL}_j\) in the two splits. All reported
results used the same training--test split as in the main text.

For cell \(j\), \(\mathrm{EOL}_j\) is the first cycle at which discharge capacity
falls to 80\% of its initial value. We used the first usable HPPC record after the
preprocessing and chemistry filters used to form the modeled sample and retained
the same early pulse window for signature and handcrafted features. The retained
    path was \(\bm X_j=\{\bm X_{t,j}:t\in[0,T]\}\), where
\begin{equation}
    \bm X_{t,j}
    =
    \left(
    t,\mathrm{Cur}_{t,j},\mathrm{Vol}_{t,j},
    \dot{\mathrm{Cur}}_{t,j},\dot{\mathrm{Vol}}_{t,j}
    \right)^\top .
    \label{eq:app-battery-path}
\end{equation}
The derivative channels record local current and voltage changes and hence retain
transient behavior that is lost in pulse-level averages.
The main specification used \(K=3\) and \(T=120\) seconds. The resulting signature feature
vector \(\bm s_K(\bm X_j)\) contained \(5+5^2+5^3=155\) non-constant coordinates. The
handcrafted baseline followed the feature-engineering philosophy of \citet{Paulson2022},
but was restricted to quantities available from the retained HPPC pulse. It consisted of 21
pulse-level summaries computed over the same \(120\)-second window, including the number of
retained observations, pulse duration, absolute and signed current moments, voltage
moments, voltage range, terminal voltage change, current sign-change counts, positive and
negative current fractions, absolute current throughput, energy-like throughput, the number
of nonzero-current segments, and segment-level summaries of duration, current magnitude,
and voltage change.

We fitted OLS, ridge, and LASSO to each representation and evaluated MAE and RMSE in
cycles together with \(R^2\), reporting both training and test values. The common split,
common response, and common pulse window made each signature--handcrafted comparison
a feature-representation comparison conditional on the learner. The training columns
diagnose interpolation and shrinkage, whereas the held-out test columns assess
generalization.

\begin{table}[h!]
\centering
\caption{Summary of the retained battery end-of-life (EOL) sample.}
\label{tab:battery_split_summary}
\small
\renewcommand{\arraystretch}{0.65}
\begin{tabular}{lrrrrrr}
\toprule
\multirow{2}{*}{\textbf{Split}}
& \multirow{2}{*}{\(\bm n\)}
& \multicolumn{5}{c}{\textbf{End-of-Life Summary}} \\
\cmidrule(lr){3-7}
& & \textbf{Mean} & \textbf{Std. Dev.} & \textbf{Minimum} & \textbf{Median} & \textbf{Maximum} \\
\midrule
Training & 166 & 581.3 & 579.1 & 21 & 441 & 2313 \\
Test     & 75  & 510.0 & 487.1 & 17 & 384 & 2313 \\
\bottomrule
\end{tabular}
\end{table}

\autoref{tab:depth_sensitivity} reports sensitivity to the signature truncation depth
\(K\). For LASSO, increasing \(K\) from 1 to 3 reduces the test MAE from 362.8 to 188.9
cycles. Beyond \(K=3\), the error increases as the feature dimension grows, suggesting a
finite-sample trade-off between richer path representation and estimation stability. Ridge
is less sensitive to sparsity but does not attain the same out-of-sample accuracy as LASSO
in the main specification.

\begin{table}[h!]
\centering
\caption{Sensitivity to signature depth \(K\). MAE and RMSE denote test mean absolute
error and root mean squared error, both measured in cycles; \(R^2\) denotes the test-set
coefficient of determination.}
\label{tab:depth_sensitivity}
\small
\setlength{\tabcolsep}{7pt}
\renewcommand{\arraystretch}{0.65}
\begin{tabular}{rrccc ccc}
\toprule
\multirow{2}{*}{\(\bm K\)}
& \multirow{2}{*}{\textbf{Feature Count}}
& \multicolumn{3}{c}{\textbf{LASSO}}
& \multicolumn{3}{c}{\textbf{Ridge}} \\
\cmidrule(lr){3-5}\cmidrule(lr){6-8}
& & \textbf{MAE} & \textbf{RMSE} & \(\bm{R^2}\)
& \textbf{MAE} & \textbf{RMSE} & \(\bm{R^2}\) \\
\midrule
1 & 5   & 362.8 & 520.1 & -0.156 & 361.1 & 518.3 & -0.148 \\
2 & 30   & 234.0 & 397.3 & 0.326  & 243.1 & 467.5 & 0.066  \\
3 & 155  & \textbf{188.9} & \textbf{331.9} & \textbf{0.529} & 314.1 & 463.1 & 0.084 \\
4 & 780  & 272.3 & 435.6 & 0.189  & 277.8 & \textbf{421.4} & \textbf{0.241} \\
5 & 3905 & 253.5 & 428.0 & 0.217  & 291.1 & 433.6 & 0.197 \\
\bottomrule
\end{tabular}
\end{table}

\autoref{tab:threshold_sensitivity} reports sensitivity to the retained pulse horizon
\(T\). The early part of the HPPC pulse is highly informative: for LASSO, the lowest test
MAE in this grid is achieved at \(T=80\) seconds. Ridge attains its best performance at a
longer horizon, \(T=200\) seconds. These results are reported as robustness diagnostics.
The main specification fixes \(T=120\) seconds to match the benchmark setting used in
\autoref{sec:battery}.

\begin{table}[h!]
\centering
\caption{Sensitivity to the retained pulse horizon \(T\). MAE and RMSE denote test mean
absolute error and root mean squared error, both measured in cycles; \(R^2\) denotes the
test-set coefficient of determination.}
\label{tab:threshold_sensitivity}
\small
\setlength{\tabcolsep}{12pt}
\renewcommand{\arraystretch}{0.65}
\begin{tabular}{rccc ccc}
\toprule
\multirow{2}{*}{\(\bm T\) \textbf{(s)}}
& \multicolumn{3}{c}{\textbf{LASSO}}
& \multicolumn{3}{c}{\textbf{Ridge}} \\
\cmidrule(lr){2-4}\cmidrule(lr){5-7}
& \textbf{MAE} & \textbf{RMSE} & \(\bm{R^2}\)
& \textbf{MAE} & \textbf{RMSE} & \(\bm{R^2}\) \\
\midrule
40  & 248.8 & 410.1 & 0.281 & 265.6 & 411.8 & 0.275 \\
80  & \textbf{179.3} & \textbf{276.9} & \textbf{0.672} & 228.7 & 433.7 & 0.196 \\
120 & 188.9 & 331.9 & 0.529 & 314.1 & 463.1 & 0.084 \\
160 & 287.2 & 441.9 & 0.166 & 297.8 & 448.0 & 0.143 \\
200 & 320.0 & 463.1 & 0.084 & \textbf{185.0} & \textbf{292.5} & \textbf{0.634} \\
\bottomrule
\end{tabular}
\end{table}

\begin{figure}[h!]
\centering
\includegraphics[width=\linewidth]{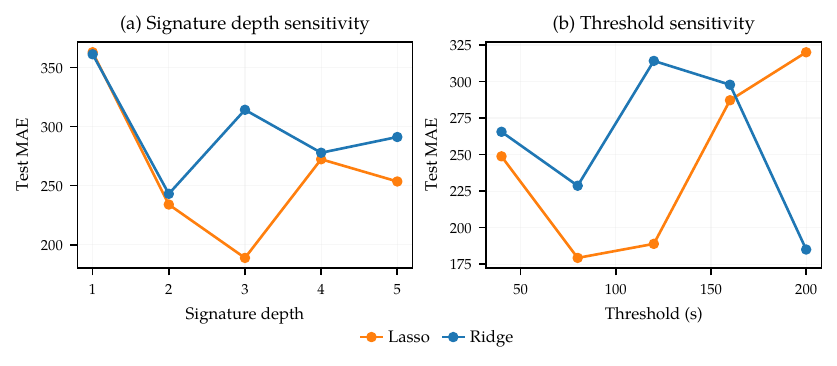}
\caption{Sensitivity of battery end-of-life (EOL) prediction accuracy to signature depth and retained pulse
horizon. Panel (a) plots test mean absolute error (MAE) against \(K\) for ridge and LASSO. Panel (b) plots test
MAE against \(T\) for ridge and LASSO.}
\label{fig:battery_sensitivity}
\end{figure}

\autoref{fig:battery_sensitivity} visualizes the two sensitivity analyses in
\autoref{tab:depth_sensitivity} and \autoref{tab:threshold_sensitivity}. The left panel
shows that the LASSO curve is minimized at \(K=3\), while the right panel shows that the
best retained horizon differs across regularization schemes. These patterns support the use
of sparse regularization when working with high-dimensional signature features.

\begin{figure}[H]
\centering
\includegraphics[width=\linewidth]{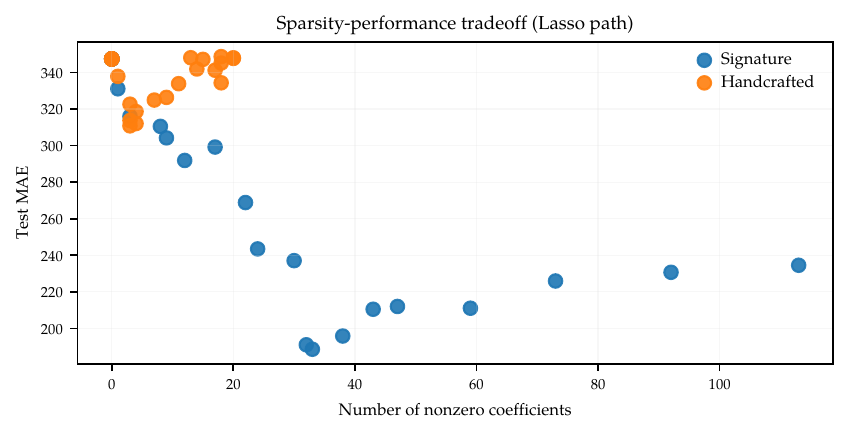}
\caption{Sparsity--performance trade-off for signature and handcrafted representations.
The horizontal axis reports the number of nonzero coefficients and the vertical axis
reports test mean absolute error (MAE).}
\label{fig:battery_lasso_path}
\end{figure}

Finally, \autoref{fig:battery_lasso_path} reports the sparsity--performance trade-off along
the LASSO regularization path. The best test accuracy is attained in a moderately sparse
region rather than in the densest model. This pattern is consistent with the sparse
signature-regression interpretation: although \(\bm s_K(\bm X_j)\) is high-dimensional,
only a subset of its coordinates appears to carry substantial predictive information for
\(Z_j=\mathrm{EOL}_j\).

\subsection{Seizure Detection Case Study}
\label{app:seizure_appendix}

This subsection provides supplementary details for the seizure detection case study in
\autoref{sec:seizure}. All appendix results were restricted to the \(22\)-subject modeled
CHB-MIT cohort used in the main text. After correcting the cohort to retain only subjects
with complete modeling outputs, the final sample contained \(620\) EDF files, \(803.2\)
recording hours, \(9{,}758\) seizure windows, and \(2{,}880{,}127\) non-seizure windows.
Thus the cohort-level seizure prevalence was only \(0.338\%\), which highlights the extreme
class imbalance of this application. All appendix results used the same patient-specific
stratified \(25\%/75\%\) train-test split as in the main text.

Each retained recording was divided into non-overlapping three-second windows. For patient
\(j\) and window \(k\), the response \(Z_{j,k}\) equaled one if the window overlapped an
annotated seizure interval and zero otherwise. Splitting and evaluation were conducted
within patient. This design avoids allowing the large number of windows from a few
long-recorded subjects to determine the cohort result and preserves patient-level
heterogeneity in the reported standard deviations.

\begin{table}[H]
\centering
\small
\caption{Modeled-cohort dataset characteristics after clean windowing. EDF denotes
European Data Format. The last column
reports the percentage of seizure windows among all retained windows for each subject.}
\label{tab:seizure_dataset}
\setlength{\tabcolsep}{13pt}
\renewcommand{\arraystretch}{0.5}
 \resizebox{\linewidth}{!}{
\begin{tabular}{lrrrrr}
\toprule
\multirow{2}{*}{\textbf{Subject}}
& \multicolumn{2}{c}{\textbf{Recording}}
& \multicolumn{2}{c}{\textbf{Window Count}}
& \multicolumn{1}{c}{\textbf{Prevalence}} \\
\cmidrule(lr){2-3}\cmidrule(lr){4-5}\cmidrule(lr){6-6}
& \textbf{EDF Files} & \textbf{Duration (h)}
& \textbf{Seizure Win.} & \textbf{Non-seizure Win.}
& \textbf{Ratio (\%)} \\
\midrule
chb01 & 42 & 40.6 & 442 & 145448 & 0.3 \\
chb02 & 36 & 35.3 & 172 & 126709 & 0.1 \\
chb03 & 38 & 38.0 & 402 & 136314 & 0.3 \\
chb05 & 39 & 39.0 & 558 & 139764 & 0.4 \\
chb06 & 18 & 66.7 & 153 & 240037 & 0.1 \\
chb07 & 19 & 67.1 & 325 & 241019 & 0.1 \\
chb08 & 20 & 20.0 & 919 & 71054 & 1.3 \\
chb09 & 19 & 67.9 & 276 & 244016 & 0.1 \\
chb10 & 25 & 50.0 & 447 & 179573 & 0.2 \\
chb11 & 35 & 34.8 & 806 & 124375 & 0.6 \\
chb13 & 33 & 33.0 & 535 & 118175 & 0.5 \\
chb14 & 26 & 26.0 & 169 & 93363 & 0.2 \\
chb15 & 40 & 40.0 & 1992 & 141924 & 1.4 \\
chb16 & 19 & 19.0 & 84 & 68258 & 0.1 \\
chb17 & 21 & 21.0 & 293 & 75283 & 0.4 \\
chb18 & 36 & 35.6 & 317 & 127884 & 0.2 \\
chb19 & 30 & 29.9 & 236 & 107444 & 0.2 \\
chb20 & 29 & 27.6 & 294 & 98998 & 0.3 \\
chb21 & 33 & 32.8 & 199 & 117916 & 0.2 \\
chb22 & 31 & 31.0 & 204 & 111339 & 0.2 \\
chb23 & 9 & 26.6 & 424 & 95154 & 0.4 \\
chb24 & 22 & 21.3 & 511 & 76080 & 0.7 \\
\midrule
\textbf{Total} & \textbf{620} & \textbf{803.2} & \textbf{9758} & \textbf{2880127} & \textbf{0.3} \\
\bottomrule
\end{tabular}
}
\end{table}

\autoref{tab:seizure_dataset} shows substantial heterogeneity across subjects in both
recording length and seizure prevalence. Some subjects, such as chb15, contribute many
more seizure windows than others, whereas for several subjects the seizure proportion is
close to \(0.1\%\). This heterogeneity motivates evaluating performance at the patient
level and then averaging across subjects, rather than pooling all windows and reporting
only a cohort-level score.

For the \(k\)-th retained window of patient \(j\), let
\(\bm x_{t,j,k}\) be the filtered multichannel EEG vector. The path supplied to the
signature transform was
\(\bm X_{j,k}=\{\bm X_{t,j,k}:t\in[0,T]\}\), with \(T=3\) seconds and
\(\bm X_{t,j,k}=(t,\bm x_{t,j,k},\mathrm d\bm x_{t,j,k}/\mathrm dt)^\top\). Channel
availability was harmonized within each patient by taking the subject-specific channel
intersection, so the resulting channel count ranged from \(18\) to \(31\), with median
\(23\). The handcrafted baseline contained line length, nonlinear energy, and Shannon
entropy for each retained channel, which yielded between \(54\) and \(93\) features
depending on the patient. The main signature specification used the depth-\(2\)
truncated signature \(\bm{s}_2(\bm X_{j,k})=S(\bm X_{j,k})^{\le2}_T\), which yielded between
\(1{,}406\) and \(4{,}032\) features. To handle this high-dimensional representation, we
fitted elastic-net-regularized logistic regression with balanced class weights. This
specification is labeled ``Signature + LR'' in the main text and the tables below. The handcrafted baselines
were logistic regression and random forest. In the main text, the operating threshold for
the signature classifier was fixed at \(0.10\), the maximizer of patient-averaged balanced
accuracy in the reported threshold sweep.

For a fixed threshold, let TP, TN, FP, and FN denote the patient-level counts. The primary
threshold-dependent rates are
\begin{equation}
\begin{aligned}
\mathrm{SENS}&=\frac{\mathrm{TP}}{\mathrm{TP}+\mathrm{FN}},
&\qquad
\mathrm{SPEC}&=\frac{\mathrm{TN}}{\mathrm{TN}+\mathrm{FP}},\\
\mathrm{B\text{-}ACC}&=\frac{\mathrm{SENS}+\mathrm{SPEC}}{2},
&
\mathrm{F1}&=\frac{2\mathrm{TP}}{2\mathrm{TP}+\mathrm{FP}+\mathrm{FN}}.
\end{aligned}
\label{eq:seizure-metrics}
\end{equation}
Balanced accuracy gives equal weight to the two classes, whereas F1 depends only on
positive-class retrieval and precision. AUC and average precision (AP) summarize ranking
performance over thresholds; AP is especially informative when prevalence is very low.
All reported metrics were computed for each subject and then averaged, so they are not
equivalent to metrics computed from one pooled cohort-level confusion matrix. The threshold
sweep is presented as an operating-point diagnostic; prospective threshold tuning should be
nested within training data.

\begin{table}[H]
\centering
\small
\caption{Threshold sensitivity across methods. ACC, B-ACC, SENS, SPEC, F1, MCC, AUC,
and AP denote accuracy, balanced accuracy, sensitivity, specificity, F1 score, Matthews
correlation coefficient, area under the receiver operating characteristic curve, and average
precision, respectively. Entries are means across subjects; all metrics except MCC are
reported as percentages.}
\label{tab:seizure_threshold_sensitivity}
\setlength{\tabcolsep}{7.25pt}
\renewcommand{\arraystretch}{0.5}
 \resizebox{\linewidth}{!}{
\begin{tabular}{lcccccccc}
\toprule
\multirow{2}{*}{\textbf{Method / Threshold}}
& \multicolumn{8}{c}{\textbf{Performance Metric}} \\
\cmidrule(lr){2-9}
& \textbf{ACC} & \textbf{B-ACC} & \textbf{SENS} & \textbf{SPEC}
& \textbf{F1} & \textbf{MCC} & \textbf{AUC} & \textbf{AP} \\
\midrule
Handcrafted + LR @ 0.01 & 77.23 & 85.81 & 94.52 & 77.10 & 4.65 & 0.124 & 95.93 & 39.31 \\
Handcrafted + LR @ 0.10 & 88.16 & 89.43 & 90.77 & 88.09 & 10.15 & 0.196 & 95.93 & 39.31 \\
Handcrafted + LR @ 0.20 & 91.51 & 90.28 & 89.09 & 91.48 & 13.26 & 0.230 & 95.93 & 39.31 \\
Handcrafted + LR @ 0.30 & 93.50 & 90.36 & 87.23 & 93.49 & 15.78 & 0.255 & 95.93 & 39.31 \\
Handcrafted + LR @ 0.40 & 94.91 & 90.22 & 85.53 & 94.92 & 18.10 & 0.276 & 95.93 & 39.31 \\
Handcrafted + LR @ 0.50 & 96.00 & 89.92 & 83.82 & 96.02 & 20.37 & 0.296 & 95.93 & 39.31 \\
\midrule
Handcrafted + RF @ 0.01 & 94.94 & 92.75 & 90.57 & 94.92 & 15.11 & 0.263 & 96.81 & 62.70 \\
Handcrafted + RF @ 0.10 & 99.66 & 82.10 & 64.42 & 99.78 & 59.08 & 0.594 & 96.81 & 62.70 \\
Handcrafted + RF @ 0.20 & 99.75 & 74.35 & 48.76 & 99.94 & 56.97 & 0.593 & 96.81 & 62.70 \\
Handcrafted + RF @ 0.30 & 99.74 & 68.96 & 37.95 & 99.98 & 49.54 & 0.541 & 96.81 & 62.70 \\
Handcrafted + RF @ 0.40 & 99.72 & 64.55 & 29.12 & 99.99 & 40.64 & 0.475 & 96.81 & 62.70 \\
Handcrafted + RF @ 0.50 & 99.70 & 60.90 & 21.81 & 100.00 & 32.03 & 0.403 & 96.81 & 62.70 \\
\midrule
Signature + LR @ 0.01 & 79.78 & 88.72 & 97.74 & 79.69 & 3.70 & 0.117 & 98.64 & 75.79 \\
Signature + LR @ 0.10 & 97.49 & 95.78 & 94.08 & 97.49 & 27.28 & 0.377 & 98.64 & 75.79 \\
Signature + LR @ 0.20 & 98.74 & 95.47 & 92.19 & 98.76 & 43.11 & 0.505 & 98.64 & 75.79 \\
Signature + LR @ 0.30 & 99.18 & 94.94 & 90.69 & 99.20 & 52.82 & 0.581 & 98.64 & 75.79 \\
Signature + LR @ 0.40 & 99.41 & 94.24 & 89.04 & 99.44 & 59.56 & 0.633 & 98.64 & 75.79 \\
Signature + LR @ 0.50 & 99.55 & 93.66 & 87.72 & 99.59 & 65.00 & 0.676 & 98.64 & 75.79 \\
\bottomrule
\end{tabular}
}
\end{table}

\autoref{tab:seizure_threshold_sensitivity} reports the full threshold sweep. Several
patterns are worth noting. First, the signature classifier achieves its highest
patient-averaged balanced accuracy at threshold \(0.10\), which justifies the operating
point used in the main text. Its F1 score, by contrast, continues to increase at higher
thresholds and peaks at \(0.50\). This is not contradictory: balanced accuracy gives equal
weight to sensitivity and specificity, whereas F1 places greater emphasis on positive
predictions and is therefore more sensitive to threshold-dependent calibration.

Second, the handcrafted random forest behaves differently. At threshold \(0.10\) it
achieves very high accuracy and specificity, but its sensitivity drops sharply. As the
threshold increases further, sensitivity deteriorates even more, and balanced accuracy
falls accordingly. This confirms that the random forest is calibrated toward non-seizure
detection, which is unsurprising in view of the severe class imbalance. By contrast, the
signature classifier maintains both high sensitivity and high specificity over a wider
threshold range.

Third, the AUC and AP columns show that the signature classifier also dominates in
threshold-free ranking performance. In particular, its average precision is substantially
higher than that of both handcrafted baselines, which is especially relevant in an
ultra-imbalanced setting where precision-recall performance is often more informative than
raw accuracy.

\begin{table}[H]
\centering
\small
\caption{Subject-wise comparison between Handcrafted + RF and Signature + LR at
threshold \(0.10\), where RF and LR denote random forest and logistic regression.
Balanced accuracy, F1 score, and average precision are reported as percentages.}
\label{tab:seizure_subjectwise}
\setlength{\tabcolsep}{18.5pt}
\renewcommand{\arraystretch}{0.5}
 \resizebox{\linewidth}{!}{
\begin{tabular}{lccccc}
\toprule
\multirow{2}{*}{\textbf{Subject}}
& \multicolumn{2}{c}{\textbf{Balanced Accuracy}}
& \multicolumn{2}{c}{\textbf{F1 Score}}
& \multicolumn{1}{c}{\textbf{Average Precision}} \\
\cmidrule(lr){2-3}\cmidrule(lr){4-5}\cmidrule(lr){6-6}
& \textbf{RF} & \textbf{Signature}
& \textbf{RF} & \textbf{Signature}
& \textbf{Signature} \\
\midrule
chb01 & 91.16 & 98.48 & 66.34 & 21.95 & 93.02 \\
chb02 & 80.58 & 98.30 & 56.63 & 46.73 & 95.06 \\
chb03 & 84.05 & 98.62 & 66.99 & 34.43 & 87.73 \\
chb05 & 94.19 & 99.35 & 77.29 & 57.62 & 97.49 \\
chb06 & 87.38 & 93.09 & 74.46 & 12.37 & 50.28 \\
chb07 & 78.87 & 99.10 & 61.04 & 21.27 & 88.96 \\
chb08 & 93.75 & 94.00 & 68.81 & 21.55 & 62.80 \\
chb09 & 85.74 & 96.53 & 76.88 & 24.39 & 81.80 \\
chb10 & 95.63 & 99.19 & 81.71 & 53.60 & 94.67 \\
chb11 & 95.74 & 98.65 & 78.89 & 33.87 & 98.12 \\
chb13 & 83.06 & 92.32 & 49.86 & 12.57 & 47.02 \\
chb14 & 81.85 & 93.66 & 60.90 & 26.63 & 79.79 \\
chb15 & 72.27 & 95.34 & 43.36 & 25.73 & 75.10 \\
chb16 & 55.54 & 82.45 & 15.73 & 8.32 & 25.70 \\
chb17 & 78.53 & 96.17 & 54.43 & 19.18 & 78.11 \\
chb18 & 70.32 & 94.98 & 44.29 & 12.20 & 57.76 \\
chb19 & 76.46 & 94.07 & 44.98 & 40.20 & 70.84 \\
chb20 & 80.04 & 98.28 & 62.44 & 18.84 & 88.35 \\
chb21 & 64.33 & 98.47 & 23.63 & 16.18 & 81.84 \\
chb22 & 84.92 & 98.17 & 63.69 & 46.76 & 92.96 \\
chb23 & 89.00 & 98.11 & 63.20 & 31.53 & 72.57 \\
chb24 & 82.77 & 89.88 & 64.29 & 14.25 & 47.49 \\
\bottomrule
\end{tabular}
}
\end{table}

\autoref{tab:seizure_subjectwise} reports subject-level heterogeneity. The most stable
pattern is the balanced-accuracy comparison: the signature classifier improves on the
handcrafted random forest for every subject. The size of the improvement varies, but in
many subjects it is substantial. At the same time, F1 scores can be lower for the
signature model at threshold \(0.10\), reflecting the tension between balanced detection
and positive predictive concentration in an imbalanced classification problem. This again
shows that no single metric fully summarizes performance, and that balanced accuracy is a
more appropriate primary criterion for the operating goal considered here.



\end{document}